\documentclass[10pt]{article}
\usepackage[latin1]{inputenc}
\usepackage{epsfig}
\usepackage{color}
\usepackage[british,english]{babel}
\usepackage{amsthm}
\usepackage{amsmath}
\usepackage{amsfonts}
\usepackage{amssymb}
\usepackage{graphicx}

\usepackage{fancyhdr}

\newtheorem{corollary}{Corollary}[section]
\newtheorem{definition}[corollary]{Definition}

\newtheorem{lemma}[corollary]{Lemma}

\newtheorem{remark}[corollary]{Remark}
\newtheorem{theorem}[corollary]{Theorem}

\date{}
\begin{document}
\title{Modeling non-Newtonian fluids in a thin domain perforated \\
 with cylinders of small diameter}\maketitle

\vskip-30pt
 \centerline{Mar\'ia ANGUIANO\footnote{Departamento de An\'alisis Matem\'atico. Facultad de Matem\'aticas.
Universidad de Sevilla, 41012 Sevilla (Spain)
anguiano@us.es}, Francisco J. SU\'AREZ-GRAU\footnote{Departamento de Ecuaciones Diferenciales y An\'alisis Num\'erico. Facultad de Matem\'aticas. Universidad de Sevilla, 41012 Sevilla (Spain) fjsgrau@us.es}}

 \renewcommand{\abstractname} {\bf Abstract}
\begin{abstract} 
We consider the flow of a generalized Newtonian fluid through a thin porous medium of height $h_\varepsilon$ perforated with $\varepsilon$-periodically distributed solid cylinders of very small diameter $\varepsilon\delta_\varepsilon$, where the small parameters $\varepsilon, \delta_\varepsilon$ and $h_\varepsilon$ are devoted to tend to zero. We assume that the fluid is described by the 3D incompressible Stokes system with a non-linear power law viscosity of flow index $1<r<2$ (shear thinning). The particular case $h_\varepsilon=\sigma_\varepsilon$, where $\sigma_\varepsilon:=\varepsilon/\delta_\varepsilon^{{2-r\over r}}\to 0$, was recently published in (Anguiano and Su\'arez-Grau, \emph{Mediterr. J. Math.} (2021) 18:175). In this paper, we generalize previous study for any $h_\varepsilon$ and we provide a more complete description on the asymptotic behavior of non-Newtonian fluids in a thin porous medium composed by cylinders of small diameter. We prove that depending on the value of $\lambda:=\lim_{\varepsilon\to 0}\sigma_\varepsilon/h_\varepsilon\in [0,+\infty]$,   there exist three types of lower-dimensional asymptotic models: a non-linear Darcy law in the case $\lambda=0$, a non-linear Brinkman-type law in the case $\lambda\in (0,+\infty)$, and a non-linear Reynolds law in the case $\lambda=+\infty$. \end{abstract}

 {\small \bf AMS classification numbers: } 76A05,  76M50, 76A20, 76S05.
 
  \section {Introduction}\label{S1} 
  
This paper is devoted to derive macroscopic laws for generalized Newtonian fluids through periodic porous medium with small thickness, which is of a great importance in engineering applications, such as oil recovery, food processing, and material processing, see for instance Prat and Aga${\rm \ddot{e}}$sse \cite{Prat}, Qin and  Hassanizadeh \cite{Qin} and Yeghiazarian {\it et al.} \cite{Yeghiazarian}.  In this paper we consider generalized Newtonian fluids with the commonly used viscosity formula known as {\it power law} model. Thus, denoting the shear rate by $\mathbb{D}[u]=(Du+D^tu)/2$, with $u$ the velocity, the viscosity, as a function of the shear rate, is given by
$$\eta_r(\mathbb{D}[u])=\mu|\mathbb{D}[u]|^{r-2},\quad 1<r<2,\quad \mu>0,$$
where $\mu>0$ is the consistency, $r$ is the flow index (shear-thinning) and the matrix norm $|\cdot |$ is defined by $|\xi|^2=Tr(\xi\xi^t)$, $\xi\in\mathbb{R}^3$,  see for instance  Barnes et al. \cite{Barnes}, Bird et al. \cite{Bird} and Saramito \cite[Chapter 2]{Saramito}, for more details. 

Let us make a recollection of some previous results  in relation to the objective of this paper.  The case of the derivation of macroscopic laws for generalized Newtonian fluids with power law vicosity through a periodic porous medium $\Omega_\varepsilon\subset\mathbb{R}^3$ with a periodic arrangement of obstacles has been consider in Bourgeat {\it et al.} \cite{Bourgeat1, Bourgeat2} (see also Allaire \cite{Allaire1, Allaire0}, Sanchez-Palencia \cite{SanchezPalencia} and Tartar \cite{Tartar} for Newtonian fluids). The parameter   $\varepsilon$ is a small parameter related to the characteristic size of the obstacles and the period of the periodic structure. Thus, starting from the 3D power law Stokes system with body forces $f$ described by
\begin{equation}
\left\{
\begin{array}
[c]{r@{\;}c@{\;}ll}%
\medskip
\displaystyle -\mu\,{\rm div}\, \left(|\mathbb{D}\left[u_{\varepsilon}\right]|^{r-2}\mathbb{D}\left[u_{\varepsilon }\right] \right)+ \nabla p_{\varepsilon } &
= &
f\quad \hbox{in }\Omega_\varepsilon,\\
\medskip
{\rm div}\, u_{\varepsilon} & = & 0 \quad \hbox{in }\Omega_\varepsilon,\\
\medskip
u_\varepsilon&=&0\quad \hbox{on }\partial\Omega_\varepsilon,
\end{array}
\right. \label{Intro1}%
\end{equation}   and letting $\varepsilon$ tends to zero, by means of  the technique of two-scale convergence method, it was derived the following nonlinear Darcy type law in an $\varepsilon$-independent domain $\Omega\subset\mathbb{R}^3$:
$$
\left\{\begin{array}{l}
\medskip 
\displaystyle u(x)={1\over \mu}\mathcal{U}(f(x)-\nabla p(x))\quad\hbox{in }\Omega,\\
\medskip
{\rm div}\,u=0\quad\hbox{in }\Omega,\quad u\cdot n=0\quad\hbox{on }\partial\Omega,
\end{array}\right.
$$
with $n$ the outside normal vector to $\partial \Omega$.  The nonlinear function $\mathcal{U}:\mathbb{R}^3\to \mathbb{R}^3$ is called {\it permeability function} of the porous medium and is defined through the solutions $(v^\xi, \pi^\xi)$, $\xi\in\mathbb{R}^3$, of  3D cell problems  of power law type depending on the geometrical structure of the domain:
\begin{equation}\label{aux_problem_intro}
\left\{\begin{array}{rl}
\displaystyle
\medskip
-{\rm div}_{z}(|\mathbb{D}_{z}[v^{\xi}]|^{r-2}\mathbb{D}_{z}(v^{\xi}))+\nabla_{z}\pi^{\xi}=\xi &\hbox{in }Y\setminus Y_s,\\
\displaystyle
\medskip
{\rm div}_{z}v^{\xi}=0&\hbox{in }Y\setminus Y_s,\\
\displaystyle
\medskip
v^{\xi}=0&\hbox{in } Y_s,\\
\displaystyle
\medskip(v^{\xi},\pi^{\xi})&Y-\hbox{periodic},
\end{array}\right.
\end{equation}
where $Y\subset \mathbb{R}^3$ is the unitary reference cell and $Y_s\subset \overline Y$ is the reference obstacle. For a complete review of the homogenization of non-Newtonian fluids through a porous medium, we refer to Mikeli\'c \cite{Mikelic2, Mikelic1}.

A more general case of periodic porous medium is when the size of the solid obstacles is assumed to be much smaller than the period (the so-called  {\it tiny } or {\it small holes} in the mathematical literature, see Allaire \cite{Allaire2, Allaire3}, Brillard \cite{Brillard}, Diening {\it et al.} \cite{Diening}, Lu and Schwarzacher \cite{Lu1} or Lu \cite{Lu2}  for Newtonian fluids). In the case of power law fluids, denoting by $\varepsilon$ the period of the periodic porous structure and considering obstacles of size $\varepsilon \delta_\varepsilon$, with $\delta_\varepsilon$ a parameter depending on $\varepsilon$ and devoted to tends to zero with $\varepsilon$,  Fratrovi\'c and Maru${\rm \check{s}}$i\'c-Paloka \cite{MarusicPaloka}  found that, according to the behavior of the parameter
\begin{equation}\label{sigma3}\widehat   \sigma_\varepsilon={\varepsilon\over \delta_\varepsilon^{3-r\over r}},
\end{equation}
when $\varepsilon$ tends to zero, there are three different regimes. Namely,  if $\widehat   \sigma_\varepsilon\to +\infty$ (the case called {\it smaller holes}) the homogenized problem has a nonlinear Stokes form and, on the contrary, if $\widehat   \sigma_\varepsilon\to 0$ (the case called {\it large holes}),  the homogenized problem is a nonlinear Darcy law. Moreover,  when $\widehat   \sigma_\varepsilon\to \widehat   \sigma_0>0$ (called {\it critical case}), it was   developed  by obtaining  the following Brinkman problem by means of $\Gamma$-convergence:
$$
\left\{\begin{array}{rl}
\medskip 
\displaystyle
-\mu \,{\rm div}\left(| \mathbb{D}[u]|^{r-2}\mathbb{D}[u]\right)+{\mu\over \widehat   \sigma_0^r} \mathcal{G}(u)+\nabla p=f& \hbox{in }\Omega,\\
\medskip
{\rm div}\,u=0&\hbox{in }\Omega,\\
\medskip
\displaystyle
u=0&\hbox{on }\partial\Omega,
\end{array}\right.
$$
where the nonlinear function $\mathcal{G}:\mathbb{R}^3\to \mathbb{R}^3$, called the {\it drag force function} of the porous medium,  provides the value of the drag force on the reference obstacle $Y_s$ by the exterior fluid flow (see  also  Fratrovi\'c and Maru${\rm \check{s}}$i\'c-Paloka \cite{MarusicPaloka2} for more details):
\begin{equation}\label{aux_problem_intro2}
\left\{\begin{array}{rl}
\displaystyle
\medskip
-{\rm div}_{z}(|\mathbb{D}_{z}[v^{\xi}]|^{r-2}\mathbb{D}_{z}(v^{\xi}))+\nabla_{z}\pi^{\xi}=\xi &\hbox{in }\mathbb{R}^3\setminus Y_s,\\
\displaystyle
\medskip
{\rm div}_{z}v^{\xi}=0&\hbox{in }\mathbb{R}^3\setminus Y_s,\\
\displaystyle
\medskip
v^{\xi}=0&\hbox{on } \partial Y_s,\\
\displaystyle
\lim_{|z|\to+\infty}v^\xi=\xi.
\end{array}\right.
\end{equation}

On the other hand, the derivation of macroscopic laws for fluids in porous domains with small thickness (the so-called {\it thin porous medium})  is attracting much attention, see for instance Almqvist {\it et al.} \cite{Almqvist}, Anguiano \cite{Anguiano2}, Anguiano and Bunoiu \cite{Ang-Bun2}, Anguiano {\it et al.} \cite{Anguiano_Bonnivard_SG}, Anguiano and Su\'arez-Grau \cite{Anguiano_SuarezGrau, Anguiano_SuarezGrau2, Anguiano_SG_NHM, Anguiano_SG_Lower, Anguiano_SG_sharp}, Fabricius {\it et al.} \cite{Fabricius, Fabricius2}, Fabricius and Gahn \cite{Fabricius3}, Forslund {\it et al.} \cite{Forslund},  Jouybari and  Lundstr$\ddot{\rm o}$m \cite{Jouybari}, Mei and Vernescu \cite{Mei}, Su\'arez-Grau \cite{SuarezGrau1, SuarezGrau2} or Zhengan and  Hongxing \cite{Zhengan}. A thin porous medium can be defined as a bounded perforated 3D domain confined between two parallel plates, where the distance between the plates is very small and the perforation consists of periodically distributed solid cylinders which connect the plates in perpendicular direction. Denoting $\omega_\varepsilon\subset \mathbb{R}^2$ a perforated domain with obstacles of size and period $\varepsilon$, we  define the three-dimensional thin porous medium by $\Omega_\varepsilon=\omega_\varepsilon \times (0, h_\varepsilon)$, where $h_\varepsilon$ is the height of the domain. It happens that in this type of domains the asymptotic behavior depends on the relationship between the period and the height, that is, on the asymptotic value of the ratio  $\varepsilon/h_\varepsilon$ when $\varepsilon$ tends to zero. Thus, in the case of generalized Newtonian fluid with power law viscosity,  see Anguiano and Su\'arez-Grau \cite{Anguiano_SuarezGrau},  different lower-dimensional types of the nonlinear Darcy law in $\omega\subset\mathbb{R}^2$ were derived by using an adaptation of the unfolding method, which can be written as  follows:
$$
\left\{\begin{array}{l}
\medskip\displaystyle
U_{av}'(x')={1\over \mu}\mathcal{U}_\lambda\left(f'(x')-\nabla_{x'} p(x')\right),\quad U_{av,3}\equiv 0\quad \hbox{in }\omega,\\
\medskip
{\rm div}_{x'}U_{av}'=0\quad\hbox{in }\omega,\quad U_{av}'\cdot n'=0\quad\hbox{on }\partial\omega.
\end{array}\right.
$$
where  $U_{av}=(U_{av}',U_{av,3})$, with $U_{av}'=(U_{av,1},U_{av,2})$, is the average in height of the limit velocity, $f'=(f_1,f_2)$  is the horizontal body forces,  $x'=(x_1,x_2)\in \omega$ denotes the horizontal variables, and $n'$ is the outside normal vector to $\partial \omega$. The permeability function $\mathcal{U}_\lambda:\mathbb{R}^2\to \mathbb{R}^2$ is defined through the solutions of auxiliary problems of power law type depending on the value $\lambda:=\lim_{\varepsilon\to 0}\varepsilon/h_\varepsilon\in [0,+\infty]$:
\begin{itemize}

\item[--] When $\varepsilon/h_\varepsilon\to \lambda\in (0,+\infty)$, corresponding to the case when the cylinder height is proportional to the interspatial distance (called  {\it proportionally thin porous medium}), the permeability function is calculated by solving $3D$ power law Stokes local problems, satisfied by  $(v^{\xi'}, \pi^{\xi'})$ with $\xi'\in\mathbb{R}^2$, depending on the parameter $\lambda$:
\begin{equation}\label{aux_problem_intro_PTPM}
\left\{\begin{array}{rl}
\displaystyle
\medskip
-{\rm div}_{\lambda}(|\mathbb{D}_{\lambda}[v^{\xi'}]|^{r-2}\mathbb{D}_{\lambda}(v^{\xi'}))+\nabla_{\lambda}\pi^{\xi'}=\xi' &\hbox{in }Y\setminus Y_s,\\
\displaystyle
\medskip
{\rm div}_{\lambda}v^{\xi'}=0&\hbox{in }Y\setminus Y_s,\\
\displaystyle
\medskip
v^{\xi'}=0&\hbox{in } Y_s,\\
\displaystyle
\medskip(v^{\xi'},\pi^{\xi'})&Y'-\hbox{periodic}.
\end{array}\right.
\end{equation}
Here, we denote $Y'$ is the unitary reference cell and $Y_s'\subset \overline Y'$ the reference obstacle in $\mathbb{R}^2$ and so,   $Y=Y'\times (0,1)$ denotes the unitary reference cell and $Y_s=Y_s'\times (0,1)$ the corresponding reference obstacle in $\mathbb{R}^3$. Moreover, $\mathbb{D}_\lambda[v^{\xi'}]=\mathbb{D}_{z'}[v^{\xi'}]+\lambda\partial_{z_3}[v^{\xi'}]$, $\nabla_{\lambda}v^{\xi'}=(\nabla_{z'}v^{\xi'}, \lambda\partial_{z_3}v^{\xi'})^t $ and ${\rm div}_\lambda v^{\xi'} =\partial_{z_1}v^{\xi'}_1+\partial_{z_2}v_2^{\xi'}+\lambda\partial_{z_3}v^{\xi'}_3$.

\item[--] When $\varepsilon/h_\varepsilon\to 0$, corresponding to the case when the cylinder height is much larger than the interspatial distance (called  {\it homogeneously thin porous medium}), the permeability function is calculated by solving a purely $2D$ power law Stokes local problem, satisfied by  $(v^{\xi'}, \pi^{\xi'})$ with $\xi'\in\mathbb{R}^2$:
\begin{equation}\label{aux_problem_intro_HTPM}
\left\{\begin{array}{rl}
\displaystyle
\medskip
-{\rm div}_{z'}(|\mathbb{D}_{z'}[v^{\xi'}]|^{r-2}\mathbb{D}_{z'}(v^{\xi'}))+\nabla_{z'}\pi^{\xi'}=\xi' &\hbox{in }Y'\setminus Y_s',\\
\displaystyle
\medskip
{\rm div}_{z'}v^{\xi'}=0&\hbox{in }Y'\setminus Y_s',\\
\displaystyle
\medskip
v^{\xi'}=0&\hbox{in } Y_s',\\
\displaystyle
\medskip(v^{\xi'}, \pi^{\xi'})&Y'-\hbox{periodic}.
\end{array}\right.
\end{equation}

\item[--] When $\varepsilon/h_\varepsilon\to +\infty$, corresponding to the case when the cylinder height is much smaller than the interspatial discance  (called  {\it very thin porous medium}), the permeability function is calculated by solving $2D$ Hele-Shaw local problems for $\pi^{\xi'}$ with $\xi'\in\mathbb{R}^2$:
\begin{equation}\label{aux_problem_intro_VTPM}
\left\{\begin{array}{rl}
\displaystyle
\medskip
-{\rm div}_{z'}(|\xi'+\nabla_{z'}\pi^{\xi'}|^{r-2}(\xi'+\nabla_{z'}\pi^{\xi'}))=0 &\hbox{in }Y'\setminus Y_s',\\
\displaystyle
\medskip
 (|\xi'+\nabla_{z'}\pi^{\xi'}|^{r-2}(\xi'+\nabla_{z'}\pi^{\xi'}))\cdot n'=0&\hbox{on } \partial Y_s'.
\end{array}\right.
\end{equation}
\end{itemize}
It is worth noting that the critical ratio $\varepsilon / h_\varepsilon$ appears by obtaining a sharp Poincar\'e-Korn inequality in the thin porous media $\Omega_\varepsilon$:
$$\|v\|_{L^r(\Omega_\varepsilon)^3}\leq K_\varepsilon  \|\mathbb{D}[v]\|_{L^r(\Omega_\varepsilon)^{3\times 3}}\quad \forall\, v\in W^{1,r}_0(\Omega_\varepsilon)^3,$$
where  $K_\varepsilon =C \varepsilon$ in the cases of proportionally and homogeneously thin porous medium, and $K_\varepsilon=C h_\varepsilon$ in the case of a very thin porous medium,   where $C>0$ is independent of $\varepsilon$, see \cite[Lemma 4.1]{Anguiano_SuarezGrau} for more details.

The goal of this paper is to provide a more complete description of the study the homogenization of generalized Newtonian fluid with power law viscosity in a thin porous medium. Here, we consider the case of a thin domain perforated by an array of periodically distributed cylinders with small diameter. More precisely, we consider the thin porous medium defined by $\Omega_\varepsilon=\omega_\varepsilon\times (0,h_\varepsilon)$, with $\varepsilon$ describing the period of the periodic porous structure,  $\varepsilon \delta_\varepsilon$ describing the size of the diameters of the cylinders, and $h_\varepsilon$ describing the height of the cylinders (with $\delta_\varepsilon$ and $h_\varepsilon$  parameters depending on $\varepsilon$ and devoted to tend to zero with $\varepsilon$). As a result, we prove that, according to the behavior of the parameter
\begin{equation}\label{sigma2}\sigma_\varepsilon={\varepsilon\over \delta_\varepsilon^{2-r\over r}}\quad (\hbox{assuming }\sigma_\varepsilon \to 0),
\end{equation}
 with respect to the height parameter $h_\varepsilon$,  there are three different regimes when $\varepsilon$ tends to zero depending on the value $\lambda:=\lim_{\varepsilon\to 0}\sigma_\varepsilon/h_\varepsilon\in [0,+\infty]$:
\begin{itemize}
\item[--] When  $\sigma_\varepsilon/h_\varepsilon\to \lambda\in (0,+\infty)$, then we derive  the following lower-dimensional nonlinear Brinkman-type law posed  in $\Omega=\omega\times (0,1)\subset\mathbb{R}^3$, depending on $\lambda$, satisfied by the limit velocity $u$ and limit pressure $p$:
$$
\left\{\begin{array}{rl}\medskip
\displaystyle-\mu\,\lambda^r\,2^{-{r\over 2}}\partial_{y_3}(|\partial_{y_3}u'|^{r-2}\partial_{y_3}u'(x',y_3))+\mu \mathcal{G}(u'(x',y_3))+\nabla_{x'} p(x')=f'(x')  &\hbox{in }\Omega,\\
\medskip
u_3\equiv 0, &\\
\medskip
\displaystyle
u'(x',0)=u'(x',1)=0& \hbox{in }\omega,\\
\medskip
\displaystyle
{\rm div}_{x'}\left(\int_0^1u'(x',y_3)\,dy_3\right)=0& \hbox{in }\omega,\\
\medskip
\displaystyle
 \left(\int_0^1u'(x',y_3)\,dy_3\right)\cdot n'=0 & \hbox{on }\partial\omega,
\end{array}\right.
$$
where $x'=(x_1, x_2)\in \omega$ denotes  the horizontal variables, $y_3\in (0,1)$ is the vertical variable, and  $\mathcal{G}:\mathbb{R}^2\to \mathbb{R}^2$ is the nonlinear  drag force function, which is defined through the solution of lower-dimensional auxiliary  exterior  problems of power law type depending on the geometrical structure of the obstacles (see Theorem \ref{MainTheorem}$-(i)$).

\item[--]  When  $\sigma_\varepsilon/h_\varepsilon\to 0$, then we derive  the following lower-dimensional nonlinear Darcy  law in $\omega\subset\mathbb{R}^2$:
$$
\left\{\begin{array}{rl}
\medskip
\displaystyle    U_{av}'(x')={1\over \mu^{1\over r-1}}\mathcal{G}^{-1}\left(f'(x')-\nabla_{x'}p(x')\right),\quad U_{av, 3}\equiv 0&\hbox{in }\omega,\\
\medskip
\displaystyle {\rm div}_{x'} U_{av}' (x')=0&\hbox{in }\omega,\\
\medskip
\displaystyle  U_{av}'(x') \cdot n'=0&\hbox{on }\partial\omega,
\end{array}\right.
$$
where $U_{av}$ is the average in height of the limit velocity and $\mathcal{G}^{-1}:\mathbb{R}^2\to\mathbb{R}^2$   is the inverse of the drag force function  $\mathcal{G}:\mathbb{R}^2\to \mathbb{R}^2$ (see Theorem \ref{MainTheorem}$-(ii)$).

\item[--]  When  $\sigma_\varepsilon/h_\varepsilon\to +\infty$, then we deduce that there are no effects of the microstructure in the limit, so we derive the following nonlinear and reduced Stokes-type system satisfied by the limit velocity $u$ and limit pressure $p$:
$$
\left\{\begin{array}{rl}\medskip
\displaystyle-\mu\,\lambda^r\,2^{-{r\over 2}}\partial_{y_3}(|\partial_{y_3}u'|^{r-2}\partial_{y_3}u'(x',y_3))+\nabla_{x'} p(x')=f'(x')  &\hbox{in }\Omega,\\
\medskip
u_3\equiv 0, &\\
\medskip
\displaystyle
u'(x',0)=u'(x',1)=0& \hbox{in }\omega,\\
\medskip
\displaystyle
{\rm div}_{x'}\left(\int_0^1u'(x',y_3)\,dy_3\right)=0& \hbox{in }\omega,\\
\medskip
\displaystyle
 \left(\int_0^1u'(x',y_3)\,dy_3\right)\cdot n'=0 & \hbox{on }\partial\omega,
\end{array}\right.
$$
which leads to the classical lower-dimensional Reynolds problem for power law fluids in a thin domain without obstacles (see Theorem \ref{MainTheorem}$-(iii)$).
\end{itemize}

We comment that this study can be framed into the   regime  of ``{\it large holes}" obtained in Fratrovi\'c and Maru${\rm \check{s}}$i\'c-Paloka \cite{MarusicPaloka} (i.e. the case $\widehat  \sigma_\varepsilon$ given by (\ref{sigma3}) when $\widehat  \sigma_\varepsilon\to 0$). However,  we have to take into account that, in this case, we have cylinder-shaped obstacles and, moreover, that the height of the domain is  of order smaller than one. In this sense, due to the shape of the obstacles,  the exponent ${3-r\over r}$ given in $\widehat  \sigma_\varepsilon$ changes to the exponent ${2-r\over r}$ given in $\sigma_\varepsilon$. Moreover, the case $\sigma_\varepsilon/h_\varepsilon\to 0$, i.e. when the height of the domain is greater than $\sigma_\varepsilon$, gives a nonlinear Darcy law as in  \cite{MarusicPaloka}, but with a lower-dimension due to the small height of the domain. As the   height of the domain is reduced, we pass from a Darcy regime ($\sigma_\varepsilon/h_\varepsilon\to 0$) to a Reynolds regime ( $\sigma_\varepsilon/h_\varepsilon\to +\infty$) in which  the microstructure of the thin porous medium is not detected. Between the two regimes, a  critical size appears ($\sigma_\varepsilon/h_\varepsilon\to \lambda\in (0,+\infty)$),  where a lower-dimensional Brinkman-type law is derived.
\\

We  remark that this problem has already been studied in Anguiano and Su\'arez-Grau \cite{Anguiano_SG_Lower} in the critical case by assuming $h_\varepsilon\equiv \sigma_\varepsilon$, i.e. when $\lambda=1$, which simplifies the difficulties that appear in the present study. Thus, in this paper we carry out a more complete study by adapting the combination of reduction of dimension techniques, monotonicity arguments and the version of the unfolding method depending on the parameters $\varepsilon$ and $\delta_\varepsilon$ introduced in  \cite{Anguiano_SG_Lower}, by taking into account now the new parameter $h_\varepsilon$. It is worth mentioning that we are considering the case of shear thinning power law fluids, i.e. when the flow index satisfies $1<r<2$, because the version of the unfolding method introduces such a restriction, as explained in \cite{Anguiano_SG_Lower}. With respect to the mentioned study, the following novelties are introduced:
\begin{itemize}
\item[--] The derivation of the critical ratio $\sigma_\varepsilon / h_\varepsilon$ by obtaining a  sharp Poincar\'e-Korn inequality in the thin porous medium $\Omega_\varepsilon$, see Lemma \ref{Poincare0} for more details. Thus, we prove that for every $v\in W^{1,r}_0(\Omega_\varepsilon)^3$, it holds
$$\|v\|_{L^r(\Omega_\varepsilon)^3}\leq K_\varepsilon \|\mathbb{D}[v]\|_{L^r(\Omega_\varepsilon)^{3\times 3}},$$
where  $K_\varepsilon=C\,\sigma_\varepsilon$ in the cases $\sigma_\varepsilon\approx h_\varepsilon$ or $\sigma_\varepsilon\ll h_\varepsilon$, and $K_\varepsilon=C\,h_\varepsilon$ in the case $\sigma_\varepsilon \gg h_\varepsilon$,  with $C>0$ independent of $\varepsilon$.

\item[--] An inverse of the divergence operator on thin and perforated domains, which is necessary to derive optimal estimates for the pressure, see Lemma \ref{lem:Duvjnak} for more details. More precisely, we prove that for every $g\in L^{r}(\Omega_{\varepsilon})$,there exists $\varphi=\varphi(g)\in W^{1,r}(\Omega_\varepsilon)^3$,  with $\varphi=0$ on the exterior boundary, such that
$$
{\rm div}\,\varphi=g\ \hbox{ in }\Omega_{\varepsilon}\quad\hbox{such that}\quad \|\varphi\|_{L^r(\Omega_{\varepsilon})^{3}}\leq C\|g\|_{L^r(\Omega_{\varepsilon})},\quad \|D\varphi\|_{L^r(\Omega_{\varepsilon})^{3\times 3}}\leq  C_\varepsilon \|g\|_{L^r(\Omega_{\varepsilon})},
$$
where $C_\varepsilon=C\,\sigma_\varepsilon^{-1}$ in the cases $\sigma_\varepsilon\approx h_\varepsilon$ or $\sigma_\varepsilon\ll h_\varepsilon$, and $C_\varepsilon=C\,h_\varepsilon^{-1}$ in the case $\sigma_\varepsilon \gg h_\varepsilon$,  with $C>0$ independent of $\varepsilon$.

\item[--] Some suitable compactness results corresponding to each case, which let us to derive the different homo\-genized problems, see Lemmas \ref{LemmaConvtulde} and \ref{lem:compactnesshat} for more details.
\end{itemize}

We think that this theoretical study provides a quite complete description of the asymptotic behavior of generalized Newtonian fluids with power law viscosity through a thin porous medium composed by cylinders of very small diameter. Since the obtained models are amenable for the numerical simulations, we believe that it could prove useful in the engineering practice as well.

Finally, we comment the structure of the paper. In Section \ref{sec:setting} we introduce the domain and give the main result (Theorem \ref{MainTheorem}). In Section 3,  we develop the proof of the main result in different subsections. In Subsection \ref{S2} we prove the sharp Poincar\'e-Korn inequality and then, some {\it a priori} estimates for the velocity are derived.  In Subsection \ref{S2p} we introduce the inverse of the divergence operator on a thin porous medium and then, we deduce the {\it a priori} estimates for the pressure. In Subsection \ref{sec:unf} we recall the version of the unfolding method for a domain perforated by cylinders of  small diameter. Some compactness results, which are the main keys when we will pass to the limit later, are addressed in Section \ref{sec:comp}. In Subsection \ref{sec:proofthm} we give the proof of the Theorem \ref{MainTheorem}. We finish the paper with a list of references.

\section{Setting of the problem and main result}\label{sec:setting}
{\bf Geometrical setting.} Assume that $\omega$ is a smooth, bounded, connected set in $\mathbb{R}^2$. Then, a periodic porous medium is defined by the domain $\omega$ with an associated microstructure,
of periodic cell $Y^{\prime}=[-1/2,1/2]^2$, which is made of two complementary parts: 
the fluid part $Y^{\prime}_{f}$, and the solid part $Y^{\prime}_{s}$, such that $Y^{\prime}_f  \bigcup Y^{\prime}_{s}=Y^\prime$ and $Y^{\prime}_f  \bigcap Y^{\prime}_{s}=\emptyset$. Here we assume that $Y_{s}^{\prime}$ 
is a closed subset of $\overline Y^\prime$ with a smooth boundary (with Lipschitz boundary) $\partial Y_{ s}^\prime$, such that $Y_{ s}^\prime$ is strictly included  in $\overline Y^\prime$.  We denote $Y=Y'\times (0,1)\subset\mathbb{R}^3$ and so $Y_f=Y'_f\times (0,1)$ and $Y_s=Y'_s\times (0,1)$.

To define the cylinders with small diameter, we consider the positive parametes $\varepsilon$ and $\delta_\varepsilon$,  smaller than one, where $\delta_\varepsilon$ is such that $\delta_\varepsilon \to 0$ as $\varepsilon \to 0$. Using both parameters, we define the parameter
\begin{equation}\label{sigma}
\sigma_{\varepsilon}={\varepsilon\over \delta_\varepsilon^{2-r\over r}},
\end{equation}
and we consider that it holds the following condition 
\begin{equation}\label{sigma0}
\lim_{\varepsilon\to 0}\sigma_{\varepsilon}=0\quad\hbox{or equivalently }\quad 0<\varepsilon\ll \delta_\varepsilon^{2-r\over r}.
\end{equation}
This assumption leads to a solid cylinders of very small diameter $\varepsilon\delta_\varepsilon$ satisfying
$$0<\varepsilon^{2\over 2-r}\ll \varepsilon\delta_\varepsilon,$$
which, according to \cite{MarusicPaloka}, it could be called ``{\it large diameter cylinders}".
Also, since the diameter $\varepsilon\delta_\varepsilon$ is still smaller than the inter-obstacle distance $\varepsilon$, then (\ref{sigma0}) yields
$$\varepsilon\ll \sigma_\varepsilon\ll 1.$$

To define the thickness of the cylinders, we consider the positive parameter  $h_\varepsilon$, smaller than one, where $h_\varepsilon$ tends to zero as $\varepsilon\to 0$. 
We consider a thin porous medium $\Omega_{\varepsilon}$ of thickness $h_\varepsilon$, which is perforated by solid cylinders with diameter of size $\varepsilon \delta_\varepsilon$ and distributed periodically with period $\varepsilon$.

To define the microstructure of the domain $\omega$, we set $$Y^{\prime}_{\delta_\varepsilon f}=Y^{\prime}\setminus \delta_\varepsilon  Y'_s,$$
such that the domain $\omega$ is covered by a regular mesh of size $\varepsilon$, i.e. for $k^{\prime}\in \mathbb{Z}^2$, each cell $Y^{\prime}_{k^{\prime},{\varepsilon}}={\varepsilon}k^{\prime}+{\varepsilon}Y^{\prime}$ is divided in a fluid part $Y^{\prime}_{\delta_\varepsilon f_{k^{\prime}},{\varepsilon}}$ and a solid part $ Y^{\prime}_{ \delta_\varepsilon s_{k^{\prime}},\varepsilon}$, where $ Y^{\prime}_{ \delta_\varepsilon s_{k^{\prime}},\varepsilon}$ denotes the complement in $Y'_{k',\varepsilon}$ of the set $Y^{\prime}_{\delta_\varepsilon f_{k^{\prime}},{\varepsilon}}$. We observe that $Y^{\prime}_{k^{\prime},{\varepsilon}}$ is similar to the unit cell $Y^{\prime}$ rescaled to size ${\varepsilon}$. 

As consequence, $Y$  is divided in a fluid part $Y_{\delta_\varepsilon f}$ and a solid part $Y_{\delta_\varepsilon s}$, and consequently $Y_{k^{\prime},{\varepsilon}}=Y^{\prime}_{k^{\prime},{\varepsilon}}\times (0,1)\subset \mathbb{R}^3$, which is also divided in a fluid part $Y_{\delta_\varepsilon f_{k^{\prime}},{\varepsilon}}$ and a solid part $ Y_{\delta_\varepsilon s_{k^{\prime}},\varepsilon}$. 
\begin{figure}[h!]
\begin{center}
\includegraphics[width=4cm]{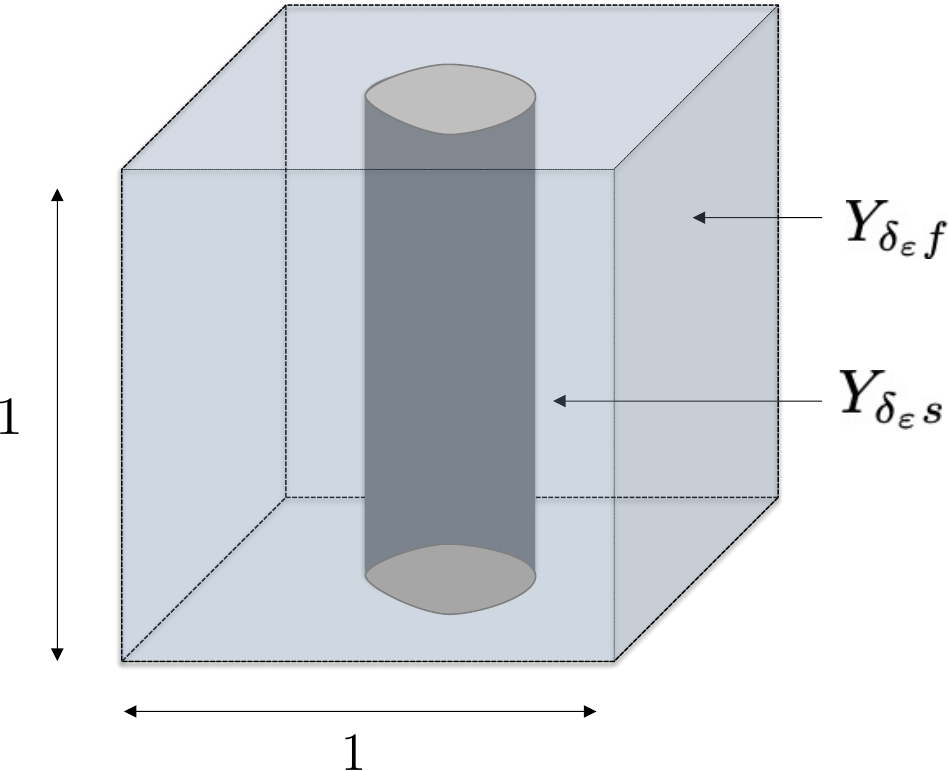}
\hspace{2.5cm}
\raisebox{.1\height}{\includegraphics[width=3.5cm]{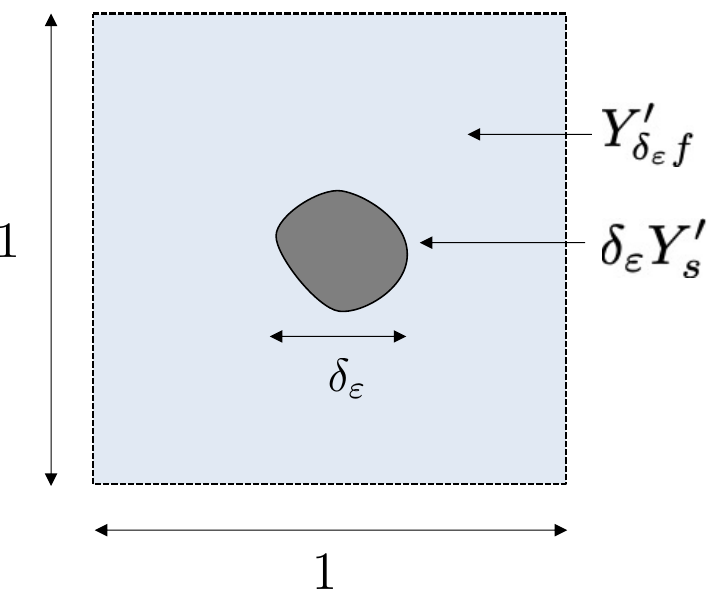}}
\end{center}
\vspace{-0.4cm}
\caption{View of the 3D reference cells  $Y$ (left) and the 2D reference cell $Y'$ (right).}
\label{fig:cell}
\end{figure}

We denote by $\tau( Y'_{ \delta_\varepsilon s_{k'},\varepsilon})$ the set of all translated images of
$   Y'_{\delta_\varepsilon s_{k'},\varepsilon}$. The set $\tau(   Y'_{\delta_\varepsilon s_{k'},\varepsilon})$ represents the solids in $\mathbb{R}^2$.
The fluid part of the bottom $\omega_{\varepsilon}\subset \mathbb{R}^2$ of the porous medium is defined by
$\omega_{\varepsilon }=\omega\backslash\bigcup_{k^{\prime}\in \mathcal{K}_{\varepsilon}}     Y^{\prime}_{\delta_\varepsilon s_{k^{\prime}},
{\varepsilon}}$, where 
$\mathcal{K}_{\varepsilon}=\{k^{\prime}\in \mathbb{Z}^2: Y^{\prime}_{k^{\prime}, {\varepsilon}} \cap \omega \neq \emptyset \}$.    

The whole fluid part $\Omega_{\varepsilon }\subset \mathbb{R}^3$ in the thin porous medium is then defined by 
\begin{equation}\label{Dominio1}
\Omega_{\varepsilon }=\{  (x',x_3)\in \mathbb{R}^2\times \mathbb{R}: x'\in\omega_{\varepsilon }\,,\  0<x_3<h_\varepsilon \}.
\end{equation}
We make the assumption that the solids  $\displaystyle \tau(  Y'_{ \delta_\varepsilon s_{k'},\varepsilon})$ do not intersect the boundary
$\partial \omega$. We define $Y_{k^{\prime},{\varepsilon}}^{h_\varepsilon}=Y'_{k',\varepsilon}\times (0,h_\varepsilon)$  and so, $Y^{h_\varepsilon}_{ \delta_\varepsilon s_{k'},\varepsilon}=  Y'_{\delta_\varepsilon  s_{k'},\varepsilon}\times (0,h_\varepsilon)$ and   $ Y^{h_\varepsilon}_{ \delta_\varepsilon f_{k'},\varepsilon}=  Y'_{\delta_\varepsilon  f_{k'},\varepsilon}\times (0,h_\varepsilon)$. Denote by $S_{\varepsilon }$
the set of the solids contained in $\Omega_{\varepsilon }$. Then, $S_{\varepsilon }$ is a finite  union of solids, i.e. $S_{\varepsilon }=\bigcup_{k^{\prime}\in \mathcal{K}_{\varepsilon}}   Y^{h_\varepsilon}_{ \delta_\varepsilon s_{k'},\varepsilon}.$ 

\begin{figure}[h!]
\begin{center}
\includegraphics[width=4cm]{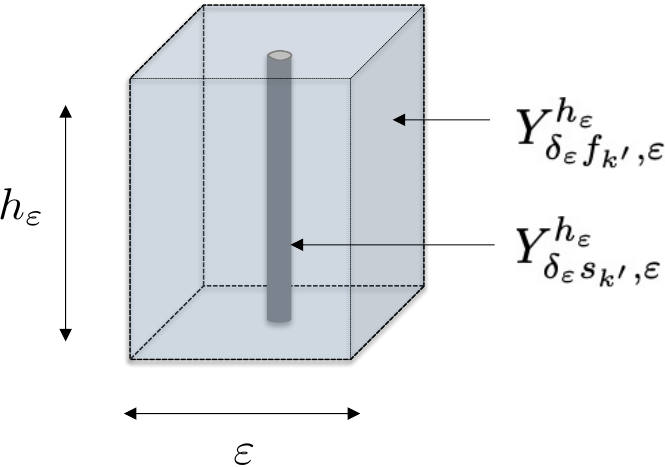}
\hspace{2.5cm}
\raisebox{.25\height}{\includegraphics[width=3.5cm]{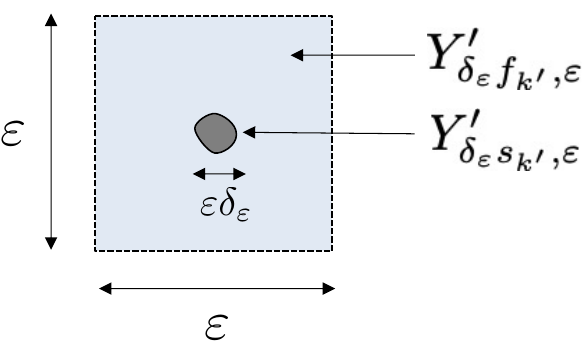}}
\end{center}
\vspace{-0.4cm}
\caption{View of the 3D rescaled cell $Y_{k',\varepsilon}^{h_\varepsilon}$ (left) and the 2D rescaled cell $Y'_{k',\varepsilon}$ (right).}
\label{fig:cellrescaled}
\end{figure}
We define the rescaled domain  $\widetilde{\Omega}_{\varepsilon }=\omega_{\varepsilon }\times (0,1)$, i.e. the domain perforated by cylinders with size one.  We observe that $\widetilde{\Omega}_{\varepsilon }=\Omega\backslash \bigcup_{k^{\prime}\in \mathcal{K}_{\varepsilon}}
   Y_{\delta_\varepsilon  s_{k^{\prime}}, {\varepsilon}},$ and we define $T_{\varepsilon }=\bigcup_{k^{\prime}\in \mathcal{K}_{\varepsilon}}
    Y_{\delta_\varepsilon s_{k^\prime}, \varepsilon}$ as the set of the solids contained in $\widetilde \Omega_{\varepsilon}$.

We also define the thin domain without microstucture by $Q_{\varepsilon }=\omega\times (0,h_\varepsilon)$ and the domain with height one without microstructure $\Omega=\omega\times (0,1)$.  
 
 \begin{figure}[h!]
\begin{center}
\includegraphics[width=8cm]{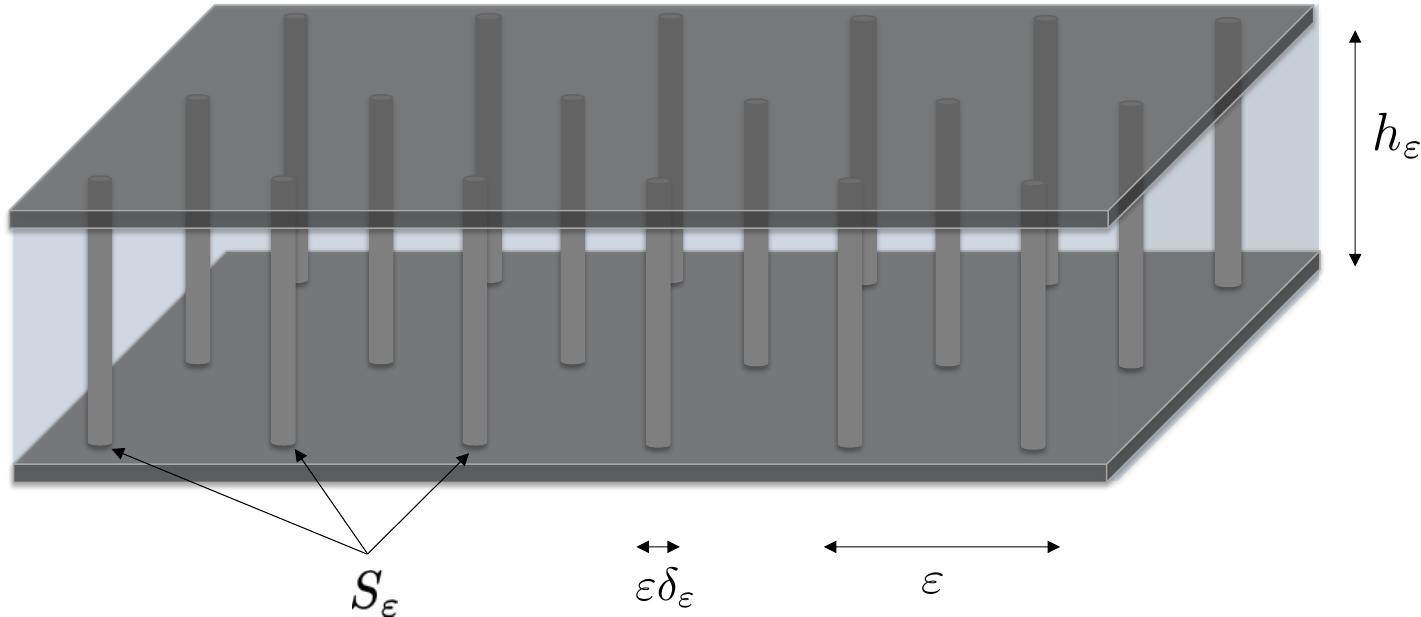}
\qquad
\raisebox{.25\height}{\includegraphics[width=8cm]{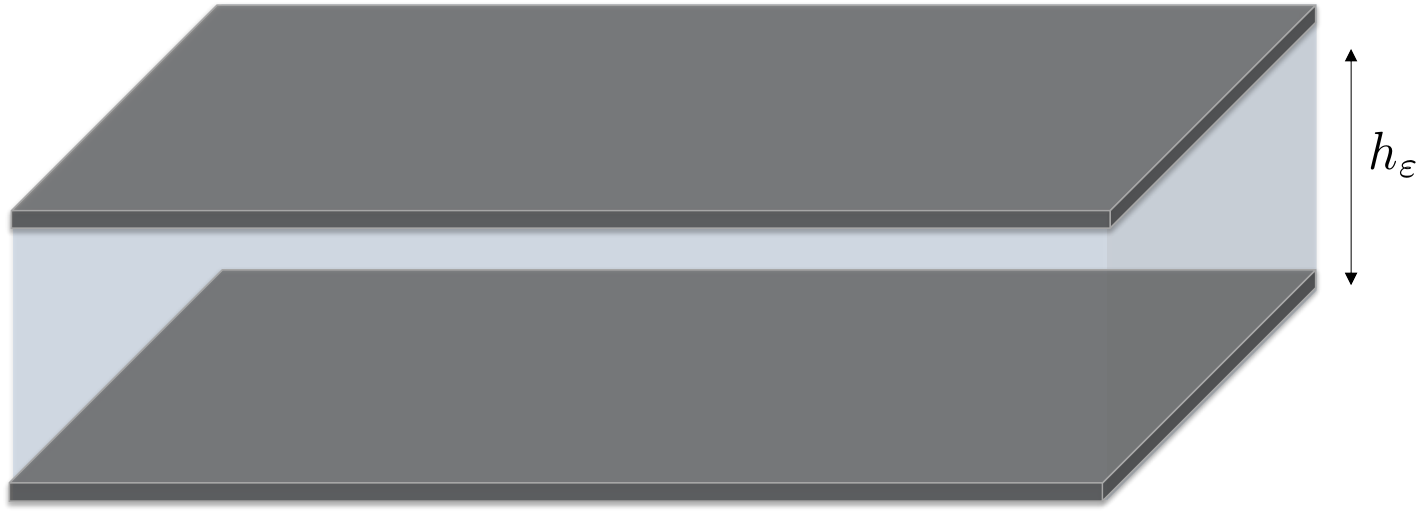}}
\end{center}
\vspace{-0.4cm}
\caption{View of the thin porous medium $\Omega_{\varepsilon }$ (left) and the domain without microstructure $Q_{\varepsilon }$ (right).}
\label{fig:cell}
\end{figure}
Following previous notation, we remark that along the paper,  the points $x\in\mathbb{R}^3$ will be decomposed as $x=(x^{\prime},x_3)$ with $x^{\prime}=(x_1,x_2)\in \mathbb{R}^2$, $x_3\in \mathbb{R}$. We also use the notation $x^{\prime}$ to denote a generic vector of $\mathbb{R}^2$.
\\

{\bf Setting of the problem.}  In the domain $\Omega_{\varepsilon }$ described previously,  we consider the following  Stokes system with power law viscosity 
\begin{equation}
\left\{
\begin{array}
[c]{r@{\;}c@{\;}ll}%
\medskip
\displaystyle -\mu\,{\rm div}\, \left(|\mathbb{D}\left[u_{\varepsilon }\right]|^{r-2}\mathbb{D}\left[u_{\varepsilon }\right] \right)+ \nabla p_{\varepsilon } &
= f&
 \text{\ in \ }\Omega_{\varepsilon },\\
\medskip
{\rm div}\, u_{\varepsilon } & = 0&  \text{\ in \ }\Omega_{\varepsilon },\\
\medskip
u_{\varepsilon }& = 0&  \text{\ on \ } \partial Q_{\varepsilon}\cup S_\varepsilon,
\end{array}
\right. \label{1}%
\end{equation}
where $u_{\varepsilon }$ is the velocity field and $p_{\varepsilon }$ is the pressure. In system (\ref{1}), we assume that the  body forces $f$ is of the usual form  in thin domain, given as follows
\begin{equation}\label{sourcef}
f(x)=(f^{\prime}(x^{\prime}),0), \text{\ a.e. \ }x\in \Omega,
\end{equation}
where $f$ is assumed in $L^{r^\prime}( \omega)^2$,  with  $r^\prime=r/(r-1)$ the conjugate exponent of $r$.  

Under previous assumptions, the classical theory gives the existence of a unique solution $(u_{\varepsilon },p_{\varepsilon })\in W_0^{1,r}(\Omega_{\varepsilon })^3\times L^{r^\prime}(\Omega_{\varepsilon})$ with $1<r<+\infty$,  see Lions \cite{Lions2}.  This solution is unique up to an additive constant for $p_{\varepsilon }$, i.e. it is unique if we consider the corresponding equivalence class $p_{\varepsilon}\in L^{r^\prime}(\Omega_{\varepsilon })/\mathbb{R}$.

Our objetive is to study the asymptotic behavior of the sequence of solution ($u_{\varepsilon }, p_{\varepsilon })$ of system (\ref{1}) when $\varepsilon$, $\delta_\varepsilon$ and $h_\varepsilon$ tend to zero. To do this, we consider  the dilatation in the variable $x_3$ given by
\begin{equation}\label{dilatacion}
y_3=\frac{x_3}{h_\varepsilon},
\end{equation}
in order to have the functions defined in the open set with fixed height $\widetilde\Omega_{\varepsilon }$. Thus, we define the dilated functions   $\widetilde {u}_{\varepsilon }\in W_0^{1,r}(\widetilde{\Omega}_{\varepsilon })^3$, $\widetilde {p}_{\varepsilon }\in L^{r^\prime}(\widetilde{\Omega}_{\varepsilon })/\mathbb{R}$ by 
$$\widetilde {u}_{\varepsilon }(x^{\prime},y_3)=u_{\varepsilon }(x^{\prime},h_\varepsilon y_3),\text{\ \ }\widetilde {p}_{\varepsilon }(x^{\prime},y_3)=p_{\varepsilon }(x^{\prime},h_\varepsilon y_3), \text{\ \ } a.e.\text{\ } (x^{\prime},y_3)\in \widetilde{\Omega}_{\varepsilon }.$$
Let us introduce some notation which will be useful in the following. For a vectorial function $v=(v',v_3)$ and a scalar function $w$, we will denote 
\begin{equation}\label{derivada_def}\mathbb{D}_{x^\prime}\left[v\right]=\frac{1}{2}\left(\begin{array}{ccc}
2\partial_{x_1}v_1&  \partial_{x_1} v_2+\partial_{x_2} v_1&  \partial_{x_1}v_3\\
\\
 \partial_{x_1} v_2+\partial_{x_2} v_1& 2\partial_{x_2}v_2 &  \partial_{x_2}v_3\\
\\
 \partial_{x_1}v_3 &  \partial_{x_2}v_3 & 0
\end{array}\right),\quad \partial_{y_3}\left[v\right]=\frac{1}{2}\left(\begin{array}{ccc}
0& 0&  \partial_{y_3}v_1\\
\\
0&0& \partial_{y_3}v_2\\
\\
 \partial_{y_3}v_1 & \partial_{y_3}v_2 & 2\partial_{y_3}v_3
\end{array}\right),
\end{equation}
  and, associated to the change of variables (\ref{dilatacion}), we introduce the operators $\mathbb{D}_{h_\varepsilon}$, $D_{h_\varepsilon}$, ${\rm div}_{h_\varepsilon}$ and $\nabla_{h_\varepsilon}$, by
\begin{equation*}
\mathbb{D}_{h_\varepsilon}\left[v\right]=\mathbb{D}_{x'}[v]+h_\varepsilon^{-1}\partial_{y_3}[v],
\end{equation*}
\begin{equation*}
(D_{h_\varepsilon}v)_{i,j}=\partial_{x_j}v_i\text{\ for \ }i=1,2,3,\ j=1,2,\quad 
(D_{h_\varepsilon}v)_{i,3}=h_\varepsilon^{-1}\partial_{y_3}v_i\text{\ for \ }i=1,2,3,
\end{equation*}
\begin{equation*}
{\rm div}_{h_\varepsilon}v={\rm div}_{x^{\prime}}v^{\prime}+ h_\varepsilon^{-1}\partial_{y_3}v_3,\quad
\nabla_{h_\varepsilon}w=(\nabla_{x^{\prime}}w, h_\varepsilon^{-1}\partial_{y_3}w)^t.
\end{equation*}
Using the transformation (\ref{dilatacion}), then the system (\ref{1}) can be rewritten as follows
\begin{equation}
\left\{
\begin{array}
[c]{r@{\;}c@{\;}ll}%
\medskip
\displaystyle -\mu\,{\rm div}_{h_\varepsilon} \left(  \left\vert \mathbb{D}_{h_\varepsilon}\left[\widetilde {u}_{\varepsilon }\right] \right\vert^{r-2}\mathbb{D}_{h_\varepsilon}\left[\widetilde {u}_{\varepsilon }\right] \right)+ \nabla_{h_\varepsilon} \widetilde {p}_{\varepsilon } &
= f&
 \text{\ in \ }\widetilde{\Omega}_{\varepsilon },\\
\medskip
{\rm div}_{h_\varepsilon} \widetilde {u}_{\varepsilon } & = 0&  \text{\ in \ }\widetilde{\Omega}_{\varepsilon },\\
\medskip
\widetilde {u}_{\varepsilon }& = 0& \text{\ on \ } \partial \widetilde{\Omega}_{\varepsilon }.
\end{array}
\right. \label{2}%
\end{equation}
The objetive of this paper is then  to describe the asymptotic behavior of the dilated sequence of solutions $(\widetilde {u}_{\varepsilon }$, $\widetilde {p}_{\varepsilon })$ of system (\ref{2}), which is described  by the theorem below.  The problem is that the sequence $(\widetilde {u}_{\varepsilon }$, $\widetilde {p}_{\varepsilon })\in W^{1,r}_0(\widetilde\Omega_{\varepsilon })^3\times L^{r'}(\widetilde\Omega_{\varepsilon })/\mathbb{R}$  is not defined in a fixed domain independent of $\varepsilon$, but is defined in $\widetilde\Omega_{\varepsilon }$, which also depends on $\varepsilon$. As usual,  to pass to the limit when $\varepsilon,  \delta_\varepsilon$  and $h_\varepsilon$ tend to zero, we use convergences in fixed Sobolev spaces (defined in $\Omega$) after obtaining optimal {\it a priori} estimates.  To do this, we need to define previously  an extension of  $(\widetilde {u}_{\varepsilon }$, $\widetilde {p}_{\varepsilon })$ to the whole domain $\Omega$. In this sense, we consider the zero extensions to the whole $\Omega$ for both velocity and pressure, which coincide with the original functions in $\widetilde\Omega_{\varepsilon }$.  For simplicity the extensions will be denoted by the same symbol.  Moreover, as recalled in the introduction, we restrict the main result to shear thinning power law fluids, i.e. for  $1<r<2$.

\begin{theorem}[Main theorem]\label{MainTheorem}
Suppose $1<r<2$ and let $\sigma_\varepsilon$ be given by (\ref{sigma}) satisfying (\ref{sigma0}). 
Depending on the values of $\sigma_\varepsilon$ and $h_\varepsilon$, we have the following cases:
\begin{itemize}
\item[(i)] In the case $\sigma_\varepsilon\approx h_\varepsilon$, with $\sigma_\varepsilon/ h_\varepsilon\to \lambda\in (0,+\infty)$, then there exist $  u\in W^{1,r}(0,1;L^r(\omega)^3)$, with $ u=0$ on $y_3=\{0,1\}$ and $u_3\equiv 0$, and $p\in L^{r'}(\Omega)/\mathbb{R}$ independent of $y_3$, such that the extension $(\widetilde u_{\varepsilon},\widetilde p_{\varepsilon})$ of the solution of (\ref{2}) satisfies the following convergences
$$\sigma_{\varepsilon}^{-{r\over r-1}}\widetilde u_{\varepsilon}\rightharpoonup  u\quad\hbox{weakly in }W^{1,r}(0,1;L^r(\omega)^3),\quad \widetilde p_{\varepsilon}\to p\quad\hbox{strongly in }L^{r'}(\Omega),$$
where $(u, p)$, with $u_3\equiv 0$, is the unique solution of the lower-dimensional nonlinear Brinkman-type  problem

\begin{equation}\label{thm:system}
\left\{\begin{array}{rl}
\medskip
\displaystyle - \mu\lambda^r 2^{-{r\over 2}}\,\partial_{y_3}\left(|\partial_{y_3}  u'|^{r-2}\partial_{y_3}  u'\right)+ \mu\, \mathcal{G}( u')+\nabla_{x'}p(x')=f'(x')&\hbox{in }\Omega,\\
\medskip
  u'=0&\hbox{on }y_3=\{0,1\},\\
\medskip
\displaystyle {\rm div}_{x'}\left(\int_0^1  u'\,dy_3\right)=0&\hbox{in }\omega,\\
\medskip
\displaystyle \left(\int_0^1  u'\,dy_3\right)\cdot n'=0&\hbox{on }\partial\omega.
\end{array}\right.
\end{equation}
Here, the drag force function $\mathcal{G}:\mathbb{R}^2\to \mathbb{R}^2$ is defined by
\begin{equation}\label{DragForceG}
\mathcal{G}(\zeta')\cdot \tau' =\int_{\mathbb{R}^2\setminus Y'_s}|\mathbb{D}_{z'}[w^{\zeta'}]|^{r-2}\mathbb{D}_{z'}[w^{\zeta'}]:\mathbb{D}_{z'}[w^{\tau'}]\,dz',\quad \forall\, \tau',\zeta'\in\mathbb{R}^2,
\end{equation}
where $w^{\xi'}$, $\xi'=\{\tau',\zeta'\}$, is the unique solution to the auxiliary exterior problem
\begin{equation}\label{LocalProblems}
\left\{\begin{array}{rl}
\medskip
\displaystyle-{\rm div}_{z'}\left(|\mathbb{D}_{z'}[w^{\xi'}]|^{r-2}\mathbb{D}_{z'}[w^{\xi'}]\right)+ \nabla_{z'}\pi^{\xi'}=0 & \hbox{in }\mathbb{R}^2\setminus Y'_s,\\
\medskip
\displaystyle
{\rm div}_{z'} w^{\xi'}=0 & \hbox{in }\mathbb{R}^2\setminus Y'_s,\\
\medskip
\displaystyle
w^{\xi'}=0& \hbox{on }\partial Y'_s,\\
\medskip
\displaystyle
\lim_{|z'|\to \infty} w^{\xi'}=\xi',&\\
\medskip
(w^{\xi'},\pi^{\xi'})\in D^{1,r}(\mathbb{R}^2\setminus Y'_s)^2\times L^{r'}(\mathbb{R}^2\setminus Y'_s)/\mathbb{R}.&
\end{array}\right.
\end{equation}
The space $D^{1,r}$  denotes the homogeneous Sobolev space given by $D^{1,r}(\mathcal{O})=\left\{v\in L^1_{loc}(\mathcal{O})\,:\, Dv\in L^r(\mathcal{O})\right\}$.

\item[(ii)] In the case $\sigma_\varepsilon\ll h_\varepsilon$, then there exist   $u\in  L^r(\Omega)^3$, with $u_3\equiv 0$, and $p\in L^{r'}(\Omega)/\mathbb{R}$ independent of $y_3$, such that the extension $(\widetilde u_{\varepsilon},\widetilde p_{\varepsilon })$ of the solution of (\ref{2}) satisfies the following convergences
$$\sigma_{\varepsilon}^{-{r\over r-1}}\widetilde u_{\varepsilon}\rightharpoonup  u\quad\hbox{weakly in }L^r(\Omega)^3,\quad \widetilde p_{\varepsilon }\to p\quad\hbox{strongly in }L^{r'}(\Omega).$$
Moreover,   the average velocity $U_{av}(x')=\int_0^1u\,dy_3$ is given by
$$U'_{av}(x') ={1\over \mu^{1\over r-1}}\mathcal{G}^{-1}\left(f'(x')-\nabla_{x'}p(x')\right),\quad \quad U_{av,3}\equiv 0\quad \hbox{in }\omega,$$
where $p\in W^{1,r'}(\omega)\cap (L^{r'}(\Omega)/\mathbb{R})$ is the unique solution of the lower-dimensional nonlinear Darcy problem 
\begin{equation}\label{thm:system2}
\left\{\begin{array}{rl}
\medskip
\displaystyle {\rm div}_{x'} U'_{av}=0&\hbox{in }\omega,\\
\medskip
\displaystyle U'_{av}\cdot n'=0&\hbox{on }\partial\omega.
\end{array}\right.
\end{equation}
Here, $\mathcal{G}^{-1}:\mathbb{R}^2\to\mathbb{R}^2$ is the inverse of the drag force function $\mathcal{G}$ defined by (\ref{DragForceG}), which is asymptotically defined, up to a suitable rescaling, by the permeability function $\mathcal{U}_{\delta_\varepsilon}:\mathbb{R}^2\to\mathbb{R}^2$,   defined by  (see Remark \ref{propG1} for more details)
\begin{equation}\label{permeabilityU}
\mathcal{U}_{\delta_\varepsilon}(\xi')=\int_{Y^{\prime}_{\delta_\varepsilon f}} w^{\xi'}(z')\,dz'\quad (\hbox{where}\quad  Y^{\prime}_{\delta_\varepsilon f}=Y^{\prime}\setminus \delta_\varepsilon   Y'_s),
\end{equation}
with $w^{\xi'}$ the solution of the purely 2D nonlinear auxiliary cell problem    
\begin{equation}\label{aux_problem}
\left\{\begin{array}{rl}
\displaystyle
\medskip
-{\rm div}_{z'}\left(|\mathbb{D}_{z'}[w^{\xi'}]|^{r-2}\mathbb{D}_{z'}[w^{\xi'}]\right)+\nabla_{z'}\pi^{\xi'}=\xi' &\hbox{in }Y^{\prime}\setminus \delta_\varepsilon   Y'_s,\\
\displaystyle
\medskip
{\rm div}_{z'}w^{\xi'}=0&\hbox{in }Y^{\prime}\setminus \delta_\varepsilon  Y'_s,\\
\displaystyle
\medskip
w^{\xi'}=0&\hbox{in } \delta_\varepsilon   Y'_s,\\
\displaystyle
\medskip(w^{\xi'},\pi^{\xi'})\in D^{1,r}_{\rm per}(Y'_{\delta_\varepsilon f})^2\times L^{r'}(Y'_{\delta_\varepsilon f})\setminus \mathbb{R}.&
\end{array}\right.
\end{equation}
The space $D^{1,r}_{\rm per}(Y'_{\delta_\varepsilon f})$ denotes the space of functions in $D^{1,r}(Y'_{\delta_\varepsilon f})$ which are $Y'$-periodic.
\item[(iii)] In the case $\sigma_\varepsilon\gg h_\varepsilon$, then, there exist $  u=(u',0)\in W^{1,r}(0,1;L^r(\omega)^3)$ with $ u=0$ on $y_3=\{0,1\}$ and $p\in L^{r'}(\Omega)/\mathbb{R}$ independent of $y_3$, such that the extension $(\widetilde u_{\varepsilon},\widetilde p_{\varepsilon})$ of the solution of (\ref{2}) satisfies the following convergences
$$h_{\varepsilon}^{-{r\over r-1}}\widetilde u_{\varepsilon}\rightharpoonup  u\quad\hbox{weakly in }W^{1,r}(0,1;L^r(\omega)^3),\quad \widetilde p_{\varepsilon }\to p\quad\hbox{strongly in }L^{r'}(\Omega).$$
Moreover,    the average velocity $U_{av}(x')=\int_0^1u\,dy_3$ is given by 
\begin{equation}\label{U_Reynolds_super}
\begin{array}{ll}
\displaystyle U'_{av}(x')={1\over  2^{r'\over 2}(r'+1)\mu^{r'-1}}\left|f'(x')-\nabla_{x'}p(x')\right|^{r'-2}(f'(x')-\nabla_{x'}p(x')),\quad U_{av,3}\equiv 0 &\hbox{in  }\omega,
\end{array} 
\end{equation}
where $p\in (L^{r'}(\omega)\setminus \mathbb{R})\cap W^{1,r}(\omega)$ is the unique solution of the lower-dimensional  non-linear Reynolds problem
\begin{equation}\label{Reynolds_super}
\left\{\begin{array}{rl}
\displaystyle
{\rm div}_{x'}U'_{av}=0&\hbox{in  }\omega,\\
\\
U'_{av}\cdot n'=0&\hbox{on  }\partial\omega.
\end{array}\right.
\end{equation}
\end{itemize}

\end{theorem}
\begin{remark}[Properties of $\mathcal{G}$]\label{RemG} We recall the following properties:
\begin{itemize}
\item[--]  For every $\xi'\in\mathbb{R}^2$, the exterior auxiliary problem (\ref{LocalProblems})  has a unique solution $(w^{\xi'},\pi^{\xi'})\in D^{1,r}(\mathbb{R}^2\setminus Y'_s)^2\times L^{r'}(\mathbb{R}^2\setminus Y'_s)/\mathbb{R}$, with $1<r<2$, see \cite[Theorem 3]{Marusic-Paloka3}. For more details concerning the homogeneous Sobolev space $D^{1,r}$ we refer to  \cite[Chapter II.6]{Galdi}.
\item[--]  According to   \cite{MarusicPaloka2}, the drag force function  $\mathcal{G}:\mathbb{R}^2\to \mathbb{R}^2$ is continuous and is strictly monotone, that is
$$\left(\mathcal{G}(\xi')-\mathcal{G}(\tau'),\xi'-\tau'\right)>0,\quad\forall\, \xi',\tau'\in\mathbb{R}^2.$$  
Moreover, it satisfies the homogeneity condition
\begin{equation}\label{homogeneitycon}
\mathcal{G}(\lambda\,\xi')=|\lambda|^{r-2}\lambda\,\mathcal{G}(\xi'),\quad \forall\,(\lambda,\xi')\in \mathbb{R}\times \mathbb{R}^2,
\end{equation}
and there exists $m,M>0$ such that for every $\xi'\in\mathbb{R}^2$ it holds
$$m|\xi'|^{r-1}\leq |\mathcal{G}(\xi')|\leq M|\xi'|^{r-1}.$$

\end{itemize}
\end{remark}

\begin{remark}[Properties of $\mathcal{G}^{-1}$]\label{propG1}
We recall the following properties:
\begin{itemize}

\item[--] For  $\xi'\in \mathbb{R}^2$, the auxiliary problem (\ref{aux_problem}) has a unique solution $(v^{\xi'},\pi^{\xi'})$ in $D^{1,r}_{\rm per}(Y'_{\delta_\varepsilon f})\times L^{r'}(Y'_{\delta_\varepsilon f})\setminus \mathbb{R}$ (recall that $Y^{\prime}_{\delta_\varepsilon f}=Y^{\prime}\setminus \delta_\varepsilon   Y'_s$), see \cite[Remark 3.3]{MarusicPaloka2}. 


\item[--] The permeability function $\mathcal{U}_{\delta_\varepsilon}:\mathbb{R}^2 \to \mathbb{R}^2$, defined by (\ref{permeabilityU}) through the solution of the cell problems (\ref{aux_problem}) posed in $Y^{\prime}\setminus \delta_\varepsilon   Y'_s$, appears after the homogenization process of problem (\ref{1}) in $\Omega_{\varepsilon}$ by passing to the limit in $\varepsilon$, but keeping $\delta_\varepsilon$ fixed. In that case, the homogenization problem fits into the case in which the cylinder arrangement period is smaller than the thickness of the domain (see the description of the  ``homogeneously thin porous media" in the introduction, which agrees in the case of $\delta_\varepsilon=1$).  Then, it can be deduced the following lower dimensional Darcy law 
\begin{equation}\label{thm:system2-remark}
\left\{\begin{array}{rl}
\medskip
\displaystyle    U_{av,\delta_\varepsilon}'(x')={1\over \mu}\mathcal{U}_{\delta_\varepsilon}(f'(x')-\nabla_{x'}p(x')),\quad  U_{av,3, \delta_\varepsilon}\equiv 0&\hbox{in }\omega,\\
\medskip
\displaystyle {\rm div}_{x'}   U'_{av,\delta_\varepsilon} =0&\hbox{in }\omega,\\
\medskip
\displaystyle  U_{av,\delta_\varepsilon}' \cdot n'=0&\hbox{on }\partial\omega.
\end{array}\right.
\end{equation}

\item[--] We recall that $\mathcal{U}_{\delta_\varepsilon}$ satisfies the following properties for every $\xi',\tau' \in\mathbb{R}^2$ (see \cite{Bourgeat2} for more details)
\begin{equation}\label{Uprop}
\mathcal{U}_{\delta_\varepsilon}(\lambda\xi')=\lambda^{r'-1}\mathcal{U}_{\delta_\varepsilon}(\xi')\quad \forall\,\lambda>0,
\end{equation}
\begin{equation}\label{Uprop2}
m|\xi'|^{r'-1}\leq |\mathcal{U}_{\delta_\varepsilon}(\xi')|\leq M|\xi'|^{r'-1}\quad\hbox{where}\quad m=\inf_{|\xi'|=1}|\mathcal{U}_{\delta_\varepsilon}(\xi')|,\quad M=\max_{|\xi'|=1}|\mathcal{U}_{\delta_\varepsilon}(\xi')|,
\end{equation}
\begin{equation}\label{Uprop3}
\left (\mathcal{U}_{\delta_\varepsilon}(\xi')-\mathcal{U}_{\delta_\varepsilon}(\tau'), \xi'-\tau'\right)>0,
\end{equation}
  and so, for every $\delta_\varepsilon>0$,  the problem (\ref{thm:system2-remark}) has a unique  solution $p\in W^{1,r'}(\omega)/\mathbb{R}.$

\item[--] Due to its strict monotonicity, function $\mathcal{G}$ is invertible, i.e.  $\mathcal{G}^{-1}(\xi')$ is well defined for $\xi'\in \mathbb{R}^2$, see \cite[Remark 4.6]{MarusicPaloka2}.
\item[--] The main result in \cite[Theorem 4.8]{MarusicPaloka2} states that if we consider $(v^{\xi'}, \pi^{\xi'})$ the solution of problem (\ref{aux_problem}), and consider  the following scaling for the solutions
$$w^{\xi'}_{\delta_\varepsilon}(y')=\delta_\varepsilon^{2-r\over r-1}w^{\xi'}(\delta_\varepsilon y'),\quad \pi^{\xi'}_{\delta_\varepsilon}(y')=\delta_\varepsilon \pi^{\xi'}(\delta_\varepsilon y'),\quad\hbox{a.e.  }y'\in \delta_\varepsilon^{-1}Y'\setminus  Y'_s,$$
which now satisfies 
\begin{equation}\label{aux_problem2}
\left\{\begin{array}{rl}
\displaystyle
\medskip
-{\rm div}_{y'}\left(|\mathbb{D}_{y'}[w^{\xi'}_{\delta_\varepsilon}]|^{r-2}\mathbb{D}_{y'}[w^{\xi'}_{\delta_\varepsilon}]\right)+\nabla_{y'}\pi^{\xi'}_{\delta_\varepsilon}=\delta_\varepsilon^2 \xi' &\hbox{in }\delta_\varepsilon^{-1} Y^{\prime}\setminus     Y'_s,\\
\displaystyle   
\medskip
{\rm div}_{y'}w^{\xi'}_{\delta_\varepsilon}=0&\hbox{in }\delta_\varepsilon^{-1}Y^{\prime}\setminus   Y'_s,\\
\displaystyle
\medskip
w^{\xi'}_{\delta_\varepsilon}=0&\hbox{in }   Y'_s,\\
\displaystyle
\medskip (w^{\xi'}_{\delta_\varepsilon},\pi^{\xi'}_{\delta_\varepsilon})& \delta_\varepsilon^{-1}Y'$-$\hbox{periodic},
\end{array}\right.
\end{equation}
then, the rescaled solutions $(w^{\xi'}_{\delta_\varepsilon},\pi^{\xi'}_{\delta_\varepsilon})$ converge   to a specific solution $(w_{\mathcal{G}^{-1}(\xi')}, \pi_{\mathcal{G}^{-1}(\xi')})$ of the exterior problem (\ref{LocalProblems}) as $\delta_\varepsilon\to 0$:
\begin{equation}\label{conv_sol_aux}
\begin{array}{l}\medskip
\displaystyle
w_{\delta_\varepsilon}^{\xi'}\rightharpoonup w_{\mathcal{G}^{-1}(\xi')}\quad\hbox{ weakly in }L^{r}_{\rm loc}(\mathbb{R}^2\setminus Y'_s)^2,\\
\medskip
\displaystyle
\pi_{\delta_\varepsilon}^{\xi'}\rightharpoonup \pi_{\mathcal{G}^{-1}(\xi')}\quad\hbox{ weakly in }L^{r'}(\mathbb{R}^2\setminus Y'_s)\setminus \mathbb{R},\\
\medskip
\displaystyle
\mathbb{D}_{z'}[w_{\delta_\varepsilon}^{\xi'}] \to \mathbb{D}_{z'}[w_{\mathcal{G}^{-1}(\xi')}]\quad\hbox{ strongly in }L^{r}(\mathbb{R}^2\setminus Y'_s).
\end{array}
\end{equation}
As consequence, it holds the continuity of the permeability function in the low-volume-fraction limit
\begin{equation}\label{lowvolume}\lim_{\delta_\varepsilon\to 0}\delta_\varepsilon^{{2-r\over r-1}}\mathcal{U}_{\delta_\varepsilon}(\xi')=\lim_{\delta_\varepsilon\to 0}\int_{\delta_\varepsilon^{-1} Y'\setminus Y'_s}\delta_\varepsilon^2 w_{\delta_\varepsilon}^{\xi'}(y')\,dy'=\lim_{\delta_\varepsilon\to 0}{1\over |\delta_\varepsilon^{-1}Y'|}\int_{\delta_\varepsilon^{-1} Y'}w_{\delta_\varepsilon}^{\xi'}(y')\,dy'=\mathcal{G}^{-1}(\xi'), 
\end{equation}
for every $\xi'\in\mathbb{R}^2$, where $\mathcal{U}_{\delta_\varepsilon}$  is the permeability function  defined by  (\ref{permeabilityU})
and $\mathcal{G}^{-1}$  is the inverse of the drag force function $\mathcal{G}$, which let us deduce that problem (\ref{thm:system2}) has a unique solution $p\in W^{1,r'}(\omega)/\mathbb{R}$.

\end{itemize}
\end{remark}

\section{Proofs}
In this section, we develop the proof of the main theorem with the following structure. In Subsection \ref{S2} we give the sharp Poincar\'e-Korn inequality and derive the {\it a priori} estimates for velocity. In Subsection \ref{S2p}, we give the inverse of the divergence operator in a thin porous medium and derive the {\it a priori} estimates for pressure. We recall the unfolding method in domains perforated with cylinders of small diameter in Subsection \ref{sec:unf}, which is necessary to capture the influence of the microstructure in the asymptotic behavior. In Subsection \ref{sec:comp} we give some compactness results, which will let us prove the main theorem, whose proof will be given in Subsection \ref{sec:proofthm}.

\subsection{{\it A priori} estimates for velocity}\label{S2}
As said previously, in this subsection we derive  {\it a priori} estimates of the solution $(u_{\varepsilon}, p_{\varepsilon})$ of system (\ref{1}). To do this, we first need a technical, which will let us identify the critical size.  
\begin{lemma}[The Poincar\'e and Korn inequalities]\label{Poincare0}
Suppose  $1<r<2$ and let $\sigma_\varepsilon$ be given by  (\ref{sigma}) satisfying (\ref{sigma0}). Then,  the following Poincar\'e inequalites hold with $C$ independent of $\varepsilon$:
\begin{itemize}
\item[--] If $\sigma_\varepsilon\ll h_\varepsilon$ or $\sigma_\varepsilon\approx h_\varepsilon$ with $\sigma_\varepsilon / h_\varepsilon\to \lambda\in (0,+\infty)$, then
\begin{eqnarray}
&\left\Vert \varphi\right\Vert_{L^r(\Omega_{\varepsilon})^3}\leq C\sigma_\varepsilon\left\Vert D\varphi\right\Vert_{L^r(\Omega_{\varepsilon})^{{3\times3}}},& \forall\,\varphi\in W_0^{1,r}(\Omega_{\varepsilon})^3.\label{Poincarethin}
\end{eqnarray}
\item[--]  If $\sigma_\varepsilon\gg h_\varepsilon$, then
\begin{eqnarray}
&\left\Vert \varphi\right\Vert_{L^r(\Omega_{\varepsilon})^3}\leq Ch_\varepsilon\left\Vert D\varphi\right\Vert_{L^r(\Omega_{\varepsilon})^{{3\times3}}},& \forall\,\varphi\in W_0^{1,r}(\Omega_{\varepsilon})^3.\label{Poincarethin_h}
\end{eqnarray}
\end{itemize}
Moreover, in every case, the following Korn inequality holds with $C$ independent of $\varepsilon$:
\begin{equation}\label{Korn}
\left\Vert D \varphi\right\Vert_{L^r(\Omega_{\varepsilon})^{3 \times 3}}\leq C\left\Vert \mathbb{D}[\varphi]\right\Vert_{L^r(\Omega_{\varepsilon})^{3\times 3}},\quad \forall\,\varphi\in W_0^{1,r}(\Omega_{\varepsilon})^3.\end{equation}
\end{lemma}
\begin{proof} First, the classical Korn inequality in a porous medium implies that (\ref{Korn}) holds. Next, to prove the Poincar\'e inequalities, we divide the proof in different cases depending on the values of $\sigma_\varepsilon$ and $h_\varepsilon$:
\begin{itemize}
\item[--] If $\sigma_\varepsilon\ll h_\varepsilon$ or $\sigma_\varepsilon\approx h_\varepsilon$, applying \cite[Lemma 4.1]{MarusicPaloka} to the bottom domain $\omega_{\varepsilon}$, we have the Poincar\'e inequality
$$\left\Vert \varphi\right\Vert^r_{L^r(\omega_{\varepsilon })^3}\leq C\sigma_\varepsilon^r\left\Vert D_{x'}\varphi\right\Vert^r_{L^r(\omega_{\varepsilon })^{3\times 2}}\leq C\sigma_\varepsilon^r\left\Vert D\varphi\right\Vert^r_{L^r(\omega_{\varepsilon })^{3\times 3}}, $$
with  $\sigma_\varepsilon$  given by (\ref{sigma}) and $C$ independent of $\varepsilon$. 
Integrating previous estimates with respect to $x_3$ between $0$ and $h_\varepsilon$, we obtain the desired estimate (\ref{Poincarethin}).

\item[--] If  $\sigma_\varepsilon\gg h_\varepsilon$, we are able to obtain a more optimal estimate. For every $\varphi(z',x_3)\in W^{1,r}(Y^{\prime}_{\delta_\varepsilon f}\times (0,h_\varepsilon))^3$, $1<r<+\infty$, with $\varphi=0$ on $\partial (Y^{\prime}_{\delta_\varepsilon s}\times (0,h_\varepsilon))$, the Friedrichs inequality in $Y^{\prime}_{\delta_\varepsilon f}\times (0,h_\varepsilon)$ states that 
\begin{equation}\label{Cell_poincare}\int_{Y^{\prime}_{\delta_\varepsilon f}\times (0,h_\varepsilon)}|\varphi|^rdz'dx_3\leq Ch_\varepsilon^r\int_{Y^{\prime}_{\delta_\varepsilon f}\times (0,h_\varepsilon)}|\partial_{x_3}\varphi|^rdz'dx_3,
\end{equation}
where the constant $C$ is independent of $\varepsilon$. For every $k'\in \mathcal{K}_{\varepsilon}$, we apply the change of variables
$$k'+z'={x'\over \varepsilon},\quad dz'={dx'\over \varepsilon^2},$$
which rescales (\ref{Cell_poincare}) from $Y^{\prime}_{\delta_\varepsilon f}\times (0,h_\varepsilon)$ to $Y^{\prime}_{\delta_\varepsilon f_{k^{\prime}},{\varepsilon}}\times (0,h_\varepsilon)$. Thus,  for every $\varphi(x)\in W^{1,r}(Y^{\prime}_{\delta_\varepsilon f_{k^{\prime}},{\varepsilon}}\times (0,h_\varepsilon))^3$ with $\varphi=0$ in $\partial(Y^{\prime}_{\delta_\varepsilon f_{k^{\prime}},{\varepsilon}}\times (0,h_\varepsilon))$, it holds
\begin{equation}\label{Cell_poincare2}\int_{Y^{\prime}_{\delta_\varepsilon f_{k^{\prime}},{\varepsilon}}\times (0,h_\varepsilon)}|\varphi|^rdx\leq Ch_\varepsilon^r\int_{Y^{\prime}_{\delta_\varepsilon f_{k^{\prime}},{\varepsilon}}\times (0,h_\varepsilon)}|\partial_{x_3}\varphi|^rdx\leq Ch_\varepsilon^r\int_{Y^{\prime}_{\delta_\varepsilon f_{k^{\prime}},{\varepsilon}}\times (0,h_\varepsilon)}|D_x\varphi|^rdx,
\end{equation}
with the same constant $C$ as in (\ref{Cell_poincare2}). Summing the inequalities (\ref{Cell_poincare2}) it gives the estimate (\ref{Poincarethin_h}).
\end{itemize}
\end{proof}

\begin{remark}[The Poincaré$-$Korn inequality in $\Omega_{\varepsilon}$]\label{pokorinq} 
Taking into account the results given in Lemma \ref{Poincare0}, for every function $\varphi\in W_0^{1,r}(\Omega_{\varepsilon})^3$, $1<r<2$, we have the following estimates for $C>0$ independent of $\varepsilon$:
\begin{itemize}
\item[--] If $\sigma_\varepsilon\ll h_\varepsilon$  or $\sigma_\varepsilon\approx h_\varepsilon$ with $\sigma_\varepsilon / h_\varepsilon\to \lambda\in (0,+\infty)$,  then
\begin{equation}\label{PoincareKornthin}
\left\Vert \varphi\right\Vert_{L^r(\Omega_{\varepsilon})^3}\leq C\sigma_\varepsilon\left\Vert \mathbb{D}[\varphi]\right\Vert_{L^r(\Omega_{\varepsilon})^{{3\times3}}},\quad \forall\,\varphi\in W_0^{1,r}(\Omega_{\varepsilon})^3.\ 
\end{equation}
\item[--] If $\sigma_\varepsilon\gg h_\varepsilon$, then
\begin{equation}\label{PoincareKornthin2}
\left\Vert \varphi\right\Vert_{L^r(\Omega_{\varepsilon})^3}\leq Ch_\varepsilon\left\Vert \mathbb{D}[\varphi]\right\Vert_{L^r(\Omega_{\varepsilon})^{{3\times3}}},\quad \forall\,\varphi\in W_0^{1,r}(\Omega_{\varepsilon})^3.
\end{equation}
\end{itemize}
\end{remark}

Now, we are ready to derive the {\it a priori} estimates of the solution $(u_{\varepsilon}, p_{\varepsilon})$ of problem (\ref{1}) in $\Omega_{\varepsilon}$.

\begin{lemma}[Estimates for velocity  in $\Omega_{\varepsilon}$]\label{Lemma_a3}
Suppose $1<r<2$ and let $\sigma_\varepsilon$ be given by (\ref{sigma}) satisfying  (\ref{sigma0}). Then, the solution ${u}_{\varepsilon} \in W_0^{1,r}( {\Omega}_{\varepsilon})^3$ of the problem (\ref{1}) satisfies the following estimates, depending on the values of $\sigma_\varepsilon$ and $h_\varepsilon$, for $C>0$ independent of $\varepsilon$:
\begin{itemize}
\item[a)] If $\sigma_\varepsilon\ll h_\varepsilon$ or $\sigma_\varepsilon\approx h_\varepsilon$ with $\sigma_\varepsilon/ h_\varepsilon\to \lambda\in (0,+\infty)$, then
\begin{equation}\label{a1_thin}
\left\Vert  {u}_{\varepsilon}\right\Vert_{L^r({\Omega}_{\varepsilon})^3}\leq Ch_\varepsilon^{1\over r}\sigma_\varepsilon^{{r\over r-1}},\quad \left\Vert \mathbb{D}\left[{u}_{\varepsilon}\right]\right\Vert_{L^r( {\Omega}_{\varepsilon})^{{3\times3}}}\leq Ch_\varepsilon^{1\over r}\sigma_\varepsilon^{{1\over r-1}},
\end{equation}
\begin{equation}\label{Derivada_U_thin}
\left\Vert D {u}_{\varepsilon}\right\Vert_{L^r( {\Omega}_{\varepsilon})^{{3\times3}}}\leq Ch_\varepsilon^{1\over r}\sigma_\varepsilon^{{1\over r-1}}.
\end{equation}

\item[b)] If $\sigma_\varepsilon\gg h_\varepsilon$, then
\begin{equation}\label{a1_thin3}
\left\Vert  {u}_{\varepsilon}\right\Vert_{L^r( {\Omega}_{\varepsilon})^3}\leq Ch_\varepsilon^{{2r-1\over r(r-1)}+1},\quad \left\Vert \mathbb{D}\left[ {u}_{\varepsilon}\right]\right\Vert_{L^r( {\Omega}_{\varepsilon})^{{3\times3}}}\leq Ch_\varepsilon^{{2r-1\over r(r-1)}},
\end{equation}
\begin{equation}\label{Derivada_U_thin3}
\left\Vert D {u}_{\varepsilon}\right\Vert_{L^r( {\Omega}_{\varepsilon})^{{3\times3}}}\leq Ch_\varepsilon^{{2r-1\over r(r-1)}}.
\end{equation}
\end{itemize}
\end{lemma}
\begin{proof}
 Multiplying the first equation of (\ref{1})  by ${u}_{\varepsilon}$, integrating over ${\Omega}_{\varepsilon}$ and taking into account that ${\rm div}\,u_\varepsilon=0$ in $\Omega_\varepsilon$, then we have
\begin{eqnarray}\label{a1_new}
\mu\int_{ {\Omega}_{\varepsilon }}\left\vert \mathbb{D} \left[ {u}_{\varepsilon }\right]\right\vert^{r-2}\mathbb{D} \left[ {u}_{\varepsilon }\right]:\mathbb{D} \left[ {u}_{\varepsilon }\right]dx=\int_{ {\Omega}_{\varepsilon }}f\cdot {u}_{\varepsilon }\,dx.
\end{eqnarray}
By means of  H${\rm \ddot{o}}$lder's inequality, the assumption of $f$ given in (\ref{sourcef}) and the relation $1/r'=(r-1)/r$, we deduce
\begin{eqnarray*}
\int_{ {\Omega}_{\varepsilon }}f\cdot {u}_{\varepsilon }\,dx^{\prime}dy_3\leq Ch_\varepsilon^{r-1\over r} \left\Vert {u}_{\varepsilon} \right\Vert_{L^r( {\Omega}_{\varepsilon })^3},
\end{eqnarray*}
and so we deduce that, from (\ref{a1_new}), it holds
\begin{equation}\label{a2}
\left\Vert \mathbb{D}\left[ {u}_{\varepsilon }\right]\right\Vert_{L^r( {\Omega}_{\varepsilon })^{{3\times3}}}^r\leq C h_\varepsilon^{r-1\over r} \left\Vert {u}_{\varepsilon } \right\Vert_{L^r( {\Omega}_{\varepsilon })^3}.
\end{equation}
Depending on the values of $\sigma_\varepsilon$ and $h_\varepsilon$, we have the following:
\begin{itemize}
\item[--] If $\sigma_\varepsilon\ll h_\varepsilon$ or $\sigma_\varepsilon\approx h_\varepsilon$, we use (\ref{PoincareKornthin}) in (\ref{a2}) to deduce
$$
\left\Vert \mathbb{D}\left[ {u}_{\varepsilon }\right]\right\Vert_{L^r( {\Omega}_{\varepsilon })^{{3\times3}}}^r\leq C h_\varepsilon^{r-1\over r} \sigma_\varepsilon \left\Vert \mathbb{D}\left[ {u}_{\varepsilon }\right]\right\Vert_{L^r( {\Omega}_{\varepsilon })^{{3\times3}}},
$$
and then, we get the second estimate in (\ref{a1_thin}). As consequence,
  from  Korn's inequality (\ref{Korn}), we obtain (\ref{Derivada_U_thin}). Now, using (\ref{PoincareKornthin}) and the second estimate in (\ref{a1_thin}), we deduce the first estimate in (\ref{a1_thin}).

\item[--] If $\sigma_\varepsilon\gg h_\varepsilon$, we use (\ref{PoincareKornthin2}) in (\ref{a2}) to deduce
$$
\left\Vert \mathbb{D}\left[ {u}_{\varepsilon }\right]\right\Vert_{L^r( {\Omega}_{\varepsilon })^{{3\times3}}}^r\leq C h_\varepsilon^{2r-1\over r}   \left\Vert \mathbb{D}\left[ {u}_{\varepsilon }\right]\right\Vert_{L^r( {\Omega}_{\varepsilon })^{{3\times3}}},
$$
and then, we obtain the second estimate in (\ref{a1_thin3}). The rest of velocity estimates are obtained in a similar way    to the previous case.
\end{itemize}

\end{proof}

Below, we get the estimates for dilated velocity in $\widetilde \Omega_{\varepsilon}$.
\begin{lemma}[Estimates for dilated velocity  in $\widetilde\Omega_{\varepsilon}$]\label{estOmega}
Suppose $1<r<2$ and let $\sigma_\varepsilon$ be given by (\ref{sigma}) satisfying (\ref{sigma0}). Then, 
  the solution $\widetilde {u}_{\varepsilon } \in W_0^{1,r}(\widetilde{\Omega}_{\varepsilon })^3$ of the dilated problem (\ref{2}) satisfies the following estimates, depending on the values of $\sigma_\varepsilon$ and $h_\varepsilon$, for $C>0$ independent of $\varepsilon$:

\begin{itemize}
\item[--] If $\sigma_\varepsilon\ll h_\varepsilon$, then
\begin{equation}\label{a1_thin_tilde}
\left\Vert  \widetilde {u}_{\varepsilon }\right\Vert_{L^r(\widetilde {\Omega}_{\varepsilon })^3}\leq C \sigma_\varepsilon^{{r\over r-1}},\quad \left\Vert \mathbb{D}_{x'}\left[\widetilde {u}_{\varepsilon }\right]\right\Vert_{L^r( \widetilde {\Omega}_{\varepsilon })^{{3\times 2}}}\leq C\sigma_\varepsilon^{{1\over r-1}},\quad \left\Vert \partial_{y_3}\left[\widetilde {u}_{\varepsilon }\right]\right\Vert_{L^r( \widetilde {\Omega}_{\varepsilon })^{{3 }}}\leq C\sigma_\varepsilon^{{1\over r-1}}h_\varepsilon,
\end{equation}
\begin{equation}\label{Derivada_U_thin_tilde}
\left\Vert D_{x'} \widetilde {u}_{\varepsilon }\right\Vert_{L^r( \widetilde {\Omega}_{\varepsilon })^{{3\times 2}}}\leq C\sigma_\varepsilon^{{1\over r-1}},\quad \left\Vert \partial_{y_3} \widetilde {u}_{\varepsilon }\right\Vert_{L^r( \widetilde {\Omega}_{\varepsilon })^{{3}}}\leq \sigma_\varepsilon^{{1\over r-1}}h_\varepsilon.
\end{equation}
\item[--] If $\sigma_\varepsilon\approx h_\varepsilon$ with $\sigma_\varepsilon/ h_\varepsilon\to \lambda\in (0,+\infty)$, then
\begin{equation}\label{a1_thin2_tilde}
\left\Vert  \widetilde {u}_{\varepsilon }\right\Vert_{L^r(\widetilde {\Omega}_{\varepsilon })^3}\leq C \sigma_\varepsilon^{r\over r-1},\quad \left\Vert \mathbb{D}_{x'}\left[\widetilde {u}_{\varepsilon }\right]\right\Vert_{L^r( \widetilde {\Omega}_{\varepsilon })^{{3\times2}}}\leq C \sigma_\varepsilon^{1\over r-1},\quad  \left\Vert \partial_{y_3}\left[\widetilde {u}_{\varepsilon }\right]\right\Vert_{L^r( \widetilde {\Omega}_{\varepsilon })^{{3}}}\leq C \sigma_\varepsilon^{r\over r-1},
\end{equation}
\begin{equation}\label{Derivada_U_thin2_tilde}
\left\Vert D_{x'} \widetilde {u}_{\varepsilon }\right\Vert_{L^r( \widetilde {\Omega}_{\varepsilon })^{{3\times 2}}}\leq C\sigma_\varepsilon^{1\over r-1}, \quad \left\Vert \partial_{y_3} \widetilde {u}_{\varepsilon }\right\Vert_{L^r( \widetilde {\Omega}_{\varepsilon })^{{3}}}\leq \sigma_\varepsilon^{{r\over r-1}}.
\end{equation}
\item[--] If $\sigma_\varepsilon\gg h_\varepsilon$, then
\begin{equation}\label{a1_thin3_tilde}
\left\Vert  \widetilde {u}_{\varepsilon }\right\Vert_{L^r( \widetilde {\Omega}_{\varepsilon })^3}\leq Ch_\varepsilon^{r\over r-1},\quad \left\Vert \mathbb{D}_{x'}\left[ \widetilde {u}_{\varepsilon }\right]\right\Vert_{L^r( \widetilde {\Omega}_{\varepsilon })^{{3\times 2}}}\leq Ch_\varepsilon^{1\over r-1},\quad \left\Vert \partial_{y_3}\left[ \widetilde {u}_{\varepsilon }\right]\right\Vert_{L^r( \widetilde {\Omega}_{\varepsilon })^{{3}}}\leq Ch_\varepsilon^{r\over r-1},
\end{equation}
\begin{equation}\label{Derivada_U_thin3_tilde}
\left\Vert D_{x'}\widetilde  {u}_{\varepsilon }\right\Vert_{L^r( \widetilde {\Omega}_{\varepsilon })^{{3\times 2}}}\leq Ch_\varepsilon^{1\over r-1},\quad \left\Vert \partial_{y_3}\widetilde  {u}_{\varepsilon }\right\Vert_{L^r( \widetilde {\Omega}_{\varepsilon })^{{3}}}\leq Ch_\varepsilon^{r\over r-1}.
\end{equation}
\end{itemize}
\end{lemma}
\begin{proof}
The proof follows from applying  the change of variables (\ref{dilatacion}) to estimates given in Lemma \ref{Lemma_a3}, taking into account that
$$\begin{array}{c}
\medskip
\displaystyle 
\| u_{\varepsilon  }\|_{L^r(\Omega_{\varepsilon  })^{3}}=h_\varepsilon^{{1\over r}}\|\widetilde u_{\varepsilon  }\|_{L^r(\widetilde \Omega_{\varepsilon  })^{3}},\quad \displaystyle  \| \mathbb{D}[ u_{\varepsilon }]\|_{L^r(\Omega_{\varepsilon  })^{3\times 3}}=h_\varepsilon^{{1\over r}}\|\mathbb{D}_{h_\varepsilon}[\widetilde u_{\varepsilon  }]\|_{L^r(\widetilde \Omega_{\varepsilon  })^{3\times 3}},\\
\medskip
\displaystyle
 \| D u_{\varepsilon }\|_{L^r(\Omega_{\varepsilon  })^{3\times 3}}=h_\varepsilon^{{1\over r}}\|D_{h_\varepsilon}\widetilde u_{\varepsilon}\|_{L^r(\widetilde \Omega_{\varepsilon  })^{3\times 3}}.
\end{array}$$
We remark that estimates in the critical case $\sigma_\varepsilon\approx h_\varepsilon$ have been written in terms of $\sigma_\varepsilon$.
 
\end{proof}

\begin{remark}[Extension of $\widetilde {u}_{\varepsilon }$ to the whole domain $\Omega$] 
We extend the velocity $\widetilde {u}_{\varepsilon }$ by zero to $\Omega$ and denote the extension by the same symbol. Obviously, estimates of velocity $\widetilde u_{\varepsilon }$ given in Lemma \ref{estOmega} remain valid and the extension is divergence free too. 
\end{remark}

\subsection{{\it A priori} estimates for pressure}\label{S2p}

In this subsection, we   introduce the inverse of the divergence operator for thin porous medium, which will let us derive optimal estimates for the pressure. The idea is to find, for any $g\in L^r(\Omega_\varepsilon)$, a function $\varphi=\varphi(g)$ such that ${\rm div}\varphi=g$ and the following estimate
$$\|\varphi\|_{L^r(\Omega_{\varepsilon})^{3}}\leq C\|g\|_{L^r(\Omega_{\varepsilon})},\quad 
\|D\varphi\|_{L^r(\Omega_{\varepsilon})^{3\times 3}}\leq C_\varepsilon\|g\|_{L^r(\Omega_{\varepsilon})},
$$
 with the constant $C_\varepsilon$ such that it is the inverse of the Poincaré constant, i.e. $C_\varepsilon=C\sigma_\varepsilon^{-1}$ in the cases $\sigma_\varepsilon\approx h_\varepsilon$ or $ \sigma_\varepsilon\ll h_\varepsilon$, and $C_\varepsilon=Ch_\varepsilon^{-1}$ in the case $\sigma_\varepsilon\gg h_\varepsilon$, with $C>0$ independent of $\varepsilon$.

\begin{lemma}[Inverse of the divergence operator in $\Omega_{\varepsilon}$]\label{lem:Duvjnak}
Suppose $1<r<2$ and let $\sigma_\varepsilon$ be given by (\ref{sigma}) satisfying (\ref{sigma0}).  
Then, there exists a constant $C$ independent of $\varepsilon$ such that for any $g\in L^{r}(\Omega_{\varepsilon})$, there exists $\varphi=\varphi(g)\in W^{1,r}(\Omega_\varepsilon)^3$  with $\varphi=0$ on $\partial Q_\varepsilon$ such that
\begin{equation}\label{lemDuv1}
{\rm div}\,\varphi=g\ \hbox{ in }\Omega_{\varepsilon},
\end{equation}
and, moreover, depending on the values of $\sigma_\varepsilon$ and $h_\varepsilon$, for some constant $C>0$ independent of $\varepsilon$:
\begin{itemize}
\item[--] If $\sigma_\varepsilon\ll h_\varepsilon$ or $\sigma_\varepsilon\approx h_\varepsilon$ with $\sigma_\varepsilon/ h_\varepsilon\to \lambda\in (0,+\infty)$, then
\begin{equation}\label{lemDuv2}
\|\varphi\|_{L^r(\Omega_{\varepsilon})^{3}}\leq C\|g\|_{L^r(\Omega_{\varepsilon})},\quad 
\|D\varphi\|_{L^r(\Omega_{\varepsilon})^{3\times 3}}\leq {C \sigma_\varepsilon^{-1}}\|g\|_{L^r(\Omega_{\varepsilon})}.
\end{equation}

\item[--] If $\sigma_\varepsilon\gg h_\varepsilon$, then
\begin{equation}\label{lemDuv2-2}
\|\varphi\|_{L^r(\Omega_{\varepsilon})^{3}}\leq C\|g\|_{L^r(\Omega_{\varepsilon})},\quad 
\|D\varphi\|_{L^r(\Omega_{\varepsilon})^{3\times 3}}\leq {C h_\varepsilon^{-1}}\|g\|_{L^r(\Omega_{\varepsilon})}.
\end{equation}

\end{itemize}

\end{lemma}
\begin{proof} The proof follows the lines of  \cite[Lemma 3.2]{Anguiano_SG_Lower}.  Thus, given $g\in L^r(\Omega_{\varepsilon})$, it is possible to  define an extension inside the cylinders  called $G\in L^r_0(Q_{\varepsilon})=\left\{q\in L^r(Q_{\varepsilon})\,:\, \int_{Q_{\varepsilon}}q\,dx=0\right\}$ such that
\begin{equation}\label{normDuv1}
\|G\|_{L^r(Q_{\varepsilon})}\leq C\|g\|_{L^r(\Omega_{\varepsilon})}.
\end{equation}
%
%
%
Since $G\in L^r(Q_{\varepsilon })$, it follows from \cite[Lemma 4]{Duvjak}, that there exists $\varphi\in W^{1,r}_0(Q_\varepsilon)^3$ such that it holds the divergence equation
\begin{equation}\label{lemDuv1_proof}
{\rm div}\,\varphi=G\ \hbox{ in }Q_{\varepsilon},
\end{equation}
\begin{equation}\label{lemDuv2_proof}
\|\varphi\|_{L^r(Q_{\varepsilon})^3}\leq C\|G\|_{L^r(Q_{\varepsilon})},
\end{equation}
\begin{equation}\label{lemDuv3_proof}
\|D\varphi\|_{L^r(Q_{\varepsilon})^{3\times 3}}\leq {Ch_\varepsilon^{-1}}\|G\|_{L^r(Q_{\varepsilon})}.\end{equation}
Let us consider $\varphi_{|_{\Omega_{\varepsilon}}}$: it  belongs to $W^{1,r}(\Omega_{\varepsilon})$ with $\varphi=0$ on $\partial Q_{\varepsilon}$. Then, (\ref{lemDuv1}) follows from (\ref{lemDuv1_proof}) and  estimates given in (\ref{lemDuv2}) follow from (\ref{lemDuv2_proof})  and (\ref{lemDuv4_proof}).

We observe that in the cases $\sigma_\varepsilon\ll h_\varepsilon$ or $\sigma_\varepsilon\approx h_\varepsilon$, we need a more accurate estimate. By using relation between the parameters in these cases,  then (\ref{lemDuv3_proof}) implies 
\begin{equation}\label{lemDuv4_proof}\|D\varphi\|_{L^r(Q_{\varepsilon})^{3\times 3}}\leq {C\sigma_\varepsilon^{-1}}\|G\|_{L^r(Q_{\varepsilon})}.
\end{equation}
As before, considering $\varphi_{|_{\Omega_{\varepsilon}}}$, we deduce the result, i.e. problem (\ref{lemDuv1}) and estimates (\ref{lemDuv2-2}).

\end{proof}

\begin{lemma}[Estimates for pressure  in $\Omega_{\varepsilon}$]\label{Lemma_a3p}
Suppose $1<r<2$ and let $\sigma_\varepsilon$ be given by (\ref{sigma}) satisfying  (\ref{sigma0}). Then,  the pressure $p_{\varepsilon}\in L^{r'}(\Omega_{\varepsilon})/\mathbb{R}$, with $r'=r/(r-1)$,   solution of the problem (\ref{1}),  satisfies the following estimate for some constant $C>0$ independent of $\varepsilon$:
\begin{equation}\label{estim_P_thin}
\|p_{\varepsilon}\|_{L^{r'}(\Omega_{\varepsilon})}\leq Ch_\varepsilon^{1\over r'}\,. 
\end{equation}
\end{lemma}
\begin{proof}

The proof is similar to the proof given in \cite[Lemma 3.2]{Anguiano_SG_Lower}, just taking into account that the height of the domain here is $h_\varepsilon$, instead of $\sigma_\varepsilon$.  We develop the case $\sigma_\varepsilon\ll h_\varepsilon$, because the rest of the cases can be developed similarly, just taking into account the corresponding estimates for velocity in each case. 
\\

Assume $\sigma_\varepsilon\ll h_\varepsilon$ and consider  $g=\left|p_{\varepsilon }\right|^{r'-2}p_{\varepsilon }$, which satisfies $g\in L^r(\Omega_{\varepsilon })$ due to $p_{\varepsilon }\in L^{r'}(\Omega_{\varepsilon })$. From Lemma \ref{lem:Duvjnak}, we have that there exists $\varphi=\varphi(g)\in W^{1,r}_0(\Omega_{\varepsilon})^3$  such that
\begin{equation}\label{Duv1_proofestim}
{\rm div}\,\varphi=g\ \hbox{ in }\Omega_{\varepsilon },\quad \|\varphi\|_{W^{1,q}_0(\Omega_{\varepsilon})^{3}}\leq {C \sigma_\varepsilon^{-1}}\|g\|_{L^r(\Omega_{\varepsilon })}.
\end{equation}
Now, multiplying  the first equation of (\ref{1}) by $\varphi\in W^{1,r}_0(\Omega_{\varepsilon })^3$    and integrating over ${\Omega}_{\varepsilon }$, from the second estimate in (\ref{a1_thin}) and (\ref{Poincarethin}), we  deduce
$$\begin{array}{rl}
\displaystyle\medskip
\left|\int_{\Omega_{\varepsilon }}p_{\varepsilon }\,{\rm div}\,\varphi\,dx\right| &\displaystyle 
\leq C\|\mathbb{D}[u_{\varepsilon }]\|_{L^r(\Omega_{\varepsilon })^{3\times 3}}^{r-1}\|D\varphi\|_{L^r(\Omega_{\varepsilon })^{3\times 3}}
+ Ch_\varepsilon^{r-1\over r}\|\varphi\|_{L^r(\Omega_{\varepsilon })^{3}}\\
\medskip
&\displaystyle 
\leq C\|\mathbb{D}[u_{\varepsilon }]\|_{L^r(\Omega_{\varepsilon })^{3\times 3}}^{r-1}\|D\varphi\|_{L^r(\Omega_{\varepsilon })^{3\times 3}}
+ Ch_\varepsilon^{r-1\over r}\sigma_\varepsilon\|D\varphi\|_{L^r(\Omega_{\varepsilon })^{3\times 3}}\\
&\displaystyle \leq C h_\varepsilon^{r-1\over r}\sigma_\varepsilon\|D\varphi\|_{L^r(\Omega_{\varepsilon })^{3\times 3}}.
\end{array}$$
By using   (\ref{Duv1_proofestim}), we get 
\begin{equation}\label{estim_final_p}
\displaystyle
\|p_{\varepsilon }\|_{L^{r'}(\Omega_{\varepsilon })}^{r'} \leq  Ch_\varepsilon^{r-1\over r} \|g\|_{L^r(\Omega_{\varepsilon })},
\end{equation}
and from $\|g\|_{L^r(\Omega_{\varepsilon })}=\|p_{\varepsilon }\|
_{L^{r'}(\Omega_{\varepsilon })}^{r'-1}$,  we obtain
$$\|p_{\varepsilon }\|_{L^{r'}(\Omega_{\varepsilon })}^{r'}\leq Ch_\varepsilon^{r-1\over r} \|p_{\varepsilon }\|
_{L^{r'}(\Omega_{\varepsilon })}^{r'-1}.$$
 Taking into account that $1/r'=(r-1)/r$, this gives (\ref{estim_P_thin}).
 
\end{proof}

\begin{lemma}[Estimates for  dilated pressure  in $\widetilde \Omega_\varepsilon$]\label{estOmegap}
Suppose $1<r<2$ and let $\sigma_\varepsilon$ be given by (\ref{sigma}) satisfying (\ref{sigma0}). Then, 
  the  dilated pressure  $\widetilde p_\varepsilon\in L^{r'}(\widetilde \Omega_\varepsilon)/\mathbb{R}$ satisfies the following estimate for some constant $C>0$ independent of $\varepsilon$:
\begin{equation}\label{estim_P_thin_tilde}
\|\widetilde p_{\varepsilon }\|_{L^{r'}(\widetilde \Omega_\varepsilon)}\leq C.
\end{equation}
\end{lemma}
\begin{proof}
The proof follows from applying  the change of variables (\ref{dilatacion}) to estimate given in Lemma \ref{Lemma_a3p}, taking into account that
$$ 
  \| p_{\varepsilon  }\|_{L^{r'}(\widetilde \Omega_\varepsilon )}=h_\varepsilon^{{1\over r'}}\|\widetilde p_{\varepsilon  }\|_{L^{r'}(\widetilde \Omega_\varepsilon)}.
$$ 
\end{proof}
\begin{remark}[Extension of $\widetilde {p}_{\varepsilon }$ to the whole domain $\Omega$] 
We extend the pressure $\widetilde{p}_{\varepsilon }$ by zero to $\Omega$ and denote the extension by the same symbol. Obviously, estimate of pressure $\widetilde p_{\varepsilon }$ given in Lemma \ref{estOmegap} remains valid.
\end{remark}

%
%
\subsection{Unfolding Method in domains with cylinders of small diameter}\label{sec:unf}
In this section, we first recall the version and some properties of the unfolding method adapted to  domains perforated with cylinders of diameters of size $\varepsilon \delta_\varepsilon$ distributed periodically with period $\varepsilon$ introduced in \cite{Anguiano_SG_Lower}, which is necessary to capture the influence of the microstructure of  $\widetilde{\Omega}_{\varepsilon}$ in the behavior of $(\widetilde {u}_{\varepsilon }, \widetilde p_{\varepsilon })$ by introducing the unfolded functions $(\widehat  {u}_{\varepsilon }, \widehat  p_{\varepsilon })$. Then, by using the estimates of $(\widetilde {u}_{\varepsilon }, \widetilde  p_{\varepsilon })$, we will obtain the estimates for the unfolded functions.

\begin{definition}[Unfolding operator in domains perforated with cylinders of small diameter]
For $\widetilde  \varphi\in L^q(\Omega)$, $1\leq q\leq +\infty$, we define $\widehat  \varphi \in L^q(\omega\times \mathbb{R}^2\times (0,1))$ by
\begin{eqnarray}\label{hat}
\widehat {\varphi}(x^{\prime},z',y_3)=\left\{\begin{array}{ll}
\displaystyle\widetilde  \varphi\left( {\varepsilon}\kappa\left(\frac{x^{\prime}}{{\varepsilon}} \right)+{\varepsilon}\delta_\varepsilon z^{\prime},y_3 \right),&\hbox{if }(x^{\prime},z',y_3)\in \omega\times {1\over \delta_\varepsilon}Y'\times (0,1),\\
 \\
0&\hbox{otherwise,}
 \end{array}\right.
\end{eqnarray}
where the function $\kappa:\mathbb{R}^2\to \mathbb{Z}^2$ is defined  by
\begin{equation*}
\kappa(x^\prime)=k^\prime\iff x'\in Y^{\prime}_{k^{\prime},1}\,,\quad \forall\, k'\in \mathbb{Z}^2.
\end{equation*}
\end{definition}
\begin{remark}\label{remarkCV}We make the following comments:
\begin{itemize}
\item[--] For $\delta_\varepsilon=1$ we are in presence of the adaptation of the unfolding operator for domains with cylinders introduced in \cite[Subsection 4.2]{Anguiano_SuarezGrau}.

\item[--] The function $\kappa$ is well defined up to a set of zero measure in $\mathbb{R}^2$ (the set $\cup_{k^\prime\in \mathbb{Z}^2}\partial Y^{\prime}_{k^{\prime},1}$). Moreover, for every $\varepsilon>0$, we have 
$$\kappa\left(\frac{x^\prime}{\varepsilon}\right)=k^\prime \iff x^\prime\in Y^{\prime}_{k^{\prime},\varepsilon}.$$

\item[--] For $k^{\prime}\in \mathcal{K}_{\varepsilon}$, the restriction of $\widehat {\varphi}$ to $Y^{\prime}_{k^{\prime},{\varepsilon}}\times {1\over \delta_\varepsilon}Y'\times (0,1)$ does not depend on $x^{\prime}$, whereas as a function of $z'$ it is obtained from $\widetilde {\varphi}$ by using the changes of variables  $\delta_\varepsilon z^{\prime}=y'$ and 
\begin{equation}\label{CV}
y^{\prime}=\frac{x^{\prime}-{\varepsilon}k^{\prime}}{{\varepsilon}}, 
\end{equation}
which transform $Y_{k^{\prime},{\varepsilon}}$ into ${1\over \delta_\varepsilon}Y'\times (0,1)$.
\end{itemize}
\end{remark}

\begin{theorem}[Properties of the unfolding operator]\label{estimates_unfolding}
We have the following properties of $\widehat  \varphi$:
\begin{enumerate}
\item Suppose $1\leq q\leq +\infty$. For every $\widetilde  \varphi \in L^q(\Omega)$, it holds
\begin{equation}\label{estimate_varphi1}
\left\Vert \widehat {\varphi}\right\Vert_{L^q(\omega\times \mathbb{R}^2\times (0,1))}\leq \delta_\varepsilon^{-{2\over q}}\left\Vert \widetilde {\varphi}\right\Vert_{L^q(\Omega)}.
\end{equation}
\item Suppose $1\leq q\leq +\infty$. For every $\widetilde  \varphi \in W^{1,q}(\Omega)$, it holds
\begin{equation}\label{estimate_varphi2}
\left\Vert \mathbb{D}_{z^{\prime}}\!\left[\widehat {\varphi}\right]\right\Vert_{L^q(\omega\times {1\over \delta_\varepsilon}Y'\times (0,1))^2}\!\leq\! \sigma_\varepsilon\left\Vert \mathbb{D}_{x^{\prime}}\!\left[\widetilde {\varphi}\right]\right\Vert_{L^q(\Omega)^2},
\end{equation}
\begin{equation}\label{estimate_varphi3}
\left\Vert \partial_{y_3}\!\left[\widehat {\varphi}\right]\right\Vert_{L^q(\omega\times {1\over \delta_\varepsilon}Y'\times (0,1))}\!\leq\! \delta^{-{2\over q}}\left\Vert \partial_{y_3}\!\left[\widetilde {\varphi}\right]\right\Vert_{L^q(\Omega)}.
\end{equation}
\item Suppose $1\leq q < 2$ and let $\mathcal{O}$ be a bounded open set in $\mathbb{R}^2$. For every $\widetilde  \varphi \in W^{1,q}(\Omega)$, it holds
\begin{equation}\label{estimate_varphi2_mean}
\left\Vert \widehat  \varphi-\bar \varphi\right\Vert_{L^q(\Omega; L^{q^*}(\mathbb{R}^2))}\!\leq\! C\sigma_\varepsilon\left\Vert \mathbb{D}_{x^{\prime}}\!\left[\widetilde {\varphi}\right]\right\Vert_{L^q(\Omega)^2 },
\end{equation}
\begin{equation}\label{estimate_varphi2_mean2}
\|\widehat  \varphi\|_{L^q(\omega\times \mathcal{O}\times (0,1))}\leq C|\mathcal{O}|^{1\over 2}\sigma_{\varepsilon}\|\mathbb{D}_{x'}[\widetilde  \varphi]\|_{L^q(\Omega)^2 }+|\mathcal{O}|^{1\over r}\|\widetilde  \varphi\|_{L^q(\Omega)},
\end{equation}
where $q^*={2q\over 2-q}$ be the associated Sobolev exponent, $C$ denotes the Sobolev-Poincar\'e-Wirtinger constant for $W^{1,q}(Y')$ and $\bar\varphi \in L^q(\Omega)$ is the local average defined by
\begin{equation}\label{avoperator}
\bar \varphi(x',y_3)={1\over \varepsilon^{2}}\int_{\varepsilon\kappa\left({x'\over \varepsilon}\right)+\varepsilon Y'}\widetilde  \varphi(\tau',y_3)\,d\tau'=
\displaystyle\delta^2_\varepsilon\int_{{1\over \delta_\varepsilon}Y'}\widehat  \varphi(x', \tau', y_3)\,d\tau', \quad \forall\,\widetilde \varphi\in L^q(\Omega).
\end{equation}
\end{enumerate}
\end{theorem}
\begin{definition}
From extensions of the dilated velocity and pressure $(\widetilde  u_{\varepsilon }, \widetilde  p_{\varepsilon })$, we define the unfolded velocity and pressure $(\widehat  u_{\varepsilon }, \widehat  p_{\varepsilon })$ by using (\ref{hat}). 
\end{definition}
\begin{lemma}[Estimates of the unfolded functions]\label{estCV}
Suppose $1<r<2$, let $\sigma_\varepsilon$ be given by (\ref{sigma}) satisfying (\ref{sigma0}), and let $\mathcal{O}$ be a bounded open set in $\mathbb{R}^2$. Then, the unfolded functions  $(\widehat  u_{\varepsilon }, \widehat  p_{\varepsilon })$ satisfy the following estimates depending on the values of $\sigma_\varepsilon$ and $h_\varepsilon$:
\begin{itemize}
\item[--] If $\sigma_\varepsilon\ll h_\varepsilon$, then
\begin{equation}\label{unfolding3}
\left\Vert \widehat {u}_{\varepsilon }\right\Vert_{L^r(\omega\times \mathcal{O}\times (0,1))^3}\leq C \sigma_\varepsilon^{r\over r-1}\,,
\end{equation}
\begin{equation}\label{unfolding1}
\left\Vert \mathbb{D}_{z^{\prime}}\!\left[\widehat {u}_{\varepsilon }\right]\right\Vert_{L^r(\omega\times {1\over \delta_\varepsilon}Y'\times (0,1))^{3\times2}}\!\leq\! C\sigma_\varepsilon^{r\over r-1},\quad \left\Vert \partial_{y_3}\!\left[\widehat {u}_{\varepsilon }\right]\right\Vert_{L^r(\omega\times {1\over \delta_\varepsilon}Y'\times (0,1))^3}\!\leq\! C\sigma_\varepsilon^{1\over r-1}h_\varepsilon \delta_\varepsilon^{-{2\over r}},
\end{equation}
\begin{equation}\label{unfolding2}
\left\Vert D_{z^{\prime}}\widehat {u}_{\varepsilon }\right\Vert_{L^r(\omega\times {1\over \delta_\varepsilon}Y'\times (0,1))^{{3\times2}}}\!\leq\! C\sigma_\varepsilon^{r\over r-1}, \quad\left\Vert \partial_{y_3}\widehat {u}_{\varepsilon }\right\Vert_{L^r(\omega\times {1\over \delta_\varepsilon}Y'\times (0,1))^3}\!\leq\! C\sigma_\varepsilon^{1\over r-1}h_\varepsilon \delta_\varepsilon^{-{2\over r}},
\end{equation}
\begin{equation}\label{unfolding2_mean}
\left\Vert \widehat  u_{\varepsilon }-\bar u_{\varepsilon }\right\Vert_{L^r(\Omega; L^{r^*}(\mathbb{R}^2)^3)}\!\leq\! C\sigma_\varepsilon^{r\over r-1},
\end{equation}

\item[--] If $\sigma_\varepsilon\approx h_\varepsilon$ with $\sigma_\varepsilon/ h_\varepsilon\to \lambda\in (0,+\infty)$, it holds
\begin{equation}\label{unfolding3_2}
\left\Vert \widehat {u}_{\varepsilon }\right\Vert_{L^r(\omega\times \mathcal{O}\times (0,1))^3}\leq C\sigma_\varepsilon^{r\over r-1}\,,
\end{equation}
\begin{equation}\label{unfolding1_2}
\left\Vert \mathbb{D}_{z^{\prime}}\!\left[\widehat {u}_{\varepsilon }\right]\right\Vert_{L^r(\omega\times {1\over \delta_\varepsilon}Y'\times (0,1))^{3\times2}}\!\leq\! C\sigma_\varepsilon^{r\over r-1},\quad \left\Vert \partial_{y_3}\!\left[\widehat {u}_{\varepsilon }\right]\right\Vert_{L^r(\omega\times {1\over \delta_\varepsilon}Y'\times (0,1))^3}\!\leq\! C\sigma_\varepsilon^{r\over r-1}\delta_\varepsilon^{-{2\over r}},
\end{equation}
\begin{equation}\label{unfolding2_2}
\left\Vert D_{z^{\prime}}\widehat {u}_{\varepsilon }\right\Vert_{L^r(\omega\times {1\over \delta_\varepsilon}Y'\times (0,1))^{{3\times2}}}\!\leq\! C\sigma_\varepsilon^{r\over r-1}, \quad\left\Vert \partial_{y_3}\widehat {u}_{\varepsilon }\right\Vert_{L^r(\omega\times {1\over \delta_\varepsilon}Y'\times (0,1))^3}\!\leq\! C\sigma_\varepsilon^{r\over r-1}  \delta_\varepsilon^{-{2\over r}},
\end{equation}
\begin{equation}\label{unfolding2_mean_2}
\left\Vert \widehat  u_{\varepsilon }-\bar u_{\varepsilon }\right\Vert_{L^r(\Omega; L^{r^*}(\mathbb{R}^2)^3)}\!\leq\! C\sigma_\varepsilon^{r\over r-1},
\end{equation}

\item[--] If $\sigma_\varepsilon\gg h_\varepsilon$, then
\begin{equation}\label{unfolding1_3}
\left\Vert \mathbb{D}_{z^{\prime}}\!\left[\widehat {u}_{\varepsilon }\right]\right\Vert_{L^r(\omega\times {1\over \delta_\varepsilon}Y'\times (0,1))^{3\times2}}\!\leq\! Ch_\varepsilon^{1\over r-1}\sigma_\varepsilon,\quad \left\Vert \partial_{y_3}\!\left[\widehat {u}_{\varepsilon }\right]\right\Vert_{L^r(\omega\times {1\over \delta_\varepsilon}Y'\times (0,1))^3}\!\leq\! Ch_\varepsilon^{r\over r-1}\delta_\varepsilon^{-{2\over r}},
\end{equation}
\begin{equation}\label{unfolding2_3}
\left\Vert D_{z^{\prime}}\widehat {u}_{\varepsilon }\right\Vert_{L^r(\omega\times {1\over \delta_\varepsilon}Y'\times (0,1))^{{3\times2}}}\!\leq\! Ch_\varepsilon^{1\over r-1}\sigma_\varepsilon, \quad\left\Vert \partial_{y_3}\widehat {u}_{\varepsilon }\right\Vert_{L^r(\omega\times {1\over \delta_\varepsilon}Y'\times (0,1))^3}\!\leq\! Ch_\varepsilon^{r\over r-1}  \delta_\varepsilon^{-{2\over r}},
\end{equation}
\end{itemize}
Moreover, in every case it holds
\begin{equation}\label{unfolding4}
\left\Vert \widehat {p}_{\varepsilon} \right\Vert_{L^{r'}(\omega\times \mathbb{R}^2\times (0,1))/\mathbb{R}}\leq C \delta_\varepsilon^{-{2\over r'}}.
\end{equation}
\end{lemma}
\begin{proof} We describe  the case $\sigma_\varepsilon \ll h_\varepsilon$ (the rest of the cases are similar, so we omit it).
Taking into account estimates (\ref{a1_thin_tilde}) in (\ref{estimate_varphi2_mean2}), we deduce (\ref{unfolding3}). We remark that if we had used (\ref{a1_thin_tilde}) and (\ref{estimate_varphi1}), we would have obtained estimate $\|\widehat  u_{\varepsilon\delta_\varepsilon}\|_{L^r(\omega\times \mathbb{R}^2\times (0,1))^3}\leq C \sigma_\varepsilon^{r\over r-1}\delta_\varepsilon^{-{2\over r}}$, which is not as sharp a (\ref{unfolding3}). 

Also, taking into account the second estimate in (\ref{a1_thin_tilde})  and (\ref{estimate_varphi2}), we get the first estimate in (\ref{unfolding1}). And using the second estimate in (\ref{a1_thin_tilde})  and (\ref{estimate_varphi3}), we get the second estimate in (\ref{unfolding1}). Consequently, from Korn's inequality (\ref{Korn}), we also have (\ref{unfolding2}). Estimate (\ref{estimate_varphi2_mean}) together with the second estimate in (\ref{a1_thin_tilde})  gives (\ref{unfolding2_mean}).

 We remark that in the case $\sigma_\varepsilon\gg h_\varepsilon$, we only give some estimates because, in that case, the microstructure will not play an important role in the homogenized model (see proof of Theorem \ref{MainTheorem}$-(iii)$), so we omit it.
 
 Finally, taking into account the estimate (\ref{estim_P_thin_tilde}) in (\ref{estimate_varphi1}) with $q=r'$, and $r^\prime=r/(r-1)$, we can deduce (\ref{unfolding4}).

\end{proof}

\subsection{Some compactness results}\label{sec:comp}
We give some compactness results concerning the behavior of the extension of the solution $(\widetilde  u_{\varepsilon },\widetilde  p_{\varepsilon })$ and the sequence of unfolded functions $(\widehat  u_{\varepsilon},\widehat  p_\varepsilon)$ satisfying the {\it a priori} estimates given in previous sections.
\begin{lemma}[Compactness results for extension of dilated functions]\label{LemmaConvtulde}Suppose $1<r<2$ and let $\sigma_\varepsilon$ be given by (\ref{sigma}) satisfying (\ref{sigma0}). Depending on the values of $\sigma_\varepsilon$ and $h_\varepsilon$:
\begin{itemize}
\item[--] If   $\sigma_\varepsilon\ll h_\varepsilon$, for a subsequence of $\varepsilon$, still denoted by $\varepsilon$, there exists  $u\in L^r(\Omega)^3$, with $u_3\equiv0$, such that
\begin{equation}\label{conv_u_tilde}
\sigma_\varepsilon^{-{r\over r-1}} \widetilde u_{\varepsilon }\rightharpoonup   u \quad\hbox{weakly in }L^r(\Omega)^3.
\end{equation}
\item[--] If $\sigma_\varepsilon\approx h_\varepsilon$, with $\sigma_\varepsilon/ h_\varepsilon\to \lambda\in (0,+\infty)$, for a subsequence of $\varepsilon$, still denoted by $\varepsilon$, there exists  $u\in W^{1,r}(0,1;L^r(\omega)^3)$, with $u=0$ on $(\omega\times \{0\})\cup (\omega\times \{1\})$ and $u_3\equiv 0$, such that 
\begin{equation}\label{conv_u_tilde2}
\sigma_\varepsilon^{-{r\over r-1}} \widetilde u_{\varepsilon }\rightharpoonup   u \quad\hbox{weakly in }W^{1,r}(0,1;L^r(\omega)^3).
\end{equation}
\item[--] If $\sigma_\varepsilon\gg h_\varepsilon$, for a subsequence of $\varepsilon$, still denoted by $\varepsilon$, there exists  $u\in W^{1,r}(0,1;L^r(\omega)^3)$, with $u=0$ on $(\omega\times \{0\})\cup (\omega\times \{1\})$ and $u_3\equiv 0$, such that 
\begin{equation}\label{conv_u_tilde23}
h_\varepsilon^{-{r\over r-1}} \widetilde u_{\varepsilon }\rightharpoonup   u \quad\hbox{weakly in }W^{1,r}(0,1;L^r(\omega)^3).
\end{equation}
\end{itemize}
Moreover, in every case it holds the following divergence property 
\begin{equation}\label{div_x_u_tilde}
{\rm div}_{x'}\left(\int_0^1   u'(x',y_3)\,dy_3\right)=0\ \hbox{ in  }\omega,\quad \left(\int_0^1  u'(x',y_3)\,dy_3\right)\cdot n=0\ \hbox{ on  }\partial \omega,
\end{equation}
and that there exist $p\in L^{r'}(\Omega)/\mathbb{R}$, independent of $y_3$, such that
\begin{equation}\label{conv_p_tilde}
 \widetilde  p_{\varepsilon }\to   p\quad\hbox{strongly in }  L^{r'}(\Omega)/\mathbb{R}.
\end{equation}
\end{lemma}
\begin{proof} We develop the case $\sigma_\varepsilon\ll h_\varepsilon$. First, we focus on the convergence of the extension of the extension of the dilated velocity. The first   estimate in (\ref{a1_thin_tilde}) implies the existence of $u\in L^r(\Omega)^3$  such that, up to a subsequence, we have the convergence
\begin{equation}\label{conv_u}\sigma_\varepsilon^{-{r\over r-1}}\widetilde u_{\varepsilon } \rightharpoonup   u\quad \hbox{in }L^r(\Omega)^3.
\end{equation}
This also implies
\begin{equation}\label{div_conv_u}
\sigma_\varepsilon^{-{r\over r-1}}{\rm div}_{x'}\widetilde u_{\varepsilon } \rightharpoonup {\rm div}_{x'}  u\quad \hbox{in }L^r(0,1;W^{-1,r'}(\omega)^3).
\end{equation}
From ${\rm div}_{h_\varepsilon} \widetilde u_{\varepsilon  }=0$ in $\Omega$, multiplying by $\sigma_\varepsilon^{-{r\over r-1}}$, we get
\begin{equation}\label{div_proof}
\sigma_\varepsilon^{-{r\over r-1}}{\rm div}_{x'}\widetilde u_{\varepsilon  }'+h_\varepsilon^{-1}\sigma_\varepsilon^{-{r\over r-1}}\partial_{y_3}\widetilde u_{\varepsilon ,3}=0\quad\hbox{ in }\Omega.
\end{equation}
Next, (\ref{conv_u}) and (\ref{div_conv_u}) combined with (\ref{div_proof}) imply that $h_\varepsilon^{-1}\sigma_\varepsilon^{-{r\over r-1}}\partial_{y_3}\widetilde u_{\varepsilon,3}$ is bounded in $L^r(0,1;W^{-1,r'}(\omega)^3)$.  This implies that $\sigma_\varepsilon^{-{r\over r-1}}\partial_{y_3}\widetilde u_{\varepsilon,3}$ tends to zero in $L^r(0,1;W^{-1,r'}(\omega)^3)$. Also, from (\ref{conv_u}), we have that $\sigma_\varepsilon^{-{r\over r-1}}\partial_{y_3}\widetilde u_{\varepsilon, 3}$ tends to $\partial_{y_3}\widetilde u_{3}$ in $W^{-1,r'}(0,1;L^r(\omega))$. From the uniqueness of the limit, we
have that $\partial_{y_3}\widetilde u_{3}=0$ which implies that $\widetilde u_3$ does not depend on $y_3$.

Next, we prove the divergence condition (\ref{div_x_u_tilde}). For this, we consider $\varphi\in C^1_0(\omega)$ as test function in ${\rm div}_{h_\varepsilon} \widetilde  u_{\varepsilon }=0$ in $\Omega$, and we get 
$$0=\int_\Omega \left(\sigma_\varepsilon^{-{r\over r-1}}{\rm div}_{x'} \widetilde  u'_{\varepsilon }+ h_\varepsilon^{-1}\sigma_\varepsilon^{-{r\over r-1}} \partial_{y_3}\widetilde  u_{\varepsilon,3}\right)\, \varphi(x')\,dx'dy_3=-\int_\Omega  \sigma_\varepsilon^{-{r\over r-1}}\widetilde  u'_{\varepsilon}\cdot \nabla_{x'}\varphi(x')\,dx'dy_3,$$
and from convergence (\ref{conv_u_tilde}), we deduce
$$\int_\Omega    u'\cdot \nabla_{x'}\varphi(x')\,dx'dy_3=0,$$
which  implies (\ref{div_x_u_tilde}). 

Now, we focus on the extension of the dilated pressure by proving (\ref{conv_p_tilde}). From estimate (\ref{estim_P_thin_tilde}) we deduce that there exists $ p\in L^{r'}(\Omega)/\mathbb{R}$ such that, up to a subsequence,  we have the following convergence
\begin{equation}\label{conv_p_tilde1}
 \widetilde p_{\varepsilon}\rightharpoonup    p\quad\hbox{weakly in } L^{r'}(\Omega)/\mathbb{R}.
\end{equation}
In order to prove that $ p$ does not depend on $y_3$, we multiply system (\ref{2}) by $h_\varepsilon \,  \varphi$ with $  \varphi\in C^\infty_0(\Omega)^3$ and  we integrate by parts. Thus, we get
\begin{equation}\label{press}
\begin{array}{l}
\medskip
\displaystyle
\int_{\Omega}h_\varepsilon \widetilde p_{\varepsilon }{\rm div}_{x'}  \varphi'\,dx'dy_3+\int_{\Omega}  \widetilde p_{\varepsilon }\,\partial_{y_3}  \varphi\,dx'dy_3\\
\medskip
\displaystyle =\mu h_\varepsilon \int_{ {\Omega}}\left\vert \mathbb{D}_{h_\varepsilon} \left[\widetilde{u}_{\varepsilon }\right]\right\vert^{r-2}\mathbb{D}_{h_\varepsilon}\left[ \widetilde {u}_{\varepsilon }\right]:\mathbb{D} \left[ \varphi \right]dx'dy_3-h_\varepsilon \int_{ {\Omega}}f'\cdot   \varphi'\,dx'dy_3.
\end{array}
\end{equation}
Taking into account estimates of the extension of the velocity  and pressure given in Lemma \ref{estOmega} (case $\sigma_\varepsilon \ll h_\varepsilon$), and the assumption of force $f$, we have
$$\begin{array}{l}
\medskip
\displaystyle \left|\mu h_\varepsilon \int_{ {\Omega}}\left\vert \mathbb{D}_{h_\varepsilon} \left[\widetilde{u}_{\varepsilon }\right]\right\vert^{r-2}\mathbb{D}_{h_\varepsilon}\left[ \widetilde {u}_{\varepsilon }\right]:\mathbb{D} \left[\varphi \right]dx'dy_3\right|\leq Ch_\varepsilon  \|\mathbb{D}_{h_\varepsilon}[\widetilde u_{\varepsilon  }]\|_{L^r(\Omega)^{3\times 3}}^{r-1}\leq Ch_\varepsilon \sigma_{\varepsilon} \to 0,\\
\medskip
\displaystyle \left|h_\varepsilon \int_{ {\Omega}}f'\cdot  \varphi'\,dx'dy_3\right|\leq Ch_\varepsilon \to 0,\\
\medskip
\displaystyle \left|h_\varepsilon \int_{\Omega}\widetilde p_{\varepsilon }{\rm div}_{x'}  \varphi'\,dx'dy_3\right|\leq Ch_\varepsilon \|\widetilde p_{\varepsilon }\|_{L^{r'}(\Omega)}\leq Ch_\varepsilon\to 0.
\end{array}
$$
Then,   passing to the limit in (\ref{press}) by means of convergence (\ref{conv_p_tilde1}), we obtain
$$\int_{\Omega}   p\,\partial_{y_3} \varphi\,dx'dy_3=0,$$
which implies that $ p$ is independent of $y_3$. Moreover, if we argue similarly as in \cite[Lemma 4.4]{Bourgeat1}, we have that the convergence of the pressure $\widetilde  p_{\varepsilon}$ is in fact strong. 

To conclude the proof, it remains to prove that $\widetilde u_3\equiv 0$. To do this, since $\widetilde u_3$ does not depend on $y_3$, we take as test function $\varphi=(0,\sigma_\varepsilon^{-{r\over r-1}}\varphi_3(x'))$ in (\ref{2}), and passing to the limit, using monotonicity arguments (follow Step 1 of the proof of Theorem \ref{MainTheorem} without applying the unfolding method, and taking into account the test function and the definition of the symmetric derivative (\ref{derivada_def}) in this case), we can deduce that $\widetilde u_3$ satisfies 
$$-{\rm div}_{x'}\left(|\mathbb{D}_{x'}[\widetilde u_3]|^{r-2}\mathbb{D}_{x'}[\widetilde u_3]\right)=0\quad \hbox{in }\Omega,$$
and this implies $\widetilde u_3\equiv 0$, which concludes the proof.

The cases $\sigma_\varepsilon\approx h_\varepsilon$ and $\sigma_\varepsilon \gg h_\varepsilon$ are similar to the case developed in \cite[Lemma 5.1]{Anguiano_SG_Lower}, so we omit it. The main difference with the previous case is that, from estimates  (\ref{a1_thin2_tilde}) and (\ref{a1_thin3_tilde}) respectively, the velocity has more regularity in the vertical variable, i.e. we deduce $\widetilde u \in W^{1,r}(0,1;L^r(\omega)^3)$, and so we can deduce $u=0$ on $(\omega\times \{0\})\cup (\omega\times \{1\})$. Moreover, since $u_3$ is independent of $y_3$ together with the boundary conditions on the top and bottom, it holds $u_3\equiv 0$, whose proof is faster than the one in the previous case.

\end{proof}

Now, we give a compactness result for the unfolded velocity $\widehat u_{\varepsilon }$ (only in the cases $\sigma_\varepsilon \ll h_\varepsilon$ and $\sigma_\varepsilon\approx h_\varepsilon$). For this, we follow \cite[Chapters 9 and 10]{Cioran-book} and consider the homogeneous Sobolev space of weakly differentiable functions defined locally on $\mathbb{R}^2$  having a gradient in $L^{r}(\mathbb{R}^2)^2$  and  zero value on the obstacle $Y'_s$, which is given by
\begin{equation}\label{KYs}
{\bf K}_{Y'_s}=\left\{\varPhi(z')\in W^{1,r}_{loc}(\mathbb{R}^2)\,:\, \nabla_{z'}\varPhi\in L^r(\mathbb{R}^2)^{2}\,\hbox{ and }\,\varPhi=0\ \hbox{ on }Y'_s\right\}.
\end{equation}
We remark that if $\varPhi\in {\bf K}_{Y'_s}$ then it has a limit at infinity denoted $\varPhi_\infty$, i.e. there exists $\varPhi_\infty\in \mathbb{R}$ such that $\lim_{|z'|\to +\infty}\varPhi(z')=\varPhi_\infty$.
In addition,  to relate the value at infinity of the limit of $\widehat u_{\varepsilon }$ with the limit of $\widetilde u_{\varepsilon }$, we consider a more general space depending on the case
\begin{itemize}
\item[--] In the case $\sigma_\varepsilon\approx h_\varepsilon$,  with $\sigma_\varepsilon/ h_\varepsilon\to \lambda\in (0,+\infty)$:
$$
{\bf L}_{Y'_s}=\Big\{\varPhi(x',z',y_3)\in  L^r(\Omega; {\bf K}_{Y'_s})\,:\, \varPhi_\infty=\varPhi(\cdot,\infty,\cdot)\in W^{1,r}(0,1;L^r(\omega))\quad \hbox{with}\quad \varPhi_\infty=0\ \hbox{on}\ y_3=\{0,1\}\Big\}.
$$
\item[--] In the case $\sigma_\varepsilon\ll h_\varepsilon$:
$$
{\bf L}_{Y'_s}=\Big\{\varPhi(x',z',y_3)\in  L^r(\Omega; {\bf K}_{Y'_s})\,:\, \varPhi_\infty=\varPhi(\cdot,\infty,\cdot)\in L^r(\Omega)\Big\}.
$$
\end{itemize}

\begin{lemma}[Compactness results for unfolded velocity]\label{lem:compactnesshat} Consider the cases  $\sigma_\varepsilon\ll h_\varepsilon$ and $\sigma_\varepsilon\approx h_\varepsilon$, with $\sigma_\varepsilon/ h_\varepsilon\to \lambda\in (0,+\infty)$. Suppose $1<r<2$,  let $\sigma_\varepsilon$ be given by (\ref{sigma}) satisfying (\ref{sigma0}) and $u$ be given in Lemma \ref{LemmaConvtulde} according to the corresponding case. Then, for a subsequence of $\varepsilon$, still denoted by $\varepsilon$, there exists $U\in {\bf L}_{Y'_s}^3$ where $U_{\infty}=  u$ and $U_{3}$ independent of $y_3$, such that 
\begin{equation}\label{conv_u_hat}
\sigma_\varepsilon^{-{r\over r-1}}\widehat  u_{\varepsilon }\rightharpoonup U\quad\hbox{weakly in }L^r(\Omega;L_{loc}^{r}(\mathbb{R}^2)^3)\,,
\end{equation}
\begin{equation}\label{conv_Du_hat}
\sigma_\varepsilon^{-{r\over r-1}} D_{z'}\widehat  u_{\varepsilon } {\textbf{1}_{{1\over \delta}Y'}}\rightharpoonup D_{z'}U\quad\hbox{ weakly in } L^r(\omega\times \mathbb{R}^2\times (0,1))^{3\times 2},
\end{equation}
\begin{equation}\label{divx_u_hat}
{\rm div}_{x'}\left(\int_0^1  U'_\infty\,dy_3\right)=0\quad \hbox{in }\omega,
\end{equation}
\begin{equation}\label{div_u_hat}
{\rm div}_{z'}U'=0\quad\hbox{in }\omega\times \mathbb{R}^2\times (0,1)\,.
\end{equation}
\end{lemma}
\begin{proof} 
 
We only consider the cases $\sigma_\varepsilon\ll h_\varepsilon$ and $\sigma_\varepsilon\approx h_\varepsilon$, because in the case $\sigma_\varepsilon\gg h_\varepsilon$, as we said before, the microstructure will not play an important role in the homogenized model (see proof of Theorem \ref{MainTheorem}).

The proof of the existence of $U\in {\bf L}_{Y'_s}^3$ where $U_{\infty}=  u$ and $U_3$ independent of $y_3$, such that (\ref{conv_u_hat}) and (\ref{conv_Du_hat})  hold,  is the same as the one given in \cite[Lemma 5.1]{Anguiano_SG_Lower}, so we omit it.

Next,  since $U_\infty=  u$, then (\ref{divx_u_hat}) holds  from the divergence condition (\ref{div_x_u_tilde}). Finally, from divergence condition ${\rm div}_{h_\varepsilon}\widetilde  u_{\varepsilon }=0$ in $ \Omega $ and  the change of variables  (\ref{CV}) (see the proof of Theorem 4.1in \cite{Anguiano_SG_Lower} for more details), we deduce 
\begin{eqnarray}\label{divergence_hat_inicial}
 (\varepsilon\delta_\varepsilon)^{-1}{\rm div}_{z'}\widehat  u'_{\varepsilon }\,{\bf 1}_{{1\over \delta_\varepsilon}Y'}+h_\varepsilon^{-1}\partial_{y_3}\widehat  u_{\varepsilon,3}\,{\bf 1}_{{1\over \delta_\varepsilon}Y'}=0\quad\hbox{ in }\omega\times \mathbb{R}^2\times (0,1).
 \end{eqnarray}
Multiplying by 
$\sigma_{\varepsilon}^{-{r\over r-1}}(\varepsilon\delta_\varepsilon)$, we get
\begin{eqnarray}\label{divergence_hat}
\sigma_{\varepsilon}^{-{r\over r-1}}{\rm div}_{z'}\widehat  u'_{\varepsilon }\,{\bf 1}_{{1\over \delta_\varepsilon}Y'}+\delta_\varepsilon^{2\over r} h_\varepsilon^{-1}\sigma_{\varepsilon}^{-{1\over r-1}}\partial_{y_3}\widehat  u_{\varepsilon,3}\,{\bf 1}_{{1\over \delta_\varepsilon}Y'}=0\quad\hbox{ in }\omega\times \mathbb{R}^2\times (0,1).
\end{eqnarray}
For the case $\sigma_\varepsilon\ll h_\varepsilon$, from  the second estimate in (\ref{unfolding2}) and convergence (\ref{conv_u_hat}), together with the fact that $\delta_\varepsilon^{2\over r} h_\varepsilon^{-1}\sigma_{\varepsilon}^{-{1\over r-1}}=\delta_\varepsilon^{2\over r} \sigma_\varepsilon h_\varepsilon^{-1}\sigma_{\varepsilon}^{-{r\over r-1}}$ and  $\delta_\varepsilon^{2\over r} \sigma_\varepsilon h_\varepsilon^{-1}\to 0$, we deduce that $\delta_\varepsilon^{2\over r} h_\varepsilon^{-1}\sigma_{\varepsilon}^{-{1\over r-1}}\partial_{y_3}\widehat  u_{\varepsilon,3}$ tends to zero, and so passing to the limit in (\ref{divergence_hat}) we get  (\ref{div_u_hat}). For the case $\sigma_\varepsilon\approx h_\varepsilon$, we proceed similarly because $\delta_\varepsilon^{2\over r} \sigma_\varepsilon h_\varepsilon^{-1}$ also tends to zero.  

\end{proof}

\subsection{Homogenized model: proof of the main theorem}\label{sec:proofthm}
In this section, we  prove the main result. For this, we recall   the following version of \cite[Lemma 10.4]{Cioran-book}  to choose an appropriate test function in the variational formulation of system (\ref{2}) and thus, be able to pass to  the limit.
\begin{lemma}\label{ftestlimit} Suppose $1\leq q<+\infty$.  Let $\varphi$ be in $\mathcal{D}(\Omega;  W^{1,q}_{loc}(\mathbb{R}^2))$ such that $\nabla_{z'}\varphi$ is in $\mathcal{D}(\Omega; L^q(\mathbb{R}^2)^2)$ and has a compact support. We set
\begin{equation}\label{vtestepdelta}\displaystyle \varphi_{\varepsilon\delta_\varepsilon}(x',y_3)=\varphi\left(x^{\prime},{1\over \delta_\varepsilon}{x^{\prime}-\varepsilon\kappa({x'\over \varepsilon})\over \varepsilon},y_3\right)\quad \hbox{in } (x',y_3)\in \Omega.
\end{equation}
By \cite[Proposition 9.2]{Cioran-book}, it  has a limit at infinity denoted by $\varphi_\infty\in \mathcal{D}(\Omega)$. If $\delta_\varepsilon$ is small enough,  the function $\varphi_{\varepsilon\delta_\varepsilon}$ belongs to  $\mathcal{D}(0,1;W^{1,q}(\omega))$ and 
\begin{equation}\label{vepdeltest}
\varphi_{\varepsilon\delta_\varepsilon}\to \varphi_\infty\quad\hbox{strongly in }L^r(\Omega).
\end{equation}
\end{lemma}

\begin{remark} \label{remftest} From the definition of $\varphi_{\varepsilon\delta_\varepsilon}$ given in (\ref{vtestepdelta}), we have
\begin{equation}\label{derivativetest}
\begin{array}{l}
\medskip
\displaystyle
\nabla_{x'} \varphi_{\varepsilon\delta_\varepsilon}(x',y_3)=\nabla_{x'}\varphi\left(x^{\prime},{1\over \delta_\varepsilon}{x^{\prime}-\varepsilon\kappa({x'\over \varepsilon})\over \varepsilon},y_3\right)+(\varepsilon\delta_\varepsilon)^{-1}\nabla_{z'}\varphi\left(x^{\prime},{1\over \delta_\varepsilon}{x^{\prime}-\varepsilon\kappa({x'\over \varepsilon})\over \varepsilon},y_3\right), 
\\
\medskip
\displaystyle\partial_{y_3}\varphi_{\varepsilon\delta_\varepsilon}(x',y_3)=\partial_{y_3}\varphi\left(x^{\prime},{1\over \delta_\varepsilon}{x^{\prime}-\varepsilon\kappa({x'\over \varepsilon})\over \varepsilon},y_3\right). 
\end{array}
\end{equation}
Moreover, applying the unfolding operator (\ref{hat}), we have
$$\widehat  \varphi_{\varepsilon\delta_\varepsilon}(x',z',y_3)=\left\{\begin{array}{ll}
\varphi(x',z',y_3) + \varTheta_{\varepsilon\delta_\varepsilon}(x',z',y_3),&\hbox{if  }(x',z',y_3)\in \omega\times  {1\over \delta_\varepsilon}Y'\times (0,1),
 \\
 \\
0&\hbox{otherwise,}
 \end{array}\right.
$$
with  $\varTheta_{\varepsilon\delta_\varepsilon}(x',z',y_3)=\varphi(\varepsilon \kappa({x'\over \varepsilon})+\varepsilon\delta_\varepsilon z', z', y_3)-\varphi(x',z',y_3)$.  Consequently
\begin{equation}\label{derivadatest}\nabla_{z'}\widehat  \varphi_{\varepsilon\delta_\varepsilon}(x',z',y_3)= \nabla_{z'}\varphi(x',z',y_3)+  \nabla_{z'}\varTheta_{\varepsilon\delta_\varepsilon}(x',z',y_3)\quad\hbox{in }\omega\times {1\over \delta_\varepsilon}Y'\times (0,1),
\end{equation}
where, from the mean value theorem applied to $\nabla_{z'}\varTheta_{\varepsilon\delta_\varepsilon}$, the fact that $|\varepsilon \kappa({x'\over \varepsilon})+\varepsilon\delta_\varepsilon z'-x'|<\varepsilon$ for  $x'\in Y'_{k',\varepsilon}$, $k'\in \mathcal{K}_{\varepsilon}$, and  $\nabla_{z'}\varphi\in \mathcal{D}(\Omega;L^{q}(\mathbb{R}^2)^2)$, it holds
\begin{equation}\label{convtestder}
\|\nabla_{z'} \varTheta_{\varepsilon\delta_\varepsilon}\|_{L^q(\omega\times \mathbb{R}^2\times (0,1))^2}\leq C\varepsilon.
\end{equation}
\end{remark}

\begin{proof}[Proof of Theorem \ref{MainTheorem}] The proof of the main result will be divided in four steps.\\

\noindent {\bf Step 1.}  Cases $\sigma_\varepsilon \ll h_\varepsilon$ and $\sigma_\varepsilon\approx h_\varepsilon$. In this step, we derive a variational inequality, which will be useful to derive their respective homogenized variational formulations. \\

We define the following set
$$\mathbb{W}=\left\{\begin{array}{rcl}
\medskip 
( v', V')\in W^{1,r}(0,1;L^r(\omega)^2)\times {\bf L}_{Y'_s}^2& :&  V'_\infty(x',y_3)= v'(x',y_3)\quad\hbox{a.e. in }(x',y_3)\in\Omega,\\
\medskip
 {\rm div}_{z'}V'=0\hbox{ in }\omega\times \mathbb{R}^2\times (0,1), &&\displaystyle {\rm div}_{x'}\left(\int_0^1v'\,dy_3\right)=0\hbox{ in }\omega,\quad \left(\int_0^1v'\,dy_3\right)\cdot n'=0\hbox{ on }\partial\omega,
\end{array}\right\}.$$
To simplify expressions, we define  the operator  $S:\mathbb{R}^3_{{\rm sym}}\to \mathbb{R}^3_{{\rm sym}}$ by 
$$S(\xi)=|\xi|^{r-2}\xi,\quad\forall\,\xi\in\mathbb{R}^{3\times 3}_{{\rm sym}},$$
and denote by $O_\varepsilon$ a generic real sequence, which tends to zero with $\varepsilon$  and can change from line to line.

Now, let us define the test function according to Lemmas \ref{LemmaConvtulde} and \ref{lem:compactnesshat} in the cases $\sigma_\varepsilon \ll h_\varepsilon$ and $\sigma_\varepsilon \approx h_\varepsilon$. We consider $\varphi(x^{\prime},z',y_3)\in \mathcal{D}(\Omega; {\bf K}_{Y'_s}^3)$ such that $D_{z'}\varphi$ has a compact support and $\varphi_3\equiv 0$, and we define $\varphi_{\varepsilon\delta}$  by (\ref{vtestepdelta}). According to Lemma \ref{ftestlimit},  when $\delta_\varepsilon$ is small enough  the function $\varphi_{\varepsilon\delta_\varepsilon}$ belongs to $\mathcal{D}(0,1;W^{1,r}(\omega)^3)$ and its limit at infinity $\varphi_\infty\in \mathcal{D}(\Omega)^3$ satisfies the convergence (\ref{vepdeltest}).

Multiplying (\ref{2}) by $\displaystyle \varphi_{\varepsilon\delta_\varepsilon}$, taking into account the extensions of the dilated velocity and pressure, integrating by parts and using (\ref{derivativetest}), we have
\begin{equation}\label{FV_tilde}\begin{array}{l}
\medskip\displaystyle
\mu\,  \int_{\Omega}S\left(\mathbb{D}_{h_\varepsilon} \left[\widetilde {u}_{\varepsilon } \right] \right):\left(\mathbb{D}_{x^\prime}\left[\varphi\right]+ h_\varepsilon^{-1}\partial_{y_3}\left[\varphi\right]\right) dx^{\prime}dy_3+\mu ({\varepsilon\delta})^{-1}   \int_{\Omega}S\left(\mathbb{D}_{h_\varepsilon} \left[\widetilde {u}_{\varepsilon } \right] \right):\mathbb{D}_{z^\prime}\left[\varphi\right]\,dx'dy_3\\
\medskip
\displaystyle - \int_{\Omega} \widetilde {p}_{\varepsilon }\, {\rm div}_{x^\prime} \varphi^\prime\,dx^{\prime}dy_3- (\varepsilon\delta)^{-1} \int_{\Omega} \widetilde {p}_{\varepsilon }\, {\rm div}_{z^\prime} \varphi^\prime\,dx^{\prime}dy_3= \int_{\widetilde \Omega_{\varepsilon }}f'\cdot \varphi'\,dx^\prime dy_3\,,
\end{array}
\end{equation}
where, for simplify,  from now on we use the following notation:
\begin{itemize}
\item[--] $\varphi=\varphi(x^{\prime},{1\over \delta_\varepsilon}{x^{\prime}-\varepsilon\kappa({x'\over \varepsilon})\over \varepsilon},y_3)$ in the integrals in $\Omega$, 
\item[--] $\varphi=\varphi(x',z',y_3)$ in the integrals  in  $\omega\times {1\over \delta_\varepsilon}Y'\times (0,1)$ obtained after applying the changes of variables (\ref{CV}).
\end{itemize}
Below, we  analyze every term in (\ref{FV_tilde}):
\begin{itemize}

\item[--]  First term. Taking into account the convergences of the velocity given in Lemma \ref{LemmaConvtulde} in both cases, which are given by (\ref{conv_u_tilde}) and (\ref{conv_u_tilde2}) respectively, we rewrite this term as follows
\begin{equation}\label{integral_velocity11}
\begin{array}{l}
\medskip
\displaystyle
\mu\,  \int_{\Omega}S\left(\mathbb{D}_{x'} \left[\widetilde {u}_{\varepsilon } \right] + h_{\varepsilon}^{-1}\partial_{y_3}[\widetilde  u_{\varepsilon }] \right): \left( \mathbb{D}_{x^\prime}\left[\varphi\right] + h_{\varepsilon}^{-1}\partial_{y_3}\left[\varphi\right]\right)\,dx'dy_3\\
\medskip
\displaystyle
=\mu\,  \int_{\Omega}S\left(\sigma_{\varepsilon}\mathbb{D}_{x'}\left[\sigma_{\varepsilon}^{-{r\over r-1}}\widetilde {u}_{\varepsilon }\right] + \sigma_\varepsilon h_\varepsilon^{-1}\partial_{y_3}\left[\sigma_{\varepsilon}^{-{r\over r-1}}\widetilde {u}_{\varepsilon }\right] \right): \left(\sigma_{\varepsilon}\mathbb{D}_{x^\prime}\left[\varphi\right] + \sigma_\varepsilon h_\varepsilon^{-1}\partial_{y_3}\left[\varphi\right]\right)\,dx'dy_3.
\end{array}
\end{equation}

\item[--] Second term.    We apply  the changes of variables  (\ref{CV}) (see the proof of Theorem 4.1in \cite{Anguiano_SG_Lower} for more details) and taking into account the convergences of the unfolded velocity given in Lemma \ref{lem:compactnesshat} in both cases together with the property of the test function given by (\ref{derivadatest}), we rewrite the second term as follows
\begin{eqnarray}
\displaystyle
&&\mu\,  (\varepsilon\delta_\varepsilon)^{-1}\int_{\Omega}S\left(\mathbb{D}_{x'} \left[\widetilde {u}_{\varepsilon } \right] + h_{\varepsilon}^{-1}\partial_{y_3}[\widetilde  u_{\varepsilon }] \right): \mathbb{D}_{z^\prime}\left[\varphi\right] \,dx'dy_3  \nonumber\\
\medskip
\displaystyle
&&=\mu\, (\varepsilon\delta_\varepsilon)^{-1}\delta_\varepsilon^2 \int_{\omega\times{1\over \delta_\varepsilon}Y'\times (0,1)}S\left( (\varepsilon\delta_\varepsilon)^{-1}\mathbb{D}_{z'}\left[\widehat {u}_{\varepsilon }\right] +h_{\varepsilon}^{-1}\partial_{y_3}\left[\widehat  {u}_{\varepsilon }\right] \right):(\mathbb{D}_{z^\prime}\left[\varphi\right] + \mathbb{D}_{z'}[\varTheta_{\varepsilon\delta_\varepsilon}])\,dx'dz'dy_3 \label{integral_velocity2_bis1}
\\
\medskip
\displaystyle
&&=\mu\, \int_{\omega\times{1\over \delta_\varepsilon}Y'\times (0,1)}S\left(\mathbb{D}_{z'}\left[\sigma_{\varepsilon}^{-{r\over r-1}}\widehat {u}_{\varepsilon }\right] + \delta_\varepsilon^{2\over r}\sigma_\varepsilon  h_\varepsilon^{-1}\partial_{y_3}\left[\sigma_{\varepsilon}^{-{r\over r-1}}\widehat  {u}_{\varepsilon }\right] \right):
(\mathbb{D}_{z^\prime}\left[\varphi\right]+ \mathbb{D}_{z'}[\varTheta_{\varepsilon\delta_\varepsilon}])\,dx'dz'dy_3.\nonumber
\end{eqnarray}
Using the H${\rm \ddot{o}}$lder inequality, estimates given in (\ref{unfolding1}) in the case $\sigma_\varepsilon \ll h_\varepsilon$ or (\ref{unfolding1_2}) in the case $\sigma_\varepsilon \approx h_\varepsilon$, and taking into account (\ref{convtestder}), we have 
$$
\left|\int_{\omega\times{1\over \delta_\varepsilon}Y'\times (0,1)}S\left(\mathbb{D}_{z'}\left[\sigma_{\varepsilon}^{-{r\over r-1}}\widehat {u}_{\varepsilon}\right] + \delta_\varepsilon^{2\over r}\sigma_\varepsilon  h_\varepsilon^{-1}\partial_{y_3}\left[\sigma_{\varepsilon}^{-{r\over r-1}}\widehat  {u}_{\varepsilon}\right] \right):
 \mathbb{D}_{z'}[\varTheta_{\varepsilon\delta_\varepsilon}]\,dx'dz'dy_3\right|\leq C \varepsilon\to 0,$$
and then, we have that (\ref{integral_velocity2_bis1}) is rewritten as follows
\begin{eqnarray}
&&\mu\,  (\varepsilon\delta_\varepsilon)^{-1}\int_{\Omega}S\left(\mathbb{D}_{x'} \left[\widetilde {u}_{\varepsilon } \right] + h_{\varepsilon}^{-1}\partial_{y_3}[\widetilde  u_{\varepsilon}] \right): \mathbb{D}_{z^\prime}\left[\varphi\right] \,dx'dy_3  \nonumber\\
\medskip
\displaystyle
&&=\mu\, \int_{\omega\times{1\over \delta_\varepsilon}Y'\times (0,1)}S\left(\mathbb{D}_{z'}\left[\sigma_{\varepsilon}^{-{r\over r-1}}\widehat {u}_{\varepsilon }\right] + \delta_\varepsilon^{2\over r}\sigma_\varepsilon  h_\varepsilon^{-1}\partial_{y_3}\left[\sigma_{\varepsilon}^{-{r\over r-1}}\widehat  {u}_{\varepsilon }\right] \right):
\mathbb{D}_{z^\prime}\left[\varphi\right]\,dx'dz'dy_3 +O_\varepsilon \label{integral_velocity2bisbis}
\\
\medskip
\displaystyle
&&=\mu\, \int_{\omega\times{1\over \delta_\varepsilon}Y'\times (0,1)} S\left(\mathbb{D}_{z'}\left[\sigma_{\varepsilon}^{-{r\over r-1}}\widehat {u}_{\varepsilon}\right] + \delta^{2\over r}\sigma_\varepsilon  h_\varepsilon^{-1}\partial_{y_3}\left[\sigma_{\varepsilon}^{-{r\over r-1}}\widehat  {u}_{\varepsilon}\right] \right)
:\left(\mathbb{D}_{z^\prime}\left[\varphi\right]+\delta^{2\over r}_\varepsilon\sigma_\varepsilon  h_\varepsilon^{-1}\partial_{y_3}\left[\varphi\right]\right)\,dx'dz'dy_3\nonumber
\\
\medskip
\displaystyle
&&-\mu\,\delta_\varepsilon^{2\over r}\sigma_\varepsilon  h_\varepsilon^{-1}\, \int_{\omega\times{1\over \delta_\varepsilon}Y'\times (0,1)}S\left(\mathbb{D}_{z'}\left[\sigma_{\varepsilon}^{-{r\over r-1}}\widehat {u}_{\varepsilon }\right] +\delta_\varepsilon^{2\over r}\sigma_\varepsilon  h_\varepsilon^{-1}\partial_{y_3}\left[\sigma_{\varepsilon}^{-{r\over r-1}}\widehat  {u}_{\varepsilon}\right] \right)
:\partial_{y_3}\left[\varphi\right] \,dx'dz'dy_3 +O_\varepsilon.\nonumber
\end{eqnarray}
Using that $\varphi\in\mathcal{D}(\Omega;W^{1,r}_{loc}(\mathbb{R}^2)^3)$, by the H${\rm \ddot{o}}$lder inequality and estimates in (\ref{unfolding1}) in the case $\sigma_\varepsilon \ll h_\varepsilon$ or (\ref{unfolding1_2}) in the case $\sigma_\varepsilon \approx h_\varepsilon$, we have
$$\left|\delta_\varepsilon^{2\over r}\sigma_\varepsilon  h_\varepsilon^{-1}\, \int_{\omega\times{1\over \delta_\varepsilon}Y'\times (0,1)}S\left(\mathbb{D}_{z'}\left[\sigma_{\varepsilon}^{-{r\over r-1}}\widehat {u}_{\varepsilon }\right] + \delta_\varepsilon^{2\over r}\sigma_\varepsilon  h_\varepsilon^{-1}\partial_{y_3}\left[\sigma_{\varepsilon}^{-{r\over r-1}}\widehat  {u}_{\varepsilon }\right]\right):\partial_{y_3}[\varphi]\,dx'dz'dy_3\right|\leq C\delta_\varepsilon^{2\over r}\sigma_\varepsilon  h_\varepsilon^{-1}\to 0.$$
Then, we have that (\ref{integral_velocity2bisbis}) reads as follows
\begin{eqnarray}
&&\mu\,  (\varepsilon\delta_\varepsilon)^{-1}\int_{\Omega}S\left(\mathbb{D}_{x'} \left[\widetilde {u}_{\varepsilon } \right] + h_{\varepsilon}^{-1}\partial_{y_3}[\widetilde  u_{\varepsilon }] \right): \mathbb{D}_{z^\prime}\left[\varphi\right] \,dx'dy_3  \label{integral_velocity2}\\
\medskip
\displaystyle
&&=\mu\, \int_{\omega\times{1\over \delta_\varepsilon}Y'\times (0,1)} S\left(\mathbb{D}_{z'}\left[\sigma_{\varepsilon}^{-{r\over r-1}}\widehat {u}_{\varepsilon }\right] + \delta_\varepsilon^{2\over r}\sigma_\varepsilon  h_\varepsilon^{-1}\partial_{y_3}\left[\sigma_{\varepsilon}^{-{r\over r-1}}\widehat  {u}_{\varepsilon }\right] \right)
:\left(\mathbb{D}_{z^\prime}\left[\varphi\right]+\delta_\varepsilon^{2\over r}\sigma_\varepsilon  h_\varepsilon^{-1}\partial_{y_3}\left[\varphi\right]\right)\,dx'dz'dy_3+O_\varepsilon.\nonumber
\end{eqnarray}

\item[--] Third and fourth terms. By the changes of variables  (\ref{CV}) (see the proof of Theorem 4.1 in \cite{Anguiano_SG_Lower} for more details) in the third and  fourth terms in (\ref{FV_tilde}) and  taking into account (\ref{derivadatest}), we deduce 
\begin{equation}\label{form_var_pressurebis}
\begin{array}{l}\medskip
\displaystyle
- \int_{\Omega} \widetilde {p}_{\varepsilon }\, {\rm div}_{x^\prime}\varphi^\prime\,dx^{\prime}dy_3- (\varepsilon\delta_\varepsilon)^{-1} \int_{\Omega} \widetilde {p}_{\varepsilon }\, {\rm div}_{z^\prime} \varphi^\prime\,dx^{\prime}dy_3\\
\medskip
\displaystyle
=- \int_{\Omega} \widetilde {p}_{\varepsilon }\, {\rm div}_{x^\prime} \varphi^\prime\,dx^{\prime}dy_3-\delta_\varepsilon^2 (\varepsilon\delta_\varepsilon)^{-1} \int_{\omega\times \mathbb{R}^2\times (0,1)} \widehat {p}_{\varepsilon }\, \left({\rm div}_{z^\prime} \varphi^\prime+ {\rm div}_{z^\prime} \varTheta ^\prime_{\varepsilon\delta_\varepsilon}\right)\,dx^{\prime}dz'dy_3.
\end{array}
\end{equation}

Using the H${\rm \ddot{o}}$lder inequality, estimate  (\ref{unfolding4}) and taking into account (\ref{convtestder}), we have
$$\left|\delta_\varepsilon^2  (\varepsilon\delta_\varepsilon)^{-1}\int_{\omega\times \mathbb{R}^2\times (0,1)} \widehat {p}_{\varepsilon}\, {\rm div}_{z^\prime} \varTheta^\prime_{\varepsilon\delta_\varepsilon}\,dx^{\prime}\,dz'dy_3\right|\leq C\delta_\varepsilon^{2-r\over r}\to 0.
$$
Thus, we get that (\ref{form_var_pressurebis}) is rewritten as follows 
\begin{equation}\label{form_var_pressure}
\begin{array}{l}\medskip
\displaystyle
- \int_{\Omega} \widetilde {p}_{\varepsilon}\, {\rm div}_{x^\prime} \varphi^\prime\,dx^{\prime}dy_3- (\varepsilon\delta_\varepsilon)^{-1} \int_{\Omega} \widetilde {p}_{\varepsilon }\, {\rm div}_{z^\prime} \varphi^\prime\,dx^{\prime}dy_3\\
\medskip
\displaystyle
=- \int_{\Omega} \widetilde {p}_{\varepsilon }\, {\rm div}_{x^\prime} \varphi^\prime\,dx^{\prime}dy_3-\delta_\varepsilon^2 (\varepsilon\delta_\varepsilon)^{-1}\int_{\omega\times \mathbb{R}^2\times (0,1)} \widehat {p}_{\varepsilon}\,{\rm div}_{z^\prime} \varphi^\prime \,dx^{\prime}dz'dy_3+O_\varepsilon.
\end{array}
\end{equation}

\end{itemize}

Therefore, taking into account (\ref{integral_velocity11}), (\ref{integral_velocity2}) and (\ref{form_var_pressure}) in (\ref{FV_tilde}), we get the following variational formulation
\begin{equation}\label{formvarcvuprime}
\begin{array}{l}
\medskip
\displaystyle
\mu\,  \int_{\Omega}S\left(\sigma_{\varepsilon}\mathbb{D}_{x'}\left[\sigma_{\varepsilon}^{-{r\over r-1}}\widetilde {u}_{\varepsilon }\right] + \sigma_\varepsilon h_\varepsilon^{-1}\partial_{y_3}\left[\sigma_{\varepsilon}^{-{r\over r-1}}\widetilde {u}_{\varepsilon }\right] \right): \left(\sigma_{\varepsilon}\mathbb{D}_{x^\prime}\left[\varphi\right] + \sigma_\varepsilon h_\varepsilon^{-1}\partial_{y_3}\left[\varphi \right]\right)\,dx'dy_3\\
\medskip
\displaystyle
+\mu\, \int_{\omega\times{1\over \delta_\varepsilon}Y'\times (0,1)}S\left(\mathbb{D}_{z'}\left[\sigma_{\varepsilon}^{-{r\over r-1}}\widehat {u}_{\varepsilon }\right] + \delta_\varepsilon^{2\over r}\sigma_\varepsilon  h_\varepsilon^{-1}\partial_{y_3}\left[\sigma_{\varepsilon}^{-{r\over r-1}}\widehat  {u}_{\varepsilon }\right] \right):\left(\mathbb{D}_{z^\prime}\left[\varphi\right] + \delta_\varepsilon^{2\over r}\sigma_\varepsilon  h_\varepsilon^{-1}\partial_{y_3}[\varphi]\right)\,dx'dz'dy_3\\
\medskip
\displaystyle
- \int_{\Omega} \widetilde {p}_{\varepsilon }\, {\rm div}_{x^\prime} \varphi^\prime\,dx^{\prime}dy_3- \delta_\varepsilon^{2}(\varepsilon\delta_\varepsilon)^{-1}\int_{\omega\times \mathbb{R}^2\times (0,1)} \widehat {p}_{\varepsilon}\,{\rm div}_{z^\prime} \varphi^\prime \,dx^{\prime}dz'dy_3\\
\medskip
\displaystyle
= \int_{\widetilde\Omega_{\varepsilon}}f'\cdot \varphi' \,dx^\prime dy_3+ O_\varepsilon.
\end{array}
\end{equation}

Now, we consider $v \in\mathcal{D}(\Omega)^3$ with $v_3\equiv 0$  and  satisfying the divergence condition ${\rm div}_{x'}\int_0^1v'\,dy_3=0$ in $\omega$. Moreover, we consider $V\in  \mathcal{D}(\Omega;{\bf K}_{Y'_s}^3)$ such that $D_{z'}V$ has a compact support and ${\rm div}_{z'}V'=0$ in $\omega\times \mathbb{R}^2\times (0,1)$. We set $V_{\varepsilon\delta_\varepsilon}$ by (\ref{vtestepdelta}), which has a limit at infinity denoted by  $V_\infty\in\mathcal{D}(\Omega)^3$, and assume $V_\infty(x',y_3)=v(x',y_3)$ a.e. in $\Omega$. Then, we consider the following   test function $\varphi$  in (\ref{formvarcvuprime}):
\begin{itemize}
\item[--] $\varphi_{\varepsilon\delta_\varepsilon}=V_{\varepsilon\delta_\varepsilon}-\sigma_{\varepsilon}^{-{r\over r-1}}\widetilde  u_{\varepsilon}=V\left(x^{\prime},{1\over \delta_\varepsilon}{x^{\prime}-\varepsilon\kappa({x'\over \varepsilon})\over \varepsilon},y_3\right)-\sigma_{\varepsilon}^{-{r\over r-1}}\widetilde  u_{\varepsilon}$ in the integrals in  $\Omega$ and $\widetilde\Omega_{\varepsilon}$,
\item[--] $\varphi_{\varepsilon\delta_\varepsilon}=V(x',z',y_3)-\sigma_{\varepsilon}^{-{r\over r-1}}\widehat  u_{\varepsilon }$  in the integrals in $\omega\times {1\over \delta_\varepsilon}Y'\times (0,1)$ and $\omega\times \mathbb{R}^2\times (0,1)$.
\end{itemize}
Applying the H${\rm \ddot{o}}$lder inequality and using the first estimate in (\ref{unfolding2}) and estimate (\ref{unfolding4}),  the fourth term in (\ref{formvarcvuprime}) satisfies
$$\begin{array}{l}
\medskip
\displaystyle
\left|\delta_\varepsilon^{2}(\varepsilon\delta_\varepsilon)^{-1}\int_{\omega\times \mathbb{R}^2\times (0,1)} \widehat {p}_{\varepsilon }\,{\rm div}_{z^\prime} \varphi^\prime_{\varepsilon \delta_\varepsilon} \,dx^{\prime}dz'dy_3\right|=\left|\delta_\varepsilon^{2}(\varepsilon\delta_\varepsilon)^{-1}\sigma_{\varepsilon}^{-{r\over r-1}}\int_{\omega\times \mathbb{R}^2\times (0,1)} \widehat {p}_{\varepsilon }\,{\rm div}_{z^\prime}\widehat  u_{\varepsilon }^\prime \,dx^{\prime}dz'dy_3\right|\\
\medskip
\displaystyle
\leq \delta_\varepsilon\varepsilon^{-1}\delta_\varepsilon^{-{2\over r'}} \sigma_{\varepsilon}^{-{r\over r-1}}=\delta_\varepsilon^{2-r\over r}\varepsilon^{-1}\sigma_\varepsilon^{-{r\over r-1}}=\sigma_{\varepsilon}^{-1}\sigma_\varepsilon^{-{r\over r-1}}=\sigma_\varepsilon^{1\over r-1}\to 0,
\end{array}$$
and then (\ref{formvarcvuprime})  reads as follows
\begin{equation*}
\begin{array}{l}
\medskip
\displaystyle
\mu\,  \int_{\Omega}S\left(\sigma_{\varepsilon}\mathbb{D}_{x'}\left[\sigma_{\varepsilon}^{-{r\over r-1}}\widetilde {u}_{\varepsilon }\right] + \sigma_\varepsilon h_\varepsilon^{-1}\partial_{y_3}\left[\sigma_{\varepsilon}^{-{r\over r-1}}\widetilde {u}_{\varepsilon }\right] \right): \left(\sigma_{\varepsilon}\mathbb{D}_{x^\prime}\left[\varphi_{\varepsilon\delta_\varepsilon}\right] + \sigma_\varepsilon h_\varepsilon^{-1}\partial_{y_3}\left[\varphi_{\varepsilon\delta_\varepsilon} \right]\right)\,dx'dy_3\\
\medskip
\displaystyle
+\mu\, \int_{\omega\times{1\over \delta_\varepsilon}Y'\times (0,1)}S\left(\mathbb{D}_{z'}\left[\sigma_{\varepsilon}^{-{r\over r-1}}\widehat {u}_{\varepsilon }\right] + \delta_\varepsilon^{2\over r}\sigma_\varepsilon  h_\varepsilon^{-1}\partial_{y_3}\left[\sigma_{\varepsilon}^{-{r\over r-1}}\widehat  {u}_{\varepsilon }\right] \right):\left(\mathbb{D}_{z^\prime}\left[\varphi_{\varepsilon\delta_\varepsilon}\right] + \delta_\varepsilon^{2\over r}\sigma_\varepsilon  h_\varepsilon^{-1}\partial_{y_3}[\varphi_{\varepsilon\delta_\varepsilon}]\right)\,dx'dz'dy_3\\
\medskip
\displaystyle
- \int_{\Omega} \widetilde {p}_{\varepsilon }\, {\rm div}_{x^\prime} \varphi_{\varepsilon\delta_\varepsilon}^\prime\,dx^{\prime}dy_3
= \int_{\widetilde\Omega_{\varepsilon}}f'\cdot \varphi_{\varepsilon\delta_\varepsilon}' \,dx^\prime dy_3+ O_\varepsilon.
\end{array}
\end{equation*}
From this, we deduce 
\begin{equation*}
\begin{array}{l}
\medskip
\displaystyle
-\mu\,  \int_{\Omega}\left(S\left(\sigma_{\varepsilon}\mathbb{D}_{x'}\left[\sigma_{\varepsilon}^{-{r\over r-1}}\widetilde {u}_{\varepsilon }\right] + \sigma_\varepsilon h_\varepsilon^{-1}\partial_{y_3}\left[\sigma_{\varepsilon}^{-{r\over r-1}}\widetilde {u}_{\varepsilon }\right] \right)-S(\sigma_\varepsilon\mathbb{D}_{x'}[V]+\sigma_\varepsilon h_\varepsilon^{-1}\partial_{y_3}[V])\right)\\
\medskip
\displaystyle 
\qquad : \left(\sigma_{\varepsilon}\mathbb{D}_{x^\prime}\left[\varphi_{\varepsilon\delta_\varepsilon}\right] + \sigma_\varepsilon h_\varepsilon^{-1}\partial_{y_3}\left[\varphi_{\varepsilon\delta_\varepsilon} \right]\right)\,dx'dy_3\\
\medskip
\displaystyle
-\mu \int_{\Omega}S(\sigma_\varepsilon\mathbb{D}_{x'}[V]+\sigma_\varepsilon h_\varepsilon^{-1}\partial_{y_3}[V]): \left(\sigma_{\varepsilon}\mathbb{D}_{x^\prime}\left[\varphi_{\varepsilon\delta_\varepsilon}\right] + \sigma_\varepsilon h_\varepsilon^{-1}\partial_{y_3}\left[\varphi_{\varepsilon\delta_\varepsilon} \right]\right)\,dx'dy_3
\\
\medskip
\displaystyle
-\mu\, \int_{\omega\times{1\over \delta_\varepsilon}Y'\times (0,1)}\left(S\left(\mathbb{D}_{z'}\left[\sigma_{\varepsilon}^{-{r\over r-1}}\widehat {u}_{\varepsilon }\right] + \delta_\varepsilon^{2\over r}\sigma_\varepsilon  h_\varepsilon^{-1}\partial_{y_3}\left[\sigma_{\varepsilon}^{-{r\over r-1}}\widehat  {u}_{\varepsilon }\right] \right)- S\left(\mathbb{D}_{z'}\left[V\right] + \delta_\varepsilon^{2\over r}\sigma_\varepsilon  h_\varepsilon^{-1}\partial_{y_3}\left[V\right] \right)\right)\\
\medskip
\displaystyle\qquad\qquad\qquad\qquad:\left(\mathbb{D}_{z^\prime}\left[\varphi_{\varepsilon\delta_\varepsilon}\right] + \delta_\varepsilon^{2\over r}\sigma_\varepsilon  h_\varepsilon^{-1}\partial_{y_3}[\varphi_{\varepsilon \delta_\varepsilon}]\right)\,dx'dz'dy_3\\
\medskip
\displaystyle
-\mu\, \int_{\omega\times{1\over \delta_\varepsilon}Y'\times (0,1)}S\left(\mathbb{D}_{z'}\left[V\right] + \delta_\varepsilon^{2\over r}\sigma_\varepsilon  h_\varepsilon^{-1}\partial_{y_3}\left[V\right] \right):\left(\mathbb{D}_{z^\prime}\left[\varphi_{\varepsilon\delta_\varepsilon}\right] + \delta_\varepsilon^{2\over r}\sigma_\varepsilon  h_\varepsilon^{-1}\partial_{y_3}[\varphi_{\varepsilon\delta_\varepsilon}]\right)\,dx'dz'dy_3\\
\medskip
\displaystyle
+ \int_{\Omega} \widetilde {p}_{\varepsilon }\, {\rm div}_{x^\prime} \varphi_{\varepsilon\delta_\varepsilon}^\prime\,dx^{\prime}dy_3= -\int_{\widetilde\Omega_{\varepsilon}}f'\cdot \varphi_{\varepsilon\delta_\varepsilon}'\,dx^\prime dy_3+ O_\varepsilon.
\end{array}
\end{equation*}
Because the operator $S$ is monotone  (i.e. $(S(\xi)-S(\zeta)):(\xi-\zeta)\geq 0$  for any $\xi,\zeta\in \mathbb{R}^3$), we have
\begin{equation}\label{form_var_tilde_proof3}
\begin{array}{l}\medskip\displaystyle 
\mu \int_{\Omega}S(\sigma_\varepsilon\mathbb{D}_{x'}[V]+\sigma_\varepsilon h_\varepsilon^{-1}\partial_{y_3}[V]): \left(\sigma_{\varepsilon}\mathbb{D}_{x^\prime}\left[\varphi_{\varepsilon\delta_\varepsilon}\right] + \sigma_\varepsilon h_\varepsilon^{-1}\partial_{y_3}\left[\varphi_{\varepsilon\delta_\varepsilon} \right]\right)\,dx'dy_3\\
\medskip
\displaystyle
+\mu\, \int_{\omega\times{1\over \delta_\varepsilon}Y'\times (0,1)}S\left(\mathbb{D}_{z'}\left[V\right] + \delta_\varepsilon^{2\over r}\sigma_\varepsilon  h_\varepsilon^{-1}\partial_{y_3}\left[V\right] \right):\left(\mathbb{D}_{z^\prime}\left[\varphi_{\varepsilon\delta_\varepsilon}\right] + \delta_\varepsilon^{2\over r}\sigma_\varepsilon  h_\varepsilon^{-1}\partial_{y_3}[\varphi_{\varepsilon\delta_\varepsilon}]\right)\,dx'dz'dy_3\\
\medskip
\displaystyle
- \int_{\Omega} \widetilde {p}_{\varepsilon }\, {\rm div}_{x^\prime} \varphi'_{\varepsilon\delta_\varepsilon}\,dx^{\prime}dy_3\geq \int_{\widetilde\Omega_{\varepsilon}}f'\cdot \varphi'_{\varepsilon\delta_\varepsilon}\,dx^\prime dy_3+ O_\varepsilon.
\end{array}
\end{equation}
{\bf Step 2}. In this step deduce the homogenized problem (\ref{thm:system})  in the case $\sigma_\varepsilon\approx h_\varepsilon$, with $\sigma_\varepsilon / h_\varepsilon\to \lambda \in (0,+\infty)$.

For this, we  pass to the limit in (\ref{form_var_tilde_proof3}) by using that  $\sigma_\varepsilon $ and $\delta^{2\over r}_\varepsilon$ tend to zero, $\sigma_\varepsilon / h_\varepsilon\to \lambda$ and  convergences (\ref{conv_u_tilde}), (\ref{conv_p_tilde}), (\ref{conv_Du_hat}) and (\ref{vepdeltest}). Thus,  we obtain
$$
\begin{array}{l}
\medskip
\displaystyle
\mu \int_{\Omega}S(\lambda \partial_{y_3}[V_\infty']): \left(\lambda\partial_{y_3}\left[V_\infty'-u' \right]\right)\,dx'dy_3\\
\medskip
\displaystyle
+\mu\, \int_{\omega\times (\mathbb{R}^2\setminus Y'_s)\times (0,1)}S\left(\mathbb{D}_{z'}\left[V\right] \right):
\mathbb{D}_{z^\prime}\left[V- U\right] \,dx'dz'dy_3- \int_{\Omega} p\, {\rm div}_{x^\prime} (V'_\infty-  u')\,dx^{\prime}dy_3
\geq \int_{\Omega}f'\cdot (V_\infty'-  u')\,dx^\prime dy_3.
\end{array}
$$
Because $p$ is independent of $y_3$, taking into account that $V'_\infty=v'$ in $\Omega$ and ${\rm div}_{x'}\int_0^1v'\,dy_3={\rm div}_{x'}\int_0^1u'\,dy_3=0$ in $\omega$,  we have 
$$\int_{\Omega} p\, {\rm div}_{x^\prime} (V'_\infty-  u')\,dx^{\prime}dy_3=\int_{\Omega} p\, {\rm div}_{x^\prime} (v'-  u')\,dx^{\prime}dy_3=\int_{\omega} p\, {\rm div}_{x^\prime} \left(\int_0^1(v'-  u')dy_3\right)\,dx^{\prime}=0.$$
Then, we have
$$
\begin{array}{l}
\medskip
\displaystyle
\mu \int_{\Omega}S(\lambda \partial_{y_3}[V_\infty']): \left(\lambda\partial_{y_3}\left[V_\infty'-u' \right]\right)\,dx'dy_3\\
\medskip
\displaystyle
+\mu\, \int_{\omega\times (\mathbb{R}^2\setminus Y'_s)\times (0,1)}S\left(\mathbb{D}_{z'}\left[V\right] \right):
\mathbb{D}_{z^\prime}\left[V- U\right] \,dx'dz'dy_3\geq \int_{\Omega}f'\cdot (V_\infty'-  u')\,dx^\prime dy_3.
\end{array}
$$
From Minty's Lemma \cite[Chapter 3, Lemma 1.2]{Lions2} and the homogeneity property of the operator $S$, then previous inequality is equivalent  to the following variational formulation
\begin{equation}\label{casifinalformvar}
\begin{array}{l}
\medskip
\displaystyle
\mu\lambda^{r}  \int_{\Omega}S\left(\partial_{y_3}\left[ u'\right]\right):  \partial_{y_3}\left[V'_\infty\right] \,dx'dy_3
+\mu\, \int_{\omega\times   (\mathbb{R}^2\setminus Y'_s)\times (0,1)}\!\!\!\!S\left(\mathbb{D}_{z'}\left[ U\right] \right):
\mathbb{D}_{z^\prime}\left[V\right] \,dx'dz'dy_3
= \int_{ \Omega}f'\cdot V_\infty'\,dx^\prime dy_3.
\end{array}
\end{equation}
Following \cite[Proof of Theorem 2.1]{Anguiano_SG_Lower}, it is easy to prove that $U_3\equiv 0$. Then, we get that   (\ref{casifinalformvar}) reads as follows
\begin{equation}\label{casifinalformvar21}
\begin{array}{l}
\medskip
\displaystyle
\mu\lambda^r  \int_{\Omega}S\left(\partial_{y_3}\left[ u'\right]\right):  \partial_{y_3}\left[v'\right] \,dx'dy_3
+\mu\, \int_{\omega\times   (\mathbb{R}^2\setminus Y'_s)\times (0,1)}S\left(\mathbb{D}_{z'}\left[ U'\right] \right):
\mathbb{D}_{z^\prime}\left[V'\right] \,dx'dz'dy_3
= \int_{ \Omega}f'\cdot v'\,dx^\prime dy_3,
\end{array}
\end{equation}
which, by density,  it holds for every $(v', V')\in\mathbb{W}$.\\

To finish, it remains to deduce problem (\ref{thm:system}) for $u'$ identifying $U'$ in (\ref{casifinalformvar21}) by means of the drag force $\mathcal{G}$ given by (\ref{DragForceG}) and the auxiliary problems (\ref{LocalProblems}).  This is developed in Step 2 of the proof of Theorem 2.1 in \cite{Anguiano_SG_Lower}, so for the reader's convenience, we reproduce it below. For this purpose, we consider the auxiliary problems $(w^{\xi'},\pi^{\xi'})$ $\forall\, \xi'\in\mathbb{R}^2$,  defined by (\ref{LocalProblems}) and the drag force function $\mathcal{G}$ defined by (\ref{DragForceG}). Then, we take in (\ref{casifinalformvar21}) the pairs of functions $(u',U'), (v',V')\in \mathbb{W}$ in  the following form
\begin{equation}\label{identification}
U'(x',z',y_3)=w^{u'(x',y_3)}(z'),\quad V'(x',z',y_3)=w^{v'(x',y_3)}(z'),
\end{equation}
a.e. in $\omega\times (\mathbb{R}^2\setminus Y'_s)\times (0,1)$, and then, we deduce
\begin{equation}\label{casifinalformvar212}
\begin{array}{l}
\medskip
\displaystyle
\mu\lambda^r  \int_{\Omega}S\left(\partial_{y_3}\left[ u'\right]\right):  \partial_{y_3}\left[v'\right] \,dx'dy_3
+\mu\, \int_{\Omega}\mathcal{G}(u')\cdot v'\,dx'dy_3
= \int_{ \Omega}f'\cdot v'\,dx^\prime dy_3,
\end{array}
\end{equation}
for every function $v'\in \mathcal{V}$, where
\begin{equation}\label{defV}\displaystyle \mathcal{V}=\left\{
\begin{array}{l}
v'\in W^{1,r}(0,1;L^r(\omega)^2)\,:\, {\rm div}_{x'}\left(\int_0^1 v'(x',y_3)\,dy_3\right)=0\quad\hbox{in } \omega,\quad \left(\int_0^1 v'(x',y_3)\,dy_3\right)\cdot n'=0\quad\hbox{in } \partial \omega
\end{array}
\right\}.
\end{equation}
From the properties of $\mathcal{G}$ described in Remark \ref{RemG} and following \cite[Lemma 4.4]{MarusicPaloka}, the variational formulation (\ref{casifinalformvar212}) has a unique solution $u'\in \mathcal{V}$. Finally,taking into account that $S(\partial_{y_3}[u'])=2^{-{r\over 2}}S(\partial_{y_3}u')$, we  deduce that there exists $q\in L^{r'}(\omega)\setminus \mathbb{R}$, where $q$ coincides with the limit pressure $p$, such that (\ref{casifinalformvar212}) is equivalent to system  (\ref{thm:system2}) for $u'$. Since problem (\ref{thm:system2}) has a unique solution (see proof of Lemma 4.4 in \cite{MarusicPaloka}), then the entire sequence $(\widetilde  u_\varepsilon, \widetilde  p_\varepsilon)$ converges to $(u, p)$. This finishes the proof in this case.
 \\

{\bf Step 3}. In this step deduce the homogenized problem (\ref{thm:system2})  in the case $\sigma_\varepsilon\ll h_\varepsilon$ (i.e. $\sigma_\varepsilon / h_\varepsilon\to \lambda=0$).

The proof is the same as Step 2 just taking into account that, in this case, $\lambda=0$. Thus, we deduce the homogenized variational equation 
\begin{equation}\label{casifinalformvar31}
\begin{array}{l}
\medskip
\displaystyle
\mu\, \int_{\omega\times   (\mathbb{R}^2\setminus Y'_s)\times (0,1)}S\left(\mathbb{D}_{z'}\left[ U'\right] \right):
\mathbb{D}_{z^\prime}\left[V'\right] \,dx'dz'dy_3
= \int_{ \Omega}f'\cdot v'\,dx^\prime dy_3,
\end{array}
\end{equation}
 for every $(v', V')\in\mathbb{W}$.    By means of the identification (\ref{identification}), we deduce that (\ref{casifinalformvar31}) is written as follows \\
 \begin{equation}\label{casifinalformvar312}
\begin{array}{l}
\medskip
\displaystyle
\mu\, \int_{\Omega}\mathcal{G}(u')\cdot v' \,dx'dy_3
= \int_{ \Omega}f'\cdot v'\,dx^\prime dy_3,
\end{array}
\end{equation}
for every function $v'\in \mathcal{V}$, where $\mathcal{V}$ is defined in (\ref{defV}). Finally, similarly to the critical case,  we can deduce that there exists $q\in L^{r'}(\omega)\setminus \mathbb{R}$, where $q$ coincides with the limit pressure $p$. Then,  we can deduce that $(u',p)$ satisfies (\ref{casifinalformvar312}), which is equivalent to system  
 \begin{equation}\label{thm:system2_proof}
\left\{\begin{array}{rl}
\medskip
\displaystyle     \mu \mathcal{G}(u') = f'(x')-\nabla_{x'}p(x'),\quad u_3\equiv 0&\hbox{in }\Omega,\\
\medskip
\displaystyle {\rm div}_{x'}\left(\int_0^1 u'\,dy_3\right)=0&\hbox{in }\omega,\\
\medskip
\displaystyle  \left(\int_0^1 u'\,dy_3\right)\cdot n'=0&\hbox{on }\partial\omega.
\end{array}\right.
\end{equation}
where $\mathcal{G}$ is the drag force function defined by (\ref{DragForceG}). From the homogeneity condition (\ref{homogeneitycon}), we have that $\mu \mathcal{G}(u')=\mathcal{G}(\mu^{1\over r-1}u')$. Moreover, the fact that $\mathcal{G}^{-1}(\xi')$ is well defined for $\xi'\in \mathbb{R}^2$ (see \cite[Remark 4.6]{MarusicPaloka2}), then problem (\ref{thm:system2_proof})$_1$ is rewritten by 
\begin{equation}\label{thm:system2_proof2}
\left\{\begin{array}{rl}
\medskip
\displaystyle     u'(x',y_3) ={1\over \mu^{1\over r-1}}\mathcal{G}^{-1}\left(f'(x')-\nabla_{x'}p(x')\right),\quad u_3\equiv 0&\hbox{in }\Omega,\\
\medskip
\displaystyle {\rm div}_{x'}\left(\int_0^1 u'\,dy_3\right)=0&\hbox{in }\omega,\\
\medskip
\displaystyle  \left(\int_0^1 u'\,dy_3\right)\cdot n'=0&\hbox{on }\partial\omega.
\end{array}\right.
\end{equation}

From the properties and convergences relating $\mathcal{G}^{-1}$ and $\mathcal{U}_{\delta_\varepsilon}$ given in Remark \ref{propG1}, given $\xi'=f'(x')-\nabla_{x'}p(x')$ for a.e. $x'\in \omega$, we have that 
$$\mathcal{\widetilde U}_{\delta_\varepsilon}(f'(x')-\nabla_{x'}p(x'))=\delta_\varepsilon^{2-r\over r-1}\mathcal{U}_{\delta_\varepsilon}(f'(x')-\nabla_{x'}p(x')),$$
where $\mathcal{\widetilde U}_{\delta_\varepsilon}$ is the permeability function given by 
$$\mathcal{\widetilde U}_{\delta_\varepsilon}(f'(x')-\nabla_{x'}p(x'))=\int_{\delta_\varepsilon^{-1}Y'\setminus Y'_s} \delta_\varepsilon^2 w^{f'(x')-\nabla_{x'}p(x')}_{\delta_\varepsilon}\,dy',$$
where $w^{f'(x')-\nabla_{x'}p(x')}_{\delta_\varepsilon}$ is the unique solution of problem (\ref{aux_problem2}) with $\xi'=f'(x')-\nabla_{x'}p(x')$. Moreover, due to the relation  given by (\ref{lowvolume}), we have that $$\lim_{\delta_\varepsilon\to 0}\delta_\varepsilon^{2-r\over r-1}\mathcal{U}_{\delta_\varepsilon}(f'(x')-\nabla_{x'}p(x'))=\mathcal{G}^{-1}(f'(x')-\nabla_{x'}p(x')),$$
then we have that $\mathcal{G}^{-1}$ does not depend on $y_3$, so $U_{av}(x')=\int_0^1u(x',y_3)\,dy_3$ satisfies
\begin{equation}\label{thm:system2_proof21}
\left\{\begin{array}{rl}
\medskip
\displaystyle     U_{av}'(x') ={1\over \mu^{1\over r-1}}\mathcal{G}^{-1}\left(f'(x')-\nabla_{x'}p(x')\right),\quad U_{av, 3}\equiv 0&\hbox{in }\omega,\\
\medskip
\displaystyle {\rm div}_{x'}U_{av}'(x')=0&\hbox{in }\omega,\\
\medskip
\displaystyle U_{av}'(x')\cdot n'=0&\hbox{on }\partial\omega,
\end{array}\right.
\end{equation}
which has a unique solution $p\in W^{1,r'}(\omega)/\mathbb{R}$. Then, the entire sequence $(\widetilde  u_\varepsilon, \widetilde  p_\varepsilon)$ converges to $(u, p)$, which finishes the proof in this case.\\

{\bf Step 4}. In this step we develop the case $\sigma_\varepsilon\gg h_\varepsilon$. The proof follows previous steps, so we will give some remarks. We take as a starting point the variational formulation (\ref{FV_tilde}).  Below, we analyze every term:  
 \begin{itemize}
 \item First term. Taking  into account convergences given in Lemmas \ref{LemmaConvtulde} in the case $\sigma_\varepsilon\gg  h_\varepsilon$, we rewrite this term as follows
 \begin{equation}\label{integral_velocity1}
\begin{array}{l}
\medskip
\displaystyle
\mu\,  \int_{\Omega}S\left(\mathbb{D}_{x'} \left[\widetilde {u}_{\varepsilon} \right] + h_{\varepsilon}^{-1}\partial_{y_3}[\widetilde  u_{\varepsilon}] \right): \left( \mathbb{D}_{x^\prime}\left[\varphi\right] + h_{\varepsilon}^{-1}\partial_{y_3}\left[\varphi\right]\right)\,dx'dy_3\\
\medskip
\displaystyle
=\mu\,  \int_{\Omega}S\left(h_{\varepsilon}\mathbb{D}_{x'}\left[h_{\varepsilon}^{-{r\over r-1}}\widetilde {u}_{\varepsilon}\right] + \partial_{y_3}\left[h_{\varepsilon}^{-{r\over r-1}}\widetilde {u}_{\varepsilon}\right] \right): \left(h_{\varepsilon}\mathbb{D}_{x^\prime}\left[\varphi\right] +  \partial_{y_3}\left[\varphi\right]\right)\,dx'dy_3.
\end{array}
\end{equation}
 
 \item Second term. 
From the changes of variables  (\ref{CV})  and using (\ref{derivadatest}), we deduce
\begin{eqnarray}
\displaystyle
&&\mu\,  (\varepsilon\delta_\varepsilon)^{-1}\int_{\Omega}S\left(\mathbb{D}_{x'} \left[\widetilde {u}_{\varepsilon} \right] + h_{\varepsilon}^{-1}\partial_{y_3}[\widetilde  u_{\varepsilon}] \right): \mathbb{D}_{z^\prime}\left[\varphi\right] \,dx'dy_3  \nonumber\\
\medskip
\displaystyle
&&=\mu\, (\varepsilon\delta_\varepsilon)^{-1}\delta_\varepsilon^2 \int_{\omega\times{1\over \delta_\varepsilon}Y'\times (0,1)}S\left( (\varepsilon\delta_\varepsilon)^{-1}\mathbb{D}_{z'}\left[\widehat {u}_{\varepsilon}\right] +h_{\varepsilon}^{-1}\partial_{y_3}\left[\widehat  {u}_{\varepsilon}\right] \right):(\mathbb{D}_{z^\prime}\left[\varphi\right] + \mathbb{D}_{z'}[\varTheta_{\varepsilon\delta_\varepsilon}])\,dx'dz'dy_3 \label{integral_velocity2_bis}
\\
\medskip
\displaystyle
&&=h_\varepsilon\sigma_\varepsilon^{-1}\mu\, \int_{\omega\times{1\over \delta_\varepsilon}Y'\times (0,1)}S\left(\mathbb{D}_{z'}\left[ \sigma_\varepsilon^{-1}h_{\varepsilon}^{-{1\over r-1}}\widehat {u}_{\varepsilon}\right] + \partial_{y_3}\left[\delta_\varepsilon^{2\over r}h_{\varepsilon}^{-{r\over r-1}}\widehat  {u}_{\varepsilon}\right] \right):
(\mathbb{D}_{z^\prime}\left[\varphi\right]+ \mathbb{D}_{z'}[\varTheta_{\varepsilon\delta_\varepsilon}])\,dx'dz'dy_3.\nonumber
\end{eqnarray}
From the H${\rm \ddot{o}}$lder inequality, estimates (\ref{unfolding1_3}) and using (\ref{convtestder}), we get 
$$\begin{array}{l}\medskip
\displaystyle
\left|
h_\varepsilon\sigma_\varepsilon^{-1}\mu\, \int_{\omega\times{1\over \delta_\varepsilon}Y'\times (0,1)}S\left(\mathbb{D}_{z'}\left[ \sigma_\varepsilon^{-1}h_{\varepsilon}^{-{1\over r-1}}\widehat {u}_{\varepsilon}\right] + \partial_{y_3}\left[\delta_\varepsilon^{2\over r}h_{\varepsilon}^{-{r\over r-1}}\widehat  {u}_{\varepsilon}\right] \right):
\mathbb{D}_{z^\prime}\left[\varphi\right]\,dx'dz'dy_3\right|\leq h_\varepsilon \sigma_\varepsilon^{-1}\to 0,
\\
\medskip
\displaystyle
\left|h_\varepsilon\sigma_\varepsilon^{-1}\mu\, \int_{\omega\times{1\over \delta_\varepsilon}Y'\times (0,1)}S\left(\mathbb{D}_{z'}\left[ \sigma_\varepsilon^{-1}h_{\varepsilon}^{-{1\over r-1}}\widehat {u}_{\varepsilon}\right] + \partial_{y_3}\left[\delta_\varepsilon^{2\over r}h_{\varepsilon}^{-{r\over r-1}}\widehat  {u}_{\varepsilon}\right] \right):
 \mathbb{D}_{z'}[\varTheta_{\varepsilon\delta_\varepsilon}]\,dx'dz'dy_3\right|\leq C \varepsilon h_\varepsilon\sigma_\varepsilon^{-1}\to 0,
 \end{array}$$
and then,  (\ref{integral_velocity2_bis}) satisfies
\begin{eqnarray}
&&\mu\,  (\varepsilon\delta_\varepsilon)^{-1}\int_{\Omega}S\left(\mathbb{D}_{x'} \left[\widetilde {u}_{\varepsilon} \right] + h_{\varepsilon}^{-1}\partial_{y_3}[\widetilde  u_{\varepsilon}] \right): \mathbb{D}_{z^\prime}\left[\varphi\right] \,dx'dy_3  \to 0.
\end{eqnarray}

\item Third and fourth terms. Similarly to the other cases, it holds
\begin{equation}\label{form_var_pressurebis_super}
\begin{array}{l}\medskip
\displaystyle
- \int_{\Omega} \widetilde {p}_{\varepsilon}\, {\rm div}_{x^\prime}\varphi^\prime\,dx^{\prime}dy_3- (\varepsilon\delta_\varepsilon)^{-1} \int_{\Omega} \widetilde {p}_{\varepsilon\delta_\varepsilon}\, {\rm div}_{z^\prime} \varphi^\prime\,dx^{\prime}dy_3\\
\medskip
\displaystyle
=- \int_{\Omega} \widetilde {p}_{\varepsilon}\, {\rm div}_{x^\prime} \varphi^\prime\,dx^{\prime}dy_3-\delta_\varepsilon^2 (\varepsilon\delta_\varepsilon)^{-1} \int_{\omega\times \mathbb{R}^2\times (0,1)} \widehat {p}_{\varepsilon}\, \left({\rm div}_{z^\prime} \varphi^\prime+ {\rm div}_{z^\prime} \varTheta ^\prime_{\varepsilon\delta_\varepsilon}\right)\,dx^{\prime}dz'dy_3.
\end{array}
\end{equation}
and also, by the H${\rm \ddot{o}}$lder inequality, estimate (\ref{unfolding4}) and   (\ref{convtestder}), then
$$\left|\delta_\varepsilon^2  (\varepsilon\delta_\varepsilon)^{-1}\int_{\omega\times \mathbb{R}^2\times (0,1)} \widehat {p}_{\varepsilon}\, {\rm div}_{z^\prime} \varTheta ^\prime_{\varepsilon\delta_\varepsilon}\,dx^{\prime}\,dz'dy_3\right|\leq C\delta_\varepsilon^{2-r\over r}\to 0.
$$
 \end{itemize}
Therefore, from this analysis, we deduce that (\ref{FV_tilde}) is rewritten as follows
$$\begin{array}{l}
\medskip\displaystyle\mu\,  \int_{\Omega}S\left(h_{\varepsilon}\mathbb{D}_{x'}\left[h_{\varepsilon}^{-{r\over r-1}}\widetilde {u}_{\varepsilon}\right] + \partial_{y_3}\left[h_{\varepsilon}^{-{r\over r-1}}\widetilde {u}_{\varepsilon}\right] \right): \left(h_{\varepsilon}\mathbb{D}_{x^\prime}\left[\varphi\right] +  \partial_{y_3}\left[\varphi\right]\right)\,dx'dy_3- \int_{\Omega} \widetilde {p}_{\varepsilon}\, {\rm div}_{x^\prime} \varphi^\prime\,dx^{\prime}dy_3\\
\medskip
\displaystyle
=  \int_{\widetilde \Omega_{\varepsilon}}f'\cdot \varphi'\,dx^\prime dy_3+O_\varepsilon.
\end{array}
$$
Observe that there is no effects from the microstucture of the original domain, because the terms involving the unfolded functions vanish.
Therefore, passing to the limit when $\varepsilon$ tends to zero by monotonicity arguments as in the previous cases, we can deduce the homogenized variational formulation
$$\begin{array}{l}
\medskip\displaystyle\mu\,  \int_{\Omega}S\left(\partial_{y_3}\left[u'\right] \right):  \partial_{y_3}\left[v'\right] \,dx'dy_3
=  \int_{  \Omega}f'\cdot v'\,dx^\prime dy_3,
\end{array}
$$
for every $v'$ belonging to the space
$$\mathbb{\widetilde W}=\left\{\begin{array}{ll}
\medskip 
v'\in W^{1,r}(0,1;L^r(\omega)^2)\quad :& \quad  v'=0\ \hbox{on}\ y_3=\{0,1\}\\
\medskip
\displaystyle {\rm div}_{x'}\left(\int_0^1v'\,dy_3\right)=0\hbox{ in }\omega,&\displaystyle \left(\int_0^1v'\,dy_3\right)\cdot n=0\hbox{ on }\partial\omega
\end{array}\right\}.$$
Taking into account that $S(\partial_{y_3}[u'])=2^{-{r\over 2}}S(\partial_{y_3}u')$, we can deduce that there exists $q\in L^{r'}(\omega)\setminus \mathbb{R}$, where $q$ coincides with the limit pressure $p$, such that the previous variational formulation is equivalent to the system 
\begin{equation}\label{thm:system3}
\left\{\begin{array}{rl}
\medskip
\displaystyle - \mu 2^{-{r\over 2}}\,\partial_{y_3}\left(|\partial_{y_3}  u'|^{r-2}\partial_{y_3}  u'\right)+\nabla_{x'}p(x')=f'(x')&\hbox{in }\Omega,\\
\medskip
  u'=0&\hbox{on }y_3=\{0,1\},\\
\medskip
\displaystyle {\rm div}_{x'}\left(\int_0^1  u'\,dy_3\right)=0&\hbox{in }\omega,\\
\medskip
\displaystyle \left(\int_0^1  u'\,dy_3\right)\cdot n'=0&\hbox{on }\partial\omega.
\end{array}\right.
\end{equation}
It holds that $(u',p)\in W^{1,r}(0,1;L^r(\omega)^2)\times \left((L^{r'}(\omega)\setminus\mathbb{R})\cap W^{1,r'}(\omega)\right)$ is the unique solution of (\ref{thm:system3}), and so, the entire sequence $(\widetilde  u_\varepsilon, \widetilde  p_\varepsilon)$ converges to $(u, p)$ with $u_3\equiv 0$. Moreover, it can be derived  the nonlinear Reynolds problem (\ref{U_Reynolds_super})-(\ref{Reynolds_super}).  We omit the derivation of the Reynolds problem and the regularity for $p$, which can be  seen in \cite[Propositions 3.3 and 3.4]{MikelicTapiero}.

\end{proof}

\end{document}